\documentclass[11pt]{article}
\usepackage{amsmath}
\usepackage{amsmath,amssymb,amsthm}
\usepackage{graphicx}
\usepackage{enumerate}
\usepackage{natbib}
\usepackage{url}
\usepackage{mathtools}

\usepackage{algorithm}
\usepackage{algpseudocode}

\newtheorem{theorem}{Theorem}
\newtheorem{lemma}{Lemma}
\newtheorem{corollary}{Corollary}

\newtheorem{proposition}{Proposition}

\newcommand{\Pimat}{\boldsymbol{\Pi}}

\newcommand{\betavec}{\boldsymbol{\beta}}

\newcommand{\zerovec}{\boldsymbol{0}}
\newcommand{\onevec}{\boldsymbol{1}}

\newcommand{\evec}{\boldsymbol{e}}

\newcommand{\uvec}{\boldsymbol{u}}

\newcommand{\vvec}{\boldsymbol{v}}

\newcommand{\yvec}{\boldsymbol{y}}
\newcommand{\zvec}{\boldsymbol{z}}

\newcommand{\A}{\textbf{A}}
\newcommand{\C}{\textbf{C}}

\newcommand{\F}{\textbf{F}}
\newcommand{\I}{\textbf{I}}
\newcommand{\J}{\textbf{J}}
\newcommand{\M}{\textbf{M}}

\newcommand{\Lmat}{\textbf{L}}
\newcommand{\Smat}{\textbf{S}}
\newcommand{\Pmat}{\textbf{P}}
\newcommand{\Fmat}{\textbf{F}}
\newcommand{\V}{\textbf{V}}

\newcommand{\X}{\textbf{X}}
\newcommand{\Z}{\textbf{Z}}

\newcommand{\sign}{\text{sign}}

\DeclareMathOperator*{\argmin}{arg\,min}

\newcommand{\blind}{1}

\begin{document}

\def\spacingset#1{\renewcommand{\baselinestretch}%
{#1}\small\normalsize} \spacingset{1}


\if1\blind
{
  \title{\bf An Optimal Design Framework for Lasso Sign Recovery}
  \author{Jonathan W. Stallrich\\
    Department of Statistics, North Carolina State University\\
    Kade Young \\
    Department of Statistics, North Carolina State University
    \\
    Maria L. Weese\\
    Department of Information Systems \& Analytics, Miami University\\
    Byran J. Smucker\\
    Department of Statistics, Miami University, 
    \\
    and \\
    David J. Edwards\\
    Department of Statistical Sciences and Operations Research, VCU}
  \maketitle
} \fi

\if0\blind
{
  \bigskip
  \bigskip
  \bigskip
  \begin{center}
    {\LARGE\bf Optimal Supersaturated Designs for Lasso Sign Recovery}
\end{center}
  \medskip
} \fi

\bigskip
\begin{abstract}
Supersaturated designs investigate more factors than there are runs, and are often constructed under a criterion measuring a design's proximity to an unattainable orthogonal design. The most popular analysis identifies active factors by inspecting the solution path of a penalized estimator, such as the lasso. Recent criteria encouraging positive correlations between factors have been shown to produce designs with more definitive solution paths so long as the active factors have positive effects. Two open problems affecting the understanding and practicality of supersaturated designs are: (1) do optimal designs under existing criteria maximize support recovery probability across an estimator's solution path, and (2) why do designs with positively correlated columns produce more definitive solution paths when the active factors have positive sign effects? To answer these questions, we develop criteria maximizing the lasso's sign recovery probability. We prove that an orthogonal design is an ideal structure when the signs of the active factors are unknown, and a design constant small, positive correlations is ideal when the signs are assumed known. A computationally-efficient design search algorithm is proposed that first filters through optimal designs under new heuristic criteria to select the one that maximizes the lasso sign recovery probability. 
\end{abstract}

\noindent%
{\it Keywords:}  Constrained-positive $\text{Var}(s)$-criterion; Gauss Dantzig selector; Supersaturated designs; $\text{UE}(s^2)$-criterion; Variable selection
\vfill

\newpage






%


\spacingset{1.9} 

\section{Introduction}\label{sec:introduction}

A screening experiment aims to learn which $k$ out of $p$ factors most influence the response variable using as few runs, $n$, as possible. Achieving this goal with a traditional least-squares analysis is challenging because often not all possible factorial effects can estimated simultaneously. 
Still, effective screening is feasible under the effect sparsity principle and when the experimental design is carefully selected \citep{box1961,box19612,mee,xu2009recent, mee2017selecting}.  Two-level supersaturated designs, or SSDs, push screening to its limits with $n < p+1$ so that not even the $p$ main effects can be estimated simultaneously via least-squares. The recent literature on SSD analysis instead promotes penalized regression methods 
to perform screening under the main effects model, but the literature on the construction of SSDs has not kept up with this development. The present work addresses the gap between analysis and construction by developing an optimal design framework that tailors design selection to statistical properties of the lasso estimator.

Throughout this paper, we will assume a  main effects model $\yvec = \beta_0\onevec+\X\betavec + \evec$
where $\X$ is the $n \times p$ design matrix with elements $x_{ij}=\pm 1$, $\betavec=(\beta_1,\dots,\beta_p)^T$ is a sparse vector with $k < p$ nonzero elements, and $\evec \sim N(\zerovec,\sigma^2\I)$.  Without loss of generality, assume $\sigma^2=1$, making $\betavec$ the signed signal-to-noise ratios. Then the analysis goal is estimation of the support of $\betavec$, denoted $\mathcal{A}=\{j : |\beta_j| > 0\}$, although for technical reasons we will be more interested in estimating the sign vector of $\betavec$, denoted $\boldsymbol{z}=\sign(\betavec)$. We will focus our attention on optimal SSDs, although the framework developed in this paper could be applied to more general screening designs. Note that when we indicate designs are ``optimal'', we are referring to the best designs found by heuristic optimality procedures, unless otherwise indicated.

\subsection{Background and Open Problems}
The least-squares estimator for $\betavec$ is not unique for any SSD, complicating its analysis and design construction. If there were SSDs with unique least-squares estimators, the ideal $\X$ would satisfy $\Smat=\Lmat^T\Lmat=n\I_{p+1}$ where $ \Lmat = (\onevec | \X)$, as its $\hat{\beta}_j$ would have the minimum possible variance across all SSDs. 
These so-called orthogonal designs only exist when $n=0 \, (\text{mod} \ 4) > p$. For arbitrary $n>p$, a design could instead be selected by minimizing a variance-based criterion such as the $D$-criterion, $|(\Lmat^T\Lmat)^{-1}|$, or $A$-criterion, $\text{tr}[(\Lmat^T\Lmat)^{-1}]$. This optimal design framework based on least-squares is well-developed  and tractable \citep{pukelsheim2006optimal, goos2011optimal}, but is not applicable for SSDs. 

Most SSD construction methods focus on optimizing heuristic criteria that find a design whose $\Smat$ is close to $n\I_{p+1}$. \cite{Booth62} proposed the criterion $\text{E}(s^2)=2/(p(p-1))\sum_{1<i<j\leq p} s_{ij}^2$ assuming balanced designs (i.e., $\X^T\onevec=\zerovec$), which has been extensively studied \citep[see][]{Lin93, Wu93,Nguyen96,LiWu97,tang1997method,Ryan2007}. 
\cite{jones2014optimal} proposed the unconditional $\text{E}(s^2)$-criterion, $\text{UE}(s^2)=2/(p(p+1)) \sum_{0\leq i<j\leq p}s_{ij}^2$, that includes the $s_{0j}^2$ elements. Other criteria have been proposed that measure a design's quality in terms of subsets of columns \citep{deng1996measurement, deng1999resolution, Jones2009}, which is connected to a stepwise or all-subsets analysis \citep{Lin95, Abraham, westfall1998forward, liu2007construction}. \cite{CHEN199899} and \cite{SARKAR20091224} investigated support recovery properties of SSDs  under such sequential analyses. 

\cite{Li2002} and \cite{phoa2009analysis} advocated for SSDs to be analyzed under penalized estimation that induces sparse estimates. The previously mentioned SSD criteria, however, are not directly tied to statistical properties of such an analysis. There have been some attempts in the literature to address this disconnect.  \cite{MarleyWoods10}, \cite{draguljic_etal2014}, and \cite{weese_etal2015} have compared optimal SSDs under different criteria via simulation with the lasso and/or Gauss-Dantzig selector (GDS), finding little difference between the SSDs in terms of support recovery. Additionally, the Dantzig selector and its statistical properties are closely related to the lasso estimator \citep{meinshausen2007discussion, lounici2008sup, bickel2009simultaneous, james2009dasso, asif2010lasso}, but the mathematics of the lasso are more tractable, making its optimal design framework more straightforward. 

Recently, \cite{weese2017powerful} and \cite{weese2021strategies} proposed $\text{Var}(s+)$-optimal designs  that minimize
\[
   \text{Var}(s)= \text{UE}(s^2)-\text{UE}(s)^2  \text{ such that } \frac{\text{UE}^*(s^2)}{\text{UE}(s^2)}>c  \text{ and } \text{UE}(s)>0\ ,\ 
\]
where $\text{UE}^*(s^2)$ is the optimal $\text{UE}(s^2)$ value for the given $n$ and $p$, and $\text{UE}(s)=2/(p(p+1)) \sum_{0\leq i < j \leq p} s_{ij}$. The ideal structure of $\Smat$ under $\text{Var}(s+)$ is notably non-orthogonal, having average off-diagonal elements equal to a small, positive constant, making $\text{Var}(s)=0$ while satisfying the constraints on $\text{UE}(s^2)$ and $\text{UE}(s)$. Simulation studies showed $\text{Var}(s+)$-optimal designs (with $c=0.8$) performed similarly to other optimal SSDs under the GDS, but when the effect directions were known in advance (i.e. the nonzero signs of $\betavec$ were all positive), the $\text{Var}(s+)$-optimal design had better support recovery properties.  The simulation studies mentioned above only perform relative comparisons of SSDs generated by heuristic criteria and so do not answer the following open problems:
\begin{enumerate}
    \item How well do existing heuristic criteria find SSDs that maximize the probability of support recovery across all possible SSDs?
    \item Why does the $\text{Var}(s+)$-criterion identify SSDs with better support recovery probability when the active effects have positive signs?
\end{enumerate} 
The theory developed in this paper answers both of these problems.

\cite{singh2023selection} addressed problem 1 by proposing two criteria related to the restricted isometry property \citep{candes_tao2007} and irrepresentable condition \citep{gai2013model} of the GDS, which are sufficient for estimation and sign recovery consistency. For given $\mathcal{A}$ of size $k=|\mathcal{A}|$, let $\C_\mathcal{A}$ denote the correlation matrix of the columns in $\mathcal{A}$. Let $0 \leq \lambda_{\mathcal{A}1} \leq 1$ and $\lambda_{\mathcal{A}k} \geq 1$ denote the minimum and maximum eigenvalues of $\C_\mathcal{A}$, respectively, and define $\tau_\mathcal{A}=\max (1-\lambda_{\mathcal{A}1}, \lambda_{\mathcal{A}k}-1)$. \cite{singh2023selection} propose calculating the 95th percentile of all $\tau_\mathcal{A}$, denoted by $\zeta_k$. They then calculate $\zeta_k$ for $k=2,\dots,\lceil n/2 \rceil$ and define the estimation capacity criterion:
\[
\text{EC} = \text{mean}(\zeta_2,\dots,\zeta_{\lceil n/2 \rceil})\ .\    \]
Ideally, EC would be 0,  corresponding to an orthogonal design. Their second criterion assumes the sign of the active effects is always positive, denoted by $\onevec_k$. For a given $\mathcal{A}$, the following two weak sign consistency conditions are checked
\begin{align}
|\C_{\mathcal{I}\mathcal{A}}\C_{\mathcal{A}}^{-1}\onevec_k| &< \onevec_{p-k}\ ,\ \label{eq:WSC1}\\ 
|\C_{\mathcal{I}\mathcal{A}}\C_{\mathcal{A}}^{-1}\text{sign}(\C_{\mathcal{A}}^{-1}\onevec_k)| &< \onevec_{p-k}\ ,\ \label{eq:WSC2}
\end{align}
where $\C_{\mathcal{I}\mathcal{A}}$ is the matrix of cross correlations of columns between the inactive set, $\mathcal{I}=\{j : \beta_j = 0\}$, and $\mathcal{A}$. Similar to EC, this condition is checked across all $\mathcal{A}$ where $|\mathcal{A}|=k$ and  $\gamma_k$ denotes the proportion of such $\mathcal{A}$ that achieve both \eqref{eq:WSC1} and \eqref{eq:WSC2}. Their weak sign consistency criterion is defined to be
\begin{align}
    \text{WSC}=\text{mean}(\gamma_2,\dots,\gamma_{\lceil n/2 \rceil})\ .\ \label{eq:WSC}
\end{align}
A design that maximizes WSC is then highly likely to satisfy the two weak sign consistency conditions across many $\mathcal{A}$ where $|\mathcal{A}|\leq \lceil n/2 \rceil$. It would be preferred that the left hand sides of \eqref{eq:WSC1} and \eqref{eq:WSC2} equal 0, again corresponding to orthogonal designs, but nonorthogonal designs may still have a large WSC. The authors employ a hybrid elitist Pareto optimality-based coordinate exchange algorithm \citep{Caoetal2017} to identify designs on the Pareto front of the EC and WSC. We refer to their optimal designs as Dantzig-based Consistency Designs (DCDs). 

The DCD criteria require evaluations across all $\mathcal{A}$ of sizes $k=2,\dots,\lceil n/2 \rceil$. To ease the computational burden, \cite{singh2023selection} randomly sample subsets for each $k$. Furthermore, their examples initialized their search algorithms with only two designs each, a $\text{Var}(s+)$-optimal design and either a random, $\text{E}(s^2)$-optimal, or $\text{UE}(s^2)$-optimal design. The resulting DCDs initialized at a $\text{Var}(s+)$ design were shown via simulation to have better support recovery properties than the $\text{Var}(s+)$-optimal design assuming positive signs of the active effects. However, their improvement over the $\text{Var}(s+)$-optimal design may be attributed to their initialization at that design. That is, they show that the initialized $\text{Var}(s+)$-optimal design is not globally optimal for known signs, but there is no guarantee a DCD is either. Indeed, an orthogonal design is the ideal design under their criteria. It may be that the globally optimal DCD performs worse than the $\text{Var}(s+)$ design. Still, their results provide strong support for constructing SSDs under criteria directly tied to statistical properties of penalized estimators.

Others have considered the role of $\X$ in the statistical properties of the lasso, though not in the context of SSDs. \cite{zhao2006model} identified the strong irrepresentable condition (SIC) on $\X$ for establishing support recovery consistency of the lasso \citep{zhang2008sparsity, jia2015preconditioning}. 
\citet{wainwright2009sharp} and \citet{hastie2019statistical} have studied support recovery for random designs where the rows of $\X$ are independently generated from a zero-mean multivariate Normal distribution. Random designs, however, generally perform poorly for small run designs like SSDs. 
Based on the lasso's SIC, \cite{deng2013lasso} constructed SSDs from nearly orthogonal Latin hypercube designs to minimize the off-diagonals of $\Smat$ assuming factors have settings between $[0,1]$. This is essentially a construction technique for minimizing $\text{UE}(s^2)$ with more general factor settings, even though SSDs commonly assume fixed settings of $\pm 1$. \cite{phdthesis} proposed a lasso SIC criterion to construct optimal two-level SSDs, similar to the $\gamma_k$ from~\eqref{eq:WSC}. 
\cite{huang2020optimal} proposed a lasso optimal design theory that applies variance-based criteria to the approximate covariance matrix of the debiased lasso \citep{javanmard2014confidence}, which is capable of performing inference via confidence intervals. Under the framework of approximate designs (i.e., $n\to\infty$) they note their criteria are not convex and give an equivalence theorem for establishing whether a design is locally optimal. 
They then propose an algorithmic construction for generating many local optimal approximate designs and implement a rounding procedure on the approximate design's replication weights to produce an exact design. Overall, their approach requires many stages of approximation, leading to a discrepancy between the approximate and exact design's covariance matrix. 

\subsection{Contributions and Overview}

This paper addresses the two SSD open problems by developing optimality criteria that directly maximize the probability of sign recovery for the lasso estimator.

The first major contribution of this paper is a local optimal design approach to maximize the sign recovery probability assuming $\betavec$ is known and the lasso tuning parameter, $\lambda>0$, is fixed. 
We then develop criteria that relax the model assumptions and summarize the sign recovery probability across a range of $\lambda$. As optimization of the proposed criteria is challenging, the second major contribution is a criterion based on the $p \times p$ correlation matrix of $\X$. We provide evidence that the ideal correlation matrices under different heuristic criteria are optimal under the new criteria, thus making substantial progress on both open problems. In particular, we prove the surprising result that an orthogonal correlation matrix is suboptimal when $\zvec$ is known.  The final contribution is a computationally efficient construction algorithm that leverages the speed of exchange algorithms for optimizing heuristic criteria.

The paper is organized as follows. Section~\ref{sec:loc} develops exact design optimality criteria targeting the lasso's sign recovery probability under different assumptions about the model. 
Section~\ref{sec:optimal_cs} presents a pseudo-approximate optimal design framework that targets the optimal correlation matrix of $\X$ and justifies the ideal designs sought after by heuristic criteria under different assumptions about $\zvec$. Section~\ref{sec:EvalConstruct} describes a computationally efficient algorithm for constructing exact optimal designs under the proposed criteria. We then demonstrate the optimal SSDs under the new framework in Section~\ref{sec:NewDesigns}.
We conclude the paper with a discussion in Section~\ref{sec:Discussion}, describing important implications of the results and future work.

\section{Exact Local Optimality Criteria}\label{sec:loc}

The lasso is often applied to the centered and scaled design matrix 
$\Fmat=(\I-\Pmat_1)\X\V^{-1/2}$ where $\Pmat_1=n^{-1}\onevec\onevec^T$ and $\V$ is a diagonal matrix comprised of the diagonal elements of $n^{-1}\X^T(\I-\Pmat_1)\X$. 
The analysis then targets support/sign recovery of $\betavec^*=\V^{1/2}\betavec$. 
The diagonal elements of $\V^{1/2}$ are nonnegative and bounded above by 1, making $|\beta^*_j| \leq |\beta_j|$, and $\sign(\beta_j^*)=\sign(\beta_j)$ when $\V$ has all positive diagonal elements. The support, $\mathcal{A}$, is estimated by the support of the lasso estimator
\[
\hat{\betavec}^*
    =\argmin_{\betavec^*} \frac{1}{2}\betavec^{*T}\C\betavec^* -\frac{1}{n} \yvec^{T}\F\betavec^* + \lambda \sum_{j=1}^p |\beta_j^*|\ ,\ \label{eq:stdlasso}
\]
where $\C=n^{-1}\F^{T}\F$ is the correlation matrix of the columns of $\X$. 
Denote the estimated support and sign of the resulting lasso estimate, $\hat{\betavec}^*$, by $\hat{\mathcal{A}}=\{j : |\hat{\beta}_j^*| > 0\}$ and $\hat{\boldsymbol{\zvec}}=\sign(\hat{\boldsymbol{\beta}}^*)$, respectively.

We denote submatrices of $\X$ and $\F$ corresponding to column subsets $\mathcal{T} \subseteq \{1,\dots,p\}$ by $\X_\mathcal{T}$ and $\F_\mathcal{T}$, respectively. For a $p \times 1$ vector, $\vvec$, $\vvec_\mathcal{T}$ denotes the $|\mathcal{T}| \times 1$ vector with the $\mathcal{T}$ elements of $\vvec$. For all other matrices, we will consider submatrices by selecting both rows and columns with two index sets $\mathcal{U}$ and $\mathcal{T}$. That is, for a matrix $\M$,
let $\M_{\mathcal{U}\mathcal{T}}$ denote the submatrix of $\M$ with rows and columns indexed by $\mathcal{U}$ and $\mathcal{T}$, respectively. For brevity, we will denote $\M_{\mathcal{T}\mathcal{T}}=\M_\mathcal{T}$, which should not be confused with a subsetting of columns alone.

For known $\hat{\mathcal{A}}$ and $\hat{\zvec}$, the lasso KKT conditions are \citep{zhao2006model,tibshirani2012lasso}:
\begin{align}
\hat{\betavec}_{\hat{\mathcal{A}}}^* &= \frac{1}{n}\C_{\hat{\mathcal{A}}}^{-1}\F_{\hat{\mathcal{A}}}^{T}\yvec-\lambda \C_{\hat{\mathcal{A}}}^{-1}\hat{\zvec}_{\hat{\mathcal{A}}}\ ,\ \label{eq:ActiveEst}\\  \zerovec &< \hat{\Z}_{\hat{\mathcal{A}}}\hat{\betavec}_{\hat{\mathcal{A}}}^* \label{eq:SignCheck} \ ,\\\
\lambda \onevec &\geq \left|\C_{\hat{\mathcal{I}}\hat{\mathcal{A}}}\hat{\betavec}_{\hat{\mathcal{A}}}^*-\frac{1}{n}\F_{\hat{\mathcal{I}}}^{T}\yvec\right| \label{eq:InactiveCheck}\ ,\
\end{align}
where $\C_{\hat{\mathcal{A}}}^{-1}=n(\Fmat_{\hat{\mathcal{A}}}^T\Fmat_{\hat{\mathcal{A}}})^{-1}$, $\hat{\Z}_{\hat{\mathcal{A}}}=\text{Diag}(\hat{\zvec}_{\hat{\mathcal{A}}})$, and $\hat{\mathcal{I}}=\{j  :  \hat{\beta}_j^* = 0\}$. For given $\X$, $\yvec$, and $\betavec$, we can check whether the resulting $\hat{\betavec}^*$ has the true support and sign vector by setting $\hat{\mathcal{A}}=\mathcal{A}$ and $\hat{\zvec}_{\hat{\mathcal{A}}}=\zvec_{\mathcal{A}}$ and checking the KKT conditions. If they hold, $\hat{\betavec}^*$ recovers the sign, which is more stringent than recovering the support. 

The proposed criteria in this section rank designs according to the probability of sign recovery, denoted $P(\hat{\zvec} = \zvec \, | \, \X, \, \betavec, \, \lambda)$. To calculate the probability, one must specify a $\lambda >0$ and a $\betavec$, so the proposed framework resembles optimality frameworks employed for nonlinear models and maximum likelihood estimation \citep{Silvey1980,Khuri2006,YangStufken2009}.  Note that specifying $\betavec$ also leads to specifying $\mathcal{A}$, $\mathcal{I}$, and $\zvec$. The criterion $P(\hat{\zvec} = \zvec \, | \, \X, \, \betavec, \lambda)$ is a joint probability of two events corresponding to \eqref{eq:SignCheck} and \eqref{eq:InactiveCheck} assuming $\hat{\mathcal{A}}=\mathcal{A}$ and $\hat{\zvec}_{\hat{\mathcal{A}}}=\zvec_{\mathcal{A}}$. 
The first event corresponds to \eqref{eq:SignCheck} and checks whether $\sign(\hat{\betavec}^*_\mathcal{A})=\zvec_\mathcal{A}$ where $\hat{\betavec}^*_\mathcal{A}$ is from \eqref{eq:ActiveEst} assuming $\hat{\mathcal{A}}=\mathcal{A}$ and $\hat{\zvec}_{\hat{\mathcal{A}}}=\zvec_{\mathcal{A}}$:
\begin{align}
    S_\lambda 
    &\coloneqq
    \left\{\uvec < \sqrt{n} \Z_{\mathcal{A}}\V^{1/2}_\mathcal{A}\betavec_{\mathcal{A}}\right\}\ ,\ \label{eq:SignCheckUnscale}
\end{align}
where $\uvec \sim N(\lambda \sqrt{n} \, \Z_{\mathcal{A}}\C_\mathcal{A}^{-1}\zvec_{\mathcal{A}}, \Z_{\mathcal{A}}\C_\mathcal{A}^{-1}\Z_{\mathcal{A}})$. The second event corresponds to \eqref{eq:InactiveCheck}, which checks whether all $j \in \mathcal{I}$ have $\hat{\beta}^*_j=0$:
\begin{align}
I_\lambda 
\coloneqq \left\{|\vvec| \leq \lambda\sqrt{n}\onevec \right\}\ ,\ \label{eq:InactiveCheck2}
\end{align}
where $\vvec \sim N(\lambda\sqrt{n} \C_{\mathcal{I}\mathcal{A}}\C_\mathcal{A}^{-1}\zvec_{\mathcal{A}}, \C_\mathcal{I}-\C_{\mathcal{I}\mathcal{A}}\C_\mathcal{A}^{-1}\C_{\mathcal{A}\mathcal{I}})$. 
The probability of event~\eqref{eq:InactiveCheck2} depends on $\betavec$ only through $\zvec$. The probabilities of both events depend on $\X$ only through its $\C$ and $\V$, a fact we exploit in Section~3. For brevity, we will sometimes refer to these events as simply $S_\lambda$ and $I_\lambda$.


\textbf{Remark 1} There is a tradeoff between the probabilities of \eqref{eq:SignCheckUnscale} and \eqref{eq:InactiveCheck2} so both should be considered simultaneously. \cite{deng2013lasso} and \cite{phdthesis} focus on heuristic criteria based on the lasso's SIC, which optimizes \eqref{eq:InactiveCheck2} alone.

\textbf{Remark 2} 
The probability of support recovery, $P(\hat{\mathcal{A}}=\mathcal{A} \, | \, \X, \, \betavec, \, \lambda)$, requires consideration of all possible $2^{k}$ sign vectors that have the same $0$ elements as $\zvec$ but alternative $\pm 1$ elements indexed by $\mathcal{A}$. Defining $\mathcal{Z}_{\mathcal{A}}$ as the set of all such sign vectors, we have $P(\hat{\mathcal{A}}=\mathcal{A} \, | \, \X, \, \betavec, \, \lambda)=2^{-k}\sum_{\tilde{\zvec} \in \mathcal{Z}_{\mathcal{A}}} P(\hat{\zvec} = \tilde{\zvec} \, | \, \X, \, \betavec, \, \lambda)$. Each $P(\hat{\zvec} = \tilde{\zvec} \, | \, \X, \, \betavec, \, \lambda)$ may be calculated by replacing $\zvec_\mathcal{A}$ and $\zvec$ with $\tilde{\zvec}_\mathcal{A}$ and $\tilde{\zvec}$, respectively, in events~\eqref{eq:SignCheckUnscale} and \eqref{eq:InactiveCheck2}. 

\textbf{Remark 3} As $|\beta_j^*| \to \infty$ for all $j \in \mathcal{A}$, the probability of sign recovery and support recovery are equivalent. To see this, note that $\Z_{\mathcal{A}}\V^{1/2}_\mathcal{A}\betavec_{\mathcal{A}}=\V^{1/2}_\mathcal{A}|\betavec_{\mathcal{A}}|$ so $P(S_\lambda \, | \, \C_\mathcal{A}, \, \V_\mathcal{A}, \, \betavec_{\mathcal{A}}) \to 1$ since the area of integration approaches $\mathbb{R}^k$, the support of $\uvec$. For any other $\tilde{\zvec} \in \mathcal{Z}_\mathcal{A}$, $\tilde{\Z}_{\mathcal{A}}\V^{1/2}_\mathcal{A}\betavec_{\mathcal{A}}$ would have negative elements, causing $P(S_\lambda \, | \, \C_\mathcal{A}, \, \V_\mathcal{A}, \, \betavec_{\mathcal{A}}) \to 0$ because the region of integration for the $u_j$ corresponding to these negative elements approaches $(-\infty,\lim_{t\to-\infty} t)$, which has zero probability. Hence $P(\hat{\zvec} = \tilde{\zvec} \, | \, \X, \, \betavec, \, \lambda) \to 0$ for any $\tilde{\zvec} \neq \zvec$. Finally, $I_\lambda$ depends on $\betavec_\mathcal{A}$ only through $\zvec_\mathcal{A}$ so it is unaffected by increasing $|\beta_j^*|$, hence $P(\hat{\mathcal{A}}=\mathcal{A} \, | \, \X, \, \betavec) \to P(\hat{\zvec}=\zvec \, | \, \X, \, \betavec)$. 

\textbf{Remark 4} We assume $\C_{\mathcal{A}}^{-1}$ always exists, making $\uvec$ nondegenerate, but $\vvec$ will have a degenerate multivariate Normal distribution for SSDs because $\C$ cannot be full rank. The probability of this event can be calculated following some linear transformation of $\vvec$.

\subsection{Local criteria assuming known $\betavec$}\label{sec:localcrit}

For a fixed $\lambda$ and known $\betavec$, define the local optimality criterion
\begin{align}
    \phi_\lambda(\X \, | \, \betavec)=P(\hat{\zvec}=\zvec \, | \, \X, \, \betavec )=P(S_\lambda \cap I_\lambda \, | \, \C, \, \V_\mathcal{A}, \,  \betavec)\ .\ \label{eq:localcrit}
\end{align}
A $\phi_\lambda$-optimal design is $\X^* = \text{argmax}_{\X} \  \phi_\lambda (\X \,|\, \betavec)$. This approach is impractical, particularly due to its perfect knowledge of $\mathcal{A}$, but it is foundational for the more practical criteria that allow for uncertainty about $\mathcal{A}$. The following is a fundamental result about the role of $\zvec$ in $\phi_\lambda(\X \, | \, \betavec)$. Its proof and all future proofs may be found in the Supplementary Materials. 
\begin{lemma} \label{lem:symmetry}  For a given $\X$, events $S_\lambda$ and $I_\lambda$ are independent, and $\phi_\lambda(\X \, | \, \betavec)=\phi_\lambda(\X \, | \, -\betavec)$ for any $\betavec$ and its reflection, $-\betavec$. Hence an optimal $\X^*$ for $\phi_\lambda$ under $\betavec$ is also optimal for $\phi_\lambda$ under $-\betavec$.

 \end{lemma}
 \noindent
Lemma~\ref{lem:symmetry}  includes a design equivalence result between $\betavec$ and its reflection, $-\betavec$. The following equivalence theorem is more general and further simplifies the implementation of the framework:

\begin{theorem}
For a $\betavec$ where $\zvec_\mathcal{A}=\onevec$, consider $\tilde{\betavec}=\tilde{\Z}\betavec$ where $\tilde{\Z}$ is any diagonal matrix comprised of $\pm 1$. Then for any $\X$, $\phi_\lambda(\X \, | \, \betavec)=\phi_\lambda(\X\tilde{\Z} \, | \, \tilde{\betavec})$. Consequently, if $\X^*$ is optimal for $\phi_\lambda$ under a $\betavec$, then $\tilde{\X^*}=\X^*\tilde{\Z}$ is optimal for $\phi_\lambda$ under $\tilde{\betavec}.$
\label{thm:localequiv}
\end{theorem}
\noindent
A consequence of Theorem~\ref{thm:localequiv} is that, when discussing criteria involving a known sign vector, we can assume $\zvec_\mathcal{A}=\onevec$ without loss of generality.

The following criterion assumes the magnitudes of the $\beta_j$ are known and treats $\zvec_\mathcal{A}$ as an unknown quantity and calculates the expected probability of sign recovery assuming all $\tilde{\zvec} \in \mathcal{Z}_\mathcal{A}$ are equally likely. Without loss of generality, the magnitudes are expressed by $\betavec > \zerovec$ and the criterion is defined to be
\begin{align}
\phi_{\lambda}^{\pm}(\X \, | \, \betavec) = \frac{1}{2^k} \sum_{\tilde{z} \in \mathcal{Z}_\mathcal{A}} \phi_\lambda (\X \, | \, \tilde{\Z}\betavec)=\frac{1}{2^k} \sum_{\tilde{z} \in \mathcal{Z}_\mathcal{A}} \phi_\lambda (\X\tilde{\Z} \, | \, \betavec)\ .\ \label{eqn:localA_allsigns}
\end{align}
Calculating all $2^k$ probabilities can be computationally intensive, but Lemma~\ref{lem:symmetry} allows us to halve the number of computations. We state this as a corollary:
\begin{corollary}
For a fixed $\mathcal{A}$ where $|\mathcal{A}|=k$, let $\mathcal{Z}_{\mathcal{A}}^{\pm}$ denote the subset of $\mathcal{Z}_\mathcal{A}$ including all $2^{k-1}$ unique $\zvec$ up to reflection. Then $\phi_{\lambda}^{\pm}(\X \, | \, \betavec) = 2^{-(k-1)}\sum_{\tilde{z} \in \mathcal{Z}_\mathcal{A}^{\pm}} \phi_\lambda (\X\tilde{\Z} \, | \, \betavec)$.
\end{corollary}
\noindent

We now investigate designs that maximize $\phi_\lambda$ and $\phi_\lambda^{\pm}$. Knowledge of $\mathcal{A}$ when constructing an optimal design under either criteria leads to a trivial construction for the columns $\X_\mathcal{I}^*$.
\begin{proposition}
For a known $\mathcal{A}$, there exists an optimal design $\X^*$ where $\X^*_\mathcal{I}=\J$. Furthermore, the criteria $\phi_\lambda(\X^* \, | \, \betavec)$ and $\phi_\lambda^{\pm}(\X^* \, | \, \betavec)$ for such designs depend only on $S_\lambda$.
\end{proposition}
\noindent
Optimal designs for criteria with perfect knowledge about $\mathcal{A}$ completely confound the columns of the inactive factors with the intercept, making $P(I_\lambda \, | \, \C, \, \zvec)=1$. Finding the optimal design, which may originally have been supersaturated, then reduces to finding a full column rank matrix $\X_\mathcal{A}$ that maximizes probability of the $S_\lambda$ event(s). 

To demonstrate this new theory, we will now show that if $\zvec_\mathcal{A}=\onevec$ is assumed known, an orthogonal design may not maximize $P(S_\lambda \, | \, \C_\mathcal{A},\, \V_\mathcal{A},\, \betavec_\mathcal{A})$. For any $\X_\mathcal{A}$ the bounds of integration for $S_\lambda$ and the mean of $\uvec$ are proportional to $\V_\mathcal{A}^{1/2}\betavec_\mathcal{A}$ and $\lambda \C_\mathcal{A}^{-1}\onevec$, respectively. Hence, a design that maximizes $\V_\mathcal{A}^{1/2}\betavec_\mathcal{A} - \lambda \C_\mathcal{A}^{-1}\onevec$ may maximize $P(S_\lambda \, | \, \C_\mathcal{A},\, \V_\mathcal{A},\, \betavec_\mathcal{A})$, although the covariance matrix of $\uvec$ also plays a role. For any orthogonal design, this difference is $\betavec_\mathcal{A}-\lambda\onevec$.  We will try to find a $\X_\mathcal{A}$ with larger $P(S_\lambda \, | \, \C_\mathcal{A},\, \V_\mathcal{A},\, \betavec_\mathcal{A})$ by identifying a design whose 
$\V_\mathcal{A}^{1/2}$ and $\C_\mathcal{A}$ satisfy
\begin{align}
(\I - \V_\mathcal{A}^{1/2})\betavec_\mathcal{A} \leq (\I - \C_\mathcal{A}^{-1})\lambda \onevec\ .\ \label{eqn:Sbound}
\end{align}

As $(\I-\V_\mathcal{A}^{1/2})\betavec_\mathcal{A} \geq \zerovec$, we need $(\I - \C_\mathcal{A}^{-1})\lambda \onevec \geq \zerovec$. For even $n \geq 6$, construct $\X_\mathcal{A}$ by choosing $k \leq n-1$ columns from the $n \times n$ matrix
\begin{align}
\left(\begin{array}{c|c}
2\I-\J & -\J \\ \hline \J & \J-2\I
\end{array}\right)\ ,\ \label{eqn:LODbetter}
\end{align}
where all $\I$ and $\J=\onevec\onevec^T$ have size $n/2 \times n/2$. Then $\V_\mathcal{A}=(1-4/n^2)\I$. Let $k_1$ denote the number of columns chosen from the first $n/2$ columns and $k_2=k-k_1$ be the number from the remaining columns. If either $k_1$ or $k_2$ equal 0, $\C_\mathcal{A}$ will be completely symmetric with positive off-diagonal elements $c=1-4n/(n^2-4)$. Thus $(\I - \C_\mathcal{A}^{-1})\lambda \onevec \geq \zerovec$ and \eqref{eqn:Sbound} becomes
\[
\left(1-\frac{\sqrt{n^2-4
}}{n}\right)\betavec_\mathcal{A} \leq \lambda \left(1-\frac{n^2-4}{kn^2-4k(n+1)+4n}\right)\onevec\ .\
\]

\noindent
When $k_1 \geq 1$ and $k_2\geq1$,  $\C_\mathcal{A}$ and $\C_\mathcal{A}^{-1}$ have a $2 \times 2$ block partitioned form with completely symmetric, diagonal block matrices and constant off-diagonal block matrices. Then $\C_\mathcal{A}^{-1}\onevec=(\xi_1 \onevec_{k_1}^T, \xi_2 \onevec_{k_2}^T)^T=\boldsymbol{\xi}$, and $\C_\mathcal{A}\boldsymbol{\xi}=\onevec$, which gives the equations
\begin{align*}
1 &= \rho_{11} \xi_1 + \rho_{12} \xi_2\\
1 &= \rho_{21} \xi_1 + \rho_{22}\xi_2\ ,\
\end{align*}
where $\rho_{ij}>0$ is the unique row sum for the corresponding $k_i \times k_j$ block matrix of $\C_\mathcal{A}$. Defining $\tilde{k}_i=(n-2k_i)\geq 0 $, we have
\begin{equation}
\xi_1 = \frac{\tilde{k}_2(n^2-4)}{(n-2)^2(k_1\tilde{k}_2+\tilde{k}_1k_2)+\tilde{k}_1\tilde{k}_2} \geq 0\ ,\ \label{eq:delta1}
\end{equation}
with equality if and only if $k_2=n/2$. A similar expression holds for $\xi_2$ with $\tilde{k}_1$ in the numerator of \eqref{eq:delta1}, and $\xi_2=0$ if and only if $k_1=n/2$. Since the diagonal elements of $\C_\mathcal{A}$ are 1, 
$\rho_{ii} \geq 1$ and it follows 
$\onevec-\boldsymbol{\xi}>\zerovec$. Hence designs constructed in this way satisfy \eqref{eqn:Sbound}
for some combinations of $\betavec_\mathcal{A}$ and $\lambda$, but whether this implies higher $P(S_\lambda \, | \, \C_\mathcal{A},\, \V_\mathcal{A},\, \betavec_\mathcal{A})$ depends on the covariance matrix of $\uvec$.

\begin{figure}[ht]
    \centering
\includegraphics[width=0.6\textwidth]{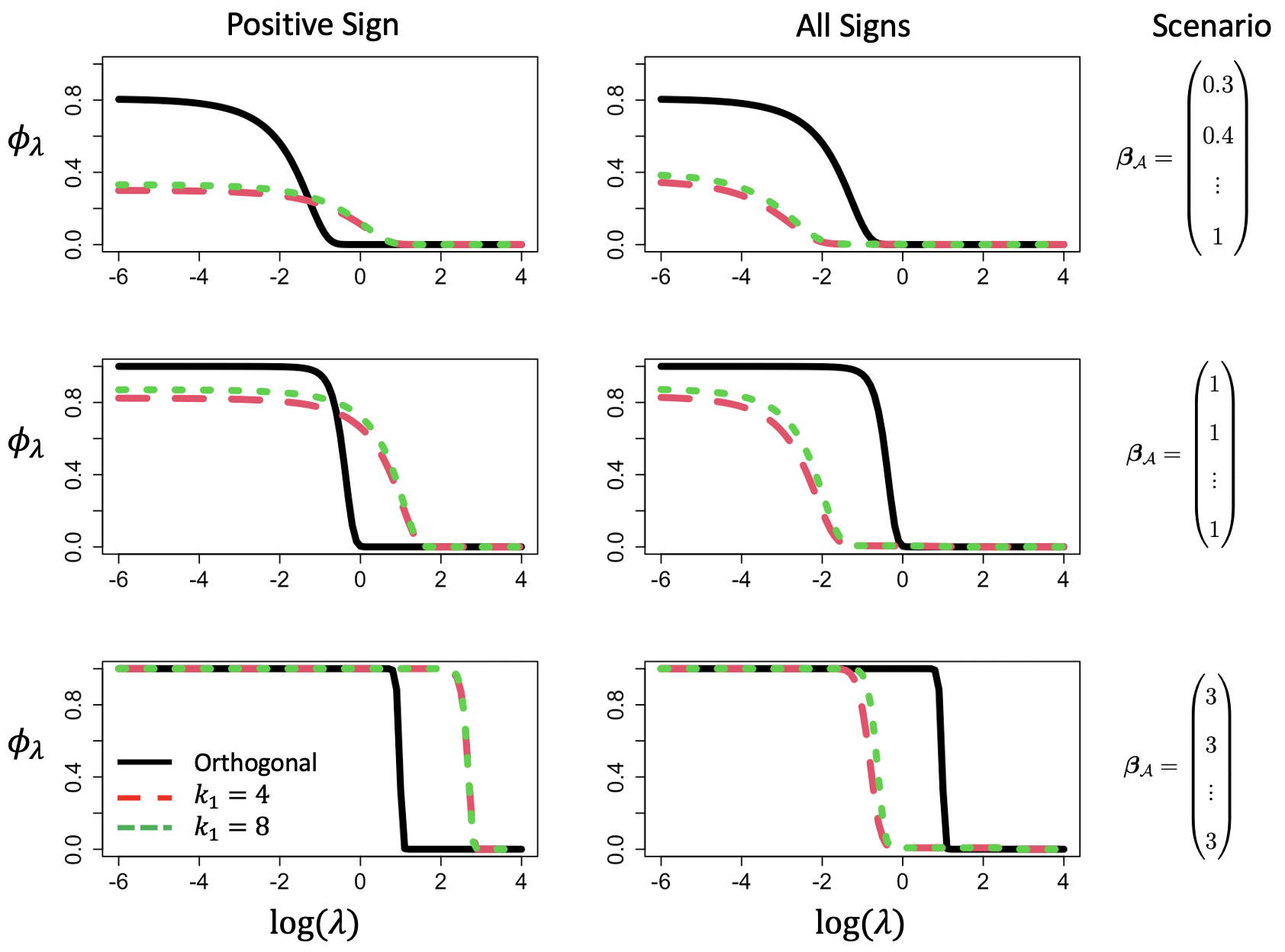}
    \caption{Probability of sign recovery for $n=16$, $p=20$, and $k=8$ for an orthogonal design and two designs constructed from selecting columns from \eqref{eqn:LODbetter} with $k_1=4$ and $8$. The left and right panels correspond to $\phi_\lambda$ and $\phi_\lambda^{\pm}$, respectively.}
    \label{fig:trivialcase}
\end{figure}

For example, suppose $n=16$ and we have $p=20$ factors but we know the first $k=8$ factors have nonzero effects and the rest are inactive.  We considered three SSDs that set all inactive factors' columns to $+1$, eliminating any ability to estimate them. The first SSD set the $8$ active factor columns to an orthogonal $\X_\mathcal{A}$ and the other two SSDs set active effects columns according to \eqref{eqn:LODbetter} with $k_1=4$ and $8$, respectively. For all designs, $\C_\mathcal{A}^{-1}\onevec=\xi \onevec$ where $\xi=1$, $0.1575$, and $0.1607$ for the orthogonal, $k_1=4$, and $k_1=8$ designs, respectively. For the $k_1=4$ and $k_1=8$ designs, \eqref{eqn:Sbound} is equivalent to $\lambda^{-1}\betavec_\mathcal{A} \leq 53.92 \times \onevec$ and $\lambda^{-1}\betavec_\mathcal{A} \leq 53.72 \times \onevec$, respectively. 
The three designs' $\phi_\lambda$ and $\phi_\lambda^{\pm}$ values were compared across a range of $\lambda$ and three scenarios for $\betavec_{\mathcal{A}}$: $(0.3,0.4,\dots,1)^T$, $\onevec$, and $3 \times \onevec$.  
The results for $\phi_\lambda$ are shown in the left panels of Figure~\ref{fig:trivialcase}.

For all scenarios, there is a range of large $\lambda$ values where all designs have $\phi_\lambda = 0$, a middle range where the two proposed designs outperform the orthogonal design, and a range of small $\lambda$ values where the orthogonal design is superior. The orthogonal design improves over the other two designs before condition~\eqref{eqn:Sbound} is violated, due to the role of the covariance matrices. The improvement is negligible when $\betavec_\mathcal{A}=3\times\onevec$, implying the covariance matrix is less important as the elements of $\betavec_\mathcal{A}$ increases. The results for $\phi_\lambda^{\pm}$ are shown in the right panels of Figure~\ref{fig:trivialcase}, and clearly favor the orthogonal design. 

\subsection{Relaxed local criteria}\label{sec:relax}
To make the criteria in Section~2.1 more practical, we introduce relaxations regarding the choice of $\lambda$ and the assumptions of the underlying model.  To reduce the dependency on $\lambda$, we propose two summary measures of $\phi_\lambda$ and $\phi_\lambda^{\pm}$. The first summary takes the maximum sign recovery probability: $\phi_{\text{max}}(\X \, | \, \betavec)=\max_{\lambda > 0} \phi_\lambda(\X \, | \, \betavec)$. However, in practice one will perform some tuning parameter selection strategy to choose a $\lambda$, and so the benefits of a design that optimizes $\phi_{\text{max}}$ will only be realized if the strategy reliably picks the corresponding $\lambda$ that maximizes $\phi_\lambda$. Hence, we would like a design that maintains large probabilities for a wide range of $\lambda$. It is common to calculate the lasso solution path with respect to $\log(\lambda)$. Such a transformation stretches the region of $\lambda \in (0,1)$, whose corresponding lasso estimates would receive the least amount of penalization. Therefore, we propose the criterion 
\[
    \phi_\Lambda(\X \, | \, \betavec)=\int_0^\infty \frac{\phi_\lambda(\X \, | \, \betavec)}{\lambda}\, d\lambda=\int_{-\infty} ^\infty \phi_{\exp(\omega)}(\X \, | \, \betavec) d\omega\ .\
\]
The definitions for $\phi_{\max}^{\pm}(\X \, | \, \betavec)$ and $\phi_\Lambda^{\pm}(\X \, | \, \betavec)$  relax assumptions about $\zvec_\mathcal{A}$ and measure the maximum expected sign recovery probability and the integrated expected sign recovery probability, assuming all $\zvec_\mathcal{A}$ are equally likely.

Turning to relaxing assumptions about the $\betavec$, 
we first assume $|\betavec_\mathcal{A}| \geq \beta \onevec$ where $\beta$ is the minimum practically significant effect size. Such a $\beta$ value is commonly chosen while performing power calculations so we will calculate the support recovery probability under the most challenging scenario of $|\betavec_\mathcal{A}| = \beta \onevec$. We will defer the discussion about incorporating uncertainty about $|\betavec_\mathcal{A}|$ to Section~\ref{sec:Discussion}. Finally, to address uncertainty about $\mathcal{A}$, we fix $k \leq n-1$ and assume all supports of this size are equally likely. Denote this set of supports by $\mathcal{A}_{k}$. 

With the above relaxed assumptions, we first introduce a criterion that incorporates uncertainty about $\mathcal{A}$, but treats $\zvec_{\mathcal{A}}$ as known. This resembles scenarios where $\text{Var}(s+)$-optimal designs have been shown to outperform $\text{E}(s^2)-$ and $\text{UE}(s^2)-$optimal SSDs in simulation studies. Assume $z_j$ depends only on whether $j \in \mathcal{A}$ or $\mathcal{I}$ as follows.  Let $\zvec^*$  be the $p \times 1$ vector comprised of $\pm 1$ for each $j$, assuming $j \in \mathcal{A}$. Then $z_j=z_j^*$ when $j \in \mathcal{A}$ and $z_j=0$ when $j \in \mathcal{I}$. Following Theorem~\ref{thm:localequiv}, we assume without loss of generality that $\zvec^*=\onevec$ to identify the optimal design and then multiply the columns of the design by their actual corresponding $z_j^*$. For a given $\mathcal{A}$, let $\A$ be the $p \times p$ diagonal matrix of all zeroes except for the diagonal elements corresponding to $\mathcal{A}$, which are set to 1. Then the true sign vector is $\zvec=\A\onevec$ and the sign-dependent criterion for a given $\beta$ is
\[
    \Phi_{\lambda}(\X \, | \, k, \beta) = \binom{p}{k}^{-1} \sum_{\mathcal{A} \in \mathcal{A}_k} \phi_{\lambda}(\X \, | \, \betavec=\beta\A\onevec)\ ,\ 
\]
being the expected probability across $\mathcal{A}_k$ for a fixed $\lambda$. The notation $\Phi$ is intended to reflect consideration of all $\mathcal{A}\in\mathcal{A}_k$. The sign-independent criterion, $\Phi_{\lambda}^{\pm}(\X \, | \, k, \beta)$, is similarly defined, considering all $2^k$ sign vectors as equally likely. The two summary measures originally defined on $\phi_\lambda$ and $\phi_\lambda^{\pm}$ are straightforward to define on these new criteria and share similar notation, i.e., $\Phi_{\text{max}}(\X \, | \, k, \beta)$ and $\Phi_{\Lambda}(\X \, | \, k, \beta)$. Computational details of the criteria may be found in the Supplementary Materials.

Optimizing any of these sign recovery criteria is both analytically and algorithmically challenging. Indeed, simply evaluating such criteria for a given $\X$ can be cumbersome. We next discuss properties of these criteria and, motivated by approximate design theory \citep{pukelsheim2006optimal}, we propose a more tractable framework centered on finding optimal $\C$ matrices.  

\section{Approximate Local Optimality Criteria} \label{sec:optimal_cs}

Popular design criteria for least-squares estimation are functions of the eigenvalues of $\M=\X^T(\I-\Pmat_1)\X$. Such criteria are difficult to optimize for fixed $n$ because of the discrete nature of the problem. A more tractable approach optimizes the criteria but with an extended domain of all symmetric, positive definite matrices, denoted $\mathcal{M}$. Such ``approximate" criteria are log concave with respect to $\mathcal{M}$, and are invariant to simultaneous row/column permutations and sign transformations of $\M$. These properties can be exploited to find optimal forms of $\M$. For an $\M \in \mathcal{M}$, define its matrix-averaged form to be $\overline{\M}=(p!)^{-1} \sum_{\Pimat \in \mathcal{P}} \Pimat \, \M \, \Pimat^T$ where $\mathcal{P}$ is the set of all $p \times p$ permutation matrices.  Then $\overline{\M} \in \mathcal{M}$ is a completely symmetric matrix, having constant diagonal elements and constant off-diagonal elements. Further averaging of $\overline{\M}$ across all $2^p$ sign transformations leads to a matrix proportional to the identity matrix. For a log concave criterion, the criterion value for $\overline{\M}$ is greater than or equal to that for $\M$. Hence the search for the optimum $\M$ may be restricted to the class of completely symmetric matrices in $\mathcal{M}$. In this class, an optimal $M^*$ can be identified, and one can search among exact designs to find ones whose $\M$'s are close to $\M^*$ \citep{Pukelsheim1992}. 

The least-squares approximate criteria framework can be extended to the lasso criteria in Section~\ref{sec:loc}. Recall the events~\eqref{eq:SignCheckUnscale} and \eqref{eq:InactiveCheck2} depend on $\X$ only through its $\C$ and $\V$ matrix. We can then extend the domain of the lasso criteria to the pair of an arbitrary correlation matrix and an arbitrary diagonal matrix with elements between 0 and 1. One can then solve for an optimal pairing, $(\C^*,\V^*)$ and rank exact designs based on some measure of proximity of their $(\C, \V)$ to $(\C^*,\V^*)$. 

For brevity, the approximate lasso criteria in this section absorb $\sqrt{n}$ into the values of $\betavec_\mathcal{A}$ and $\lambda$. The $\V$ matrix contributes exclusively to the event $S_\lambda$, and it is evident that $\V^* = \I$ maximizes $P(S_\lambda \, | \, \C_\mathcal{A}, \, \V_\mathcal{A}, \, \betavec_{\mathcal{A}})$ over all possible $\V$, as it maximizes the area of integration. 
Such a $\V$ is possible in practice, occurring if and only if the columns of $\X$ are balanced, a requirement of $\text{E}(s^2)$-optimal designs. Therefore, we assume $\V=\I$ and consider optimization of the criteria with respect to $\C$. To make this explicit, we will rewrite the previous criteria using $\C$ in place of $\X$, e.g., $\phi_\lambda(\X \, | \, \betavec)$ becomes $\phi_\lambda(\C \, | \, \betavec)$.

To simplify our investigation of these new approximate lasso criteria, we condition optimization across the space of positive definite correlation matrices.  Since such matrices cannot occur with an SSD, a more appropriate framework tailored for SSDs would require optimization across rank-deficient correlation matrices, which to our knowledge has received little attention in the statistical literature. However, similar to the motivation behind the $\text{E}(s^2)$- and $\text{UE}(s^2)$-criterion, an optimal positive definite correlation matrix still provides a hypothetical target which we use in Section~\ref{sec:EvalConstruct} to efficiently construct new SSDs to maximize lasso sign recovery probabilities.

The primary goal of this section is to identify optimal, positive definite correlation matrices under the relaxed local criteria. Such criteria are permutation invariant due to $|\betavec_\mathcal{A}|=\beta \onevec$, which suggests there exists an optimal $\C^*$ that is completely symmetric. For $p \geq k$, let $\mathcal{C}_{p,k}=\{\C = (1-c)\I_p +c\J_p \, | \, -(k-1)^{-1}<c<1\}$ denote a set of $p \times p$ completely symmetric correlation matrices. This space does include rank-deficient correlation matrices, corresponding to $-(k-1)^{-1}<c\leq -(p-1)^{-1}$. For $k=p$, we will denote this set as simply $\mathcal{C}_p$. If the approximate criteria were log concave across the space of all positive definite correlation matrices, then, similar to eigenvalue-based criteria, matrix averaging implies the optimal correlation matrix must lie in $\mathcal{C}_p$. This reduces the dimension of the optimization problem from $\binom{p}{2}$ off-diagonal elements of $\C$ to a single off-diagonal value, $c$.  Investigating log concavity for approximate criteria is challenging because $\C$ is involved in both the covariances and means for $S_\lambda$ and $I_\lambda$. Derivatives are intractable even for a fixed $\lambda$.  

To simplify the investigation of log concavity, Section 3.1 assumes a known $\betavec$ and the optimal design has $\X_\mathcal{I}=\J$. This makes the criterion equal to $P(S_\lambda \, | \, \C_\mathcal{A}, \, \V_\mathcal{A}, \, \betavec_{\mathcal{A}})$ and so we explore the log concavity of this probability with respect to $\C_\mathcal{A}$. We identify situations when such functions fail to be log concave, but provide evidence that the probability may still be optimized by $\C_\mathcal{A} \in \mathcal{C}_k$. These numerical results motivate focusing optimization of relaxed criteria across $\mathcal{C}_{p,k}$ in Section~3.2.

\subsection{Optimal $\C_{\mathcal{A}}$ for known $\betavec$}\label{sec:OptimalCA}


Following Proposition~1, approximate local criteria for a known $\betavec$ should assume $\X_\mathcal{I}=\J$ and so should be defined only on the space of positive definite $k \times k$ correlation matrices, i.e., $\phi_\lambda(\C_\mathcal{A} \, |\, \betavec_\mathcal{A})=P(S_\lambda \, | \, \C_\mathcal{A}, \, \I, \, \betavec_{\mathcal{A}})$ and $\phi_\lambda^{\pm}(\C_\mathcal{A} \, |\, \betavec_\mathcal{A})=2^{-(k-1)}\sum_{\tilde{z} \in \mathcal{Z}^{\pm}_\mathcal{A}}P(S_\lambda \, | \, \tilde{\Z}\C_\mathcal{A}\tilde{\Z}, \, \I, \, |\betavec_{\mathcal{A}}|)$. To focus on relaxed criteria, we simplified our investigation of log concavity to the case of $\betavec_\mathcal{A}=\beta \onevec_k$. While such criteria are impractical by themselves, they contribute to understanding of more useful relaxed local criteria. The purpose of this section is to investigate when such criteria are log concave and, if they are not, when they may be optimized by a $\C_\mathcal{A} \in \mathcal{C}_k$. We focus on $k=2, 3,$ and $8$, and show that $\phi_\lambda(\C_\mathcal{A} \, |\, \betavec_\mathcal{A})$ is generally optimized by a $\C_\mathcal{A} \in \mathcal{C}_k$, even when it fails to be log concave (which tends to only happen as $\beta$ and $\lambda$ approach 0). The criteria $\phi_\lambda^\pm(\C_\mathcal{A} \, | \, \betavec_\mathcal{A})$ often fail to be log concave for larger values of $\lambda$ but still we are able to find a $\C_\mathcal{A} \in \mathcal{C}_k$ that maximizes (or nearly maximizes) $\phi_\lambda^\pm(\C_\mathcal{A} \, | \, \betavec_\mathcal{A})$.



\subsubsection{Examples with $k=2$}
\begin{figure}[h]
    \centering
    \includegraphics[width=0.4\linewidth]{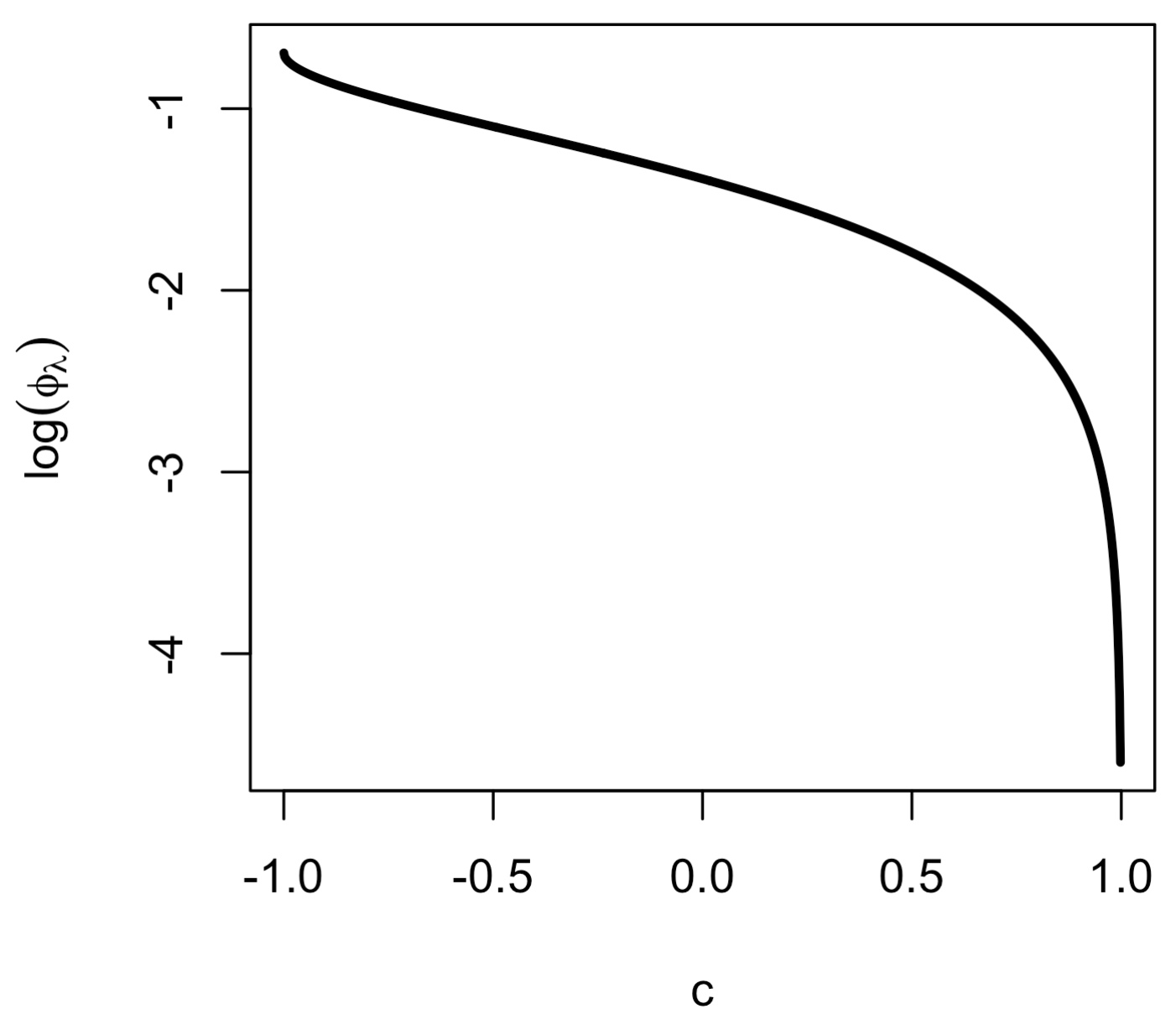}
    \caption{Plot of the logarithm of equation~\eqref{eqn:PSk2_equal_limit}. The function fails to be concave for $c<-0.4422$.}
    \label{fig:limitingphi}
\end{figure}

With $k=2$ we have an optimization problem of a single correlation value, $|c|<1$, so $\C_\mathcal{A}^*$ must be completely symmetric. Establishing log concavity in this case is essentially investigating the uniqueness of the optimal $c^*$, although log concavity is only a sufficient condition for uniqueness. In the Supplementary Materials, we show
\begin{align}
\phi_\lambda(\C_\mathcal{A} \, | \, \betavec_\mathcal{A})&=[2G(\mu(c))-1]\times G((1+c)\beta-\lambda) + 2 T\left((1+c)\beta-\lambda, \sqrt{\frac{1-c}{1+c}}\right) \label{eqn:PSk2_equal}\ ,\
\end{align}
where $\mu(c)=\sqrt{(1-c)/(1+c)}((1+c)\beta-\lambda)$; $g(\cdot)$ and $G(\cdot)$ refer to the standard Normal probability density function and cumulative distribution function, respectively; and $T(\cdot,\cdot)$ is Owen's $T$ function. There is no closed form expression for $T(\cdot,\cdot)$ except in special cases. One such case results from taking the limit:
\begin{align}
    \lim_{\beta,\lambda \to 0} \phi_\lambda(\C \, | \, \betavec) =2 T\left(0,\sqrt{\frac{1-c}{1+c}}\right)=\frac{1}{\pi}\arctan\left(\sqrt{\frac{1-c}{1+c}}\right)\ .\ \label{eqn:PSk2_equal_limit}
\end{align}
The log of \eqref{eqn:PSk2_equal_limit} is shown in Figure~\ref{fig:limitingphi} which is decreasing in $c$ and convex for $c<-0.4422$. This shows that $P(S_\lambda \, | \, \C_\mathcal{A}, \, \I, \, \betavec_{\mathcal{A}})$ is not log concave for all $\betavec_\mathcal{A}$ and $\lambda$, although the optimal $c$ is unique. 

The limiting case in \eqref{eqn:PSk2_equal_limit} is uninteresting because practical screening problems typically assume $\beta \geq 1$. Moreover, optimal designs under relaxed criteria will not set $\X_\mathcal{I}=\J$, and for such designs, $P(I_\lambda\, | \, \C, \zvec)$ will be close to 0 as $\lambda \to 0$, thereby eliminating the benefits of setting $c$ close to $-1$. Finally, for $\beta$ close to 0, sign recovery probability will be significantly smaller than the support recovery probability, and the latter should be preferred for comparing designs.

\begin{figure}[h]
    \centering
    \includegraphics[width=0.75\linewidth]{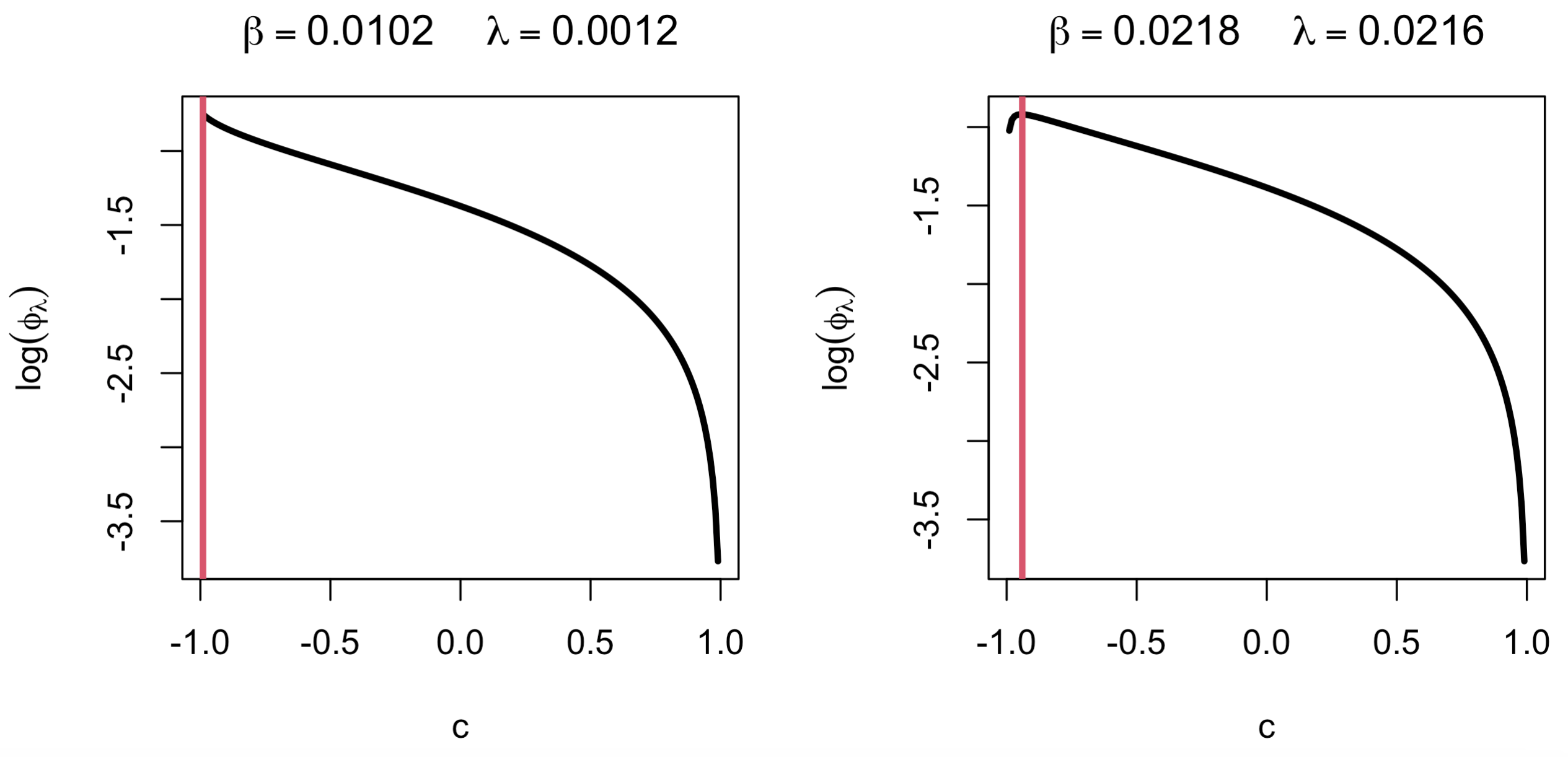}
    \caption{Plot of the logarithm of equation~\eqref{eqn:PSk2_equal} for two of the 750 maximin scenarios that violate log concavity. The red line corresponds to the $c^*$ value, being $-0.99$ and $-0.94$ for the left and right panels, respectively.}
    \label{fig:PSk2_onesign}
\end{figure}

A general assessment of log concavity may be accomplished by analyzing the second-order derivatives of $\log[\phi_\lambda(\C_\mathcal{A} \, | \, \betavec_\mathcal{A})]$. However, the derivatives are tedious and too complicated for establishing concavity. Another technique to establish log concavity of a function involving an integral is based on the Prékopa–Leindler inequality \citep{LogConcaveIntegral} which requires the integrand to be bivariate log concave. An example is provided in the Supplementary Materials to show this sufficient condition does not hold for the integrands involved in $P(S_\lambda \, | \, \C_\mathcal{A}, \, \I, \, \betavec_{\mathcal{A}}=\beta \onevec)$. Instead, we numerically investigated the log concavity conjecture under various combinations of $\betavec_\mathcal{A}$ and $\lambda$. We constructed a maximin design with 750 points to fill the space of possible $\beta$ and $\log(\lambda)$, where $\beta \in [0.001,3]$ and $\log(\lambda) \in [-10,1.2]$. For each setting, we estimated the second-order derivatives of $\log(\phi_\lambda(\C_\mathcal{A} \, | \, \betavec_\mathcal{A}))$ for a dense grid of $c \in [-0.99,0.99]$ using second-order finite differences with a step size of $h=0.01$. If the largest value was greater than 0 the corresponding criterion cannot be log concave.   Of the 750 settings considered, 49 were not log concave. These cases resembled the limiting case from \eqref{eqn:PSk2_equal_limit}, having $\lambda < 0.01$ and $\beta \leq 0.3720$. The optimal $c$ value in all of these cases was unique and close to $c^*=-0.99$. Two examples are shown in Figure~\ref{fig:PSk2_onesign}.

%



\begin{figure}[h]
    \centering
    \includegraphics[width=0.75\linewidth]{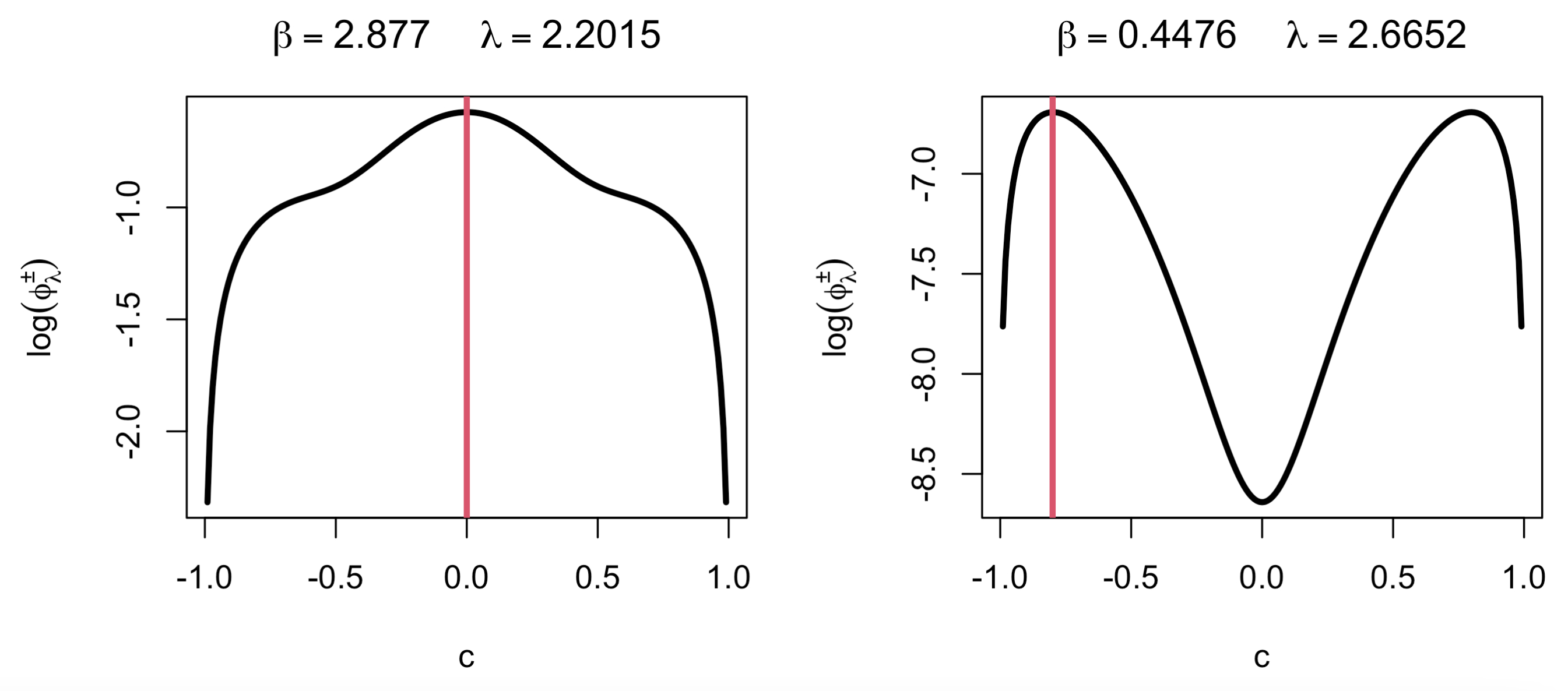}
    \caption{Plot of $\log[\phi_\lambda^\pm(\C_\mathcal{A} \, | \, \betavec_\mathcal{A})]$ for two of the 750 maximin scenarios that violate log concavity. The red line corresponds to the $c^*$ value, being $0$ and $-0.80$ for the left and right panels, respectively.}
    \label{fig:PSk2_allsigns}
\end{figure}

In the Supplementary Materials, we show $\log[\phi_\lambda^\pm(\C_\mathcal{A} \, | \, \betavec_\mathcal{A})]$ is symmetric about $c=0$. For it to be log concave, it is necessary (but not sufficient) for $c^*=0$. Of the 750 maximin settings, 81 failed to be log concave across a range of $\beta \in (0.01,3)$ and $\log(\lambda) \in (-9.854,1.198)$. Two examples that violate log concavity are shown in Figure~\ref{fig:PSk2_allsigns}, with one case having $c^*=0$. The lack of log concavity is not surprising, as the criterion is an average of functions, which can fail to be log concave even when the individual functions themselves are log concave. If log concavity were a necessary property, one could replace the average with the minimum probability across all sign vectors. We provide some discussion on this alternative definition in Section~\ref{sec:Discussion}. 

\subsubsection{Examples with $k=3$}

Maximizing $\phi_\lambda(\C_\mathcal{A} \, | \, \betavec_\mathcal{A})$ when $k=3$ requires optimization across three correlation parameters, $(c_{12}, c_{13},c_{23})$. Such criteria are permutation invariant and so, if they are log concave, then $\log[\phi_\lambda(\C_\mathcal{A} \, | \, \beta_\mathcal{A})] \leq \log[\phi_\lambda(\bar{\C}_\mathcal{A} \, | \, \beta_\mathcal{A})]$.  If there exists a $\C_\mathcal{A}$ where this does not hold, then $\phi_\lambda(\C_\mathcal{A} \, | \, \betavec_\mathcal{A})$ cannot be log concave. For each of the 750 settings from Section~3.1.1, we randomly generated 500 full-rank correlation matrices and checked how many violated $\log[\phi_\lambda(\C_\mathcal{A} \, | \, \beta_\mathcal{A})] \leq \log[\phi_\lambda(\bar{\C}_\mathcal{A} \, | \, \beta_\mathcal{A})]$.  We also recorded the best observed $\C_\mathcal{A}$ matrix and compared its criterion value to $\C_\mathcal{A}^* \in \mathcal{C}_3$, the optimal correlation matrix in $\mathcal{C}_3$.

Of the 750 settings considered, 87 found at least one $\C_\mathcal{A}$ that violated $\log[\phi_\lambda(\C_\mathcal{A} \, | \, \beta_\mathcal{A})] \leq \log[\phi_\lambda(\bar{\C}_\mathcal{A} \, | \, \beta_\mathcal{A})]$. All such settings had a $\beta < 0.628$ and $\lambda < 0.028$, similar to the scenarios that violated log concavity for $k=2$. While these 87 settings were not log concave, all of the maximum observed $\C_\mathcal{A}$ were improved upon by a $\C_\mathcal{A} \in \mathcal{C}_3$ with off-diagonal $c^*$ nearly equal to the lower bound of $-0.5$. Therefore, even in the observed cases where log concavity does not hold, it is reasonable to restrict the search for an optimal correlation matrix to $\mathcal{C}_3$.

Given what we observed for the $k=2$ case, we did not anticipate  $\phi_\lambda^\pm(\C_\mathcal{A} \, | \, \betavec_\mathcal{A})$ to be log concave even for completely symmetric $\C$, except for certain combinations of $\beta$ and $\lambda$. Our main concern though is whether there exists an optimal $\C_\mathcal{A}$ that is completely symmetric. We performed the same numerical study as for the known sign case but checked whether $\log[\phi_\lambda^\pm(\C_\mathcal{A} \, | \, \betavec_\mathcal{A})]\leq \log[\phi_\lambda^\pm(\bar{\C}_\mathcal{A} \, | \, \betavec_\mathcal{A})]$ and, if not, whether the best observed $\C_\mathcal{A}$ could be improved upon by some $\C_\mathcal{A}^* \in \mathcal{C}_3$. Only 65 of the 750 settings violated the inequality with all but one case having $\lambda > 1.3$. Of these 65 settings, $\C_\mathcal{A}^* \in \mathcal{C}_3$ improved the criterion value over the best observed $\C_\mathcal{A}$ in all but one case. This one case had $\beta=2.8651$ and $\lambda=2.7535$ and the optimal observed off-diagonals were $(c_{12},c_{13},c_{23})=(-0.0053, -0.0484, -0.5832)$ with $\phi_\lambda^\pm(\C_\mathcal{A} \, | \, \betavec_\mathcal{A})=0.1729$; the criterion value for $\C_\mathcal{A}^* \in \mathcal{C}_3$ was $0.1686$. This numerical study provides evidence that for a fixed $\lambda$, the optimal correlation matrix may not be completely symmetric, but generally there exists such a matrix that can perform nearly as well or better.


We also evaluated log concavity of $\phi_\lambda^\pm(\C_\mathcal{A} \, | \, \betavec_\mathcal{A})$ across the 750 settings restricted to $\C_\mathcal{A} \in \mathcal{C}_3$, calculating the second-order derivatives of $\log[\phi_\lambda^\pm(\C_\mathcal{A} \, | \, \betavec_\mathcal{A})]$ with second-order finite differences with a step size of $h=0.01$. Of the 750 settings, we found 77 that violated the log concavity condition but all scenarios gave a unique $c^*$. Unlike with $\phi_\lambda(\C_\mathcal{A} \, | \, \betavec_\mathcal{A})$, only three of theses 77 settings had $\lambda < 0.4$, each with an optimal $c^* = 0$ that failed to be log concave at $c$ close to 1. The remaining 74 settings had $\lambda > 0.8$ and either had $c^*\in (0.5, 0.8)$ (with a maximum second-order derivative at $c \in (0,0.241)$) or $c^*=0$ (with a maximum second-order derivative at $c \in (0.190,0.811))$. When $c^* >0$, the maximal probability was essentially driven by only the probability of recovering $\zvec=\onevec$, with all other sign recovery probabilities essentially $0$. This is a consequence of averaging over sign vectors rather than focusing on the minimum probability across the sign vectors.

\subsubsection{Examples with $k=8$}

Maximizing $\phi_\lambda(\C_\mathcal{A} \, | \, \betavec_\mathcal{A})$ when $k=8$ involves an optimization problem across 28 correlation parameters. We focus our numerical investigation to $\beta=4(=\sqrt{16}\times 1)$ and $\beta=12(=\sqrt{16} \times  3)$, which correspond to the last two panel rows in Figure~1 in Section~2.1. These two values are often encountered in simulation studies performed in SSD literature. For each $\beta$ we checked the log concavity inequality $\log[\phi_\lambda(\C_\mathcal{A} \, | \, \betavec_\mathcal{A})] \leq \log[\phi_\lambda(\bar{\C}_\mathcal{A} \, | \, \betavec_\mathcal{A})]$ across 1000 randomly generated positive definite correlation matrices for $\log(\lambda)$ between $-3$ to $3$ in steps of $0.2$. We also checked whether the maximum observed probability could be improved upon by $\C_\mathcal{A}^* \in \mathcal{C}_8$. Finally, we repeated the numerical investigation for $\phi_\lambda^\pm(\C_\mathcal{A} \, | \, \betavec_\mathcal{A})$ but for only 500 randomly generated positive definite correlation matrices due to the increased computations from averaging over $2^7$ sign vectors.

For $\beta=4$, the inequality constraint was violated for 29 out of the 31 possible $\lambda$ values. For $\log(\lambda) \leq 0.2$, no $\C_\mathcal{A} \in \mathcal{C}_8$ had a higher $\phi_\lambda(\C_\mathcal{A} \, | \, \betavec_\mathcal{A})$ than the maximum observed probability, but the differences in probabilities were negligible, with both close to 1. For $\log(\lambda)>0.2$, we were able to find a $\C_\mathcal{A} \in \mathcal{C}_8$ with a higher $\phi_\lambda(\C \, | \, \betavec)$ than the observed maximum, with meaningful differences in the two probabilities. A notable instance was for $\log(\lambda)=2$ where the maximum observed probability was $0.1282$ and the maximum probability across $\mathcal{C}_8$ was $0.7141$. Similar behavior was observed for $\beta=12$, except we found a $\C_\mathcal{A} \in \mathcal{C}_8$ that exceeded the maximum observed probability for all $\lambda$. 

Our investigation into $\phi_\lambda^\pm(\C_\mathcal{A} \, | \, \betavec_\mathcal{A})$ revealed similar results as for $k=3$, being that the log concavity condition is violated for larger $\lambda$.   For $\beta=4$, this occurred for $\log(\lambda) \in \{1.6,1.8,...,2.4\}$ but we were always able to find a $\C_\mathcal{A}^* \in \mathcal{C}_8$ that had larger $\phi_\lambda^\pm(\C_\mathcal{A} \, | \, \betavec_\mathcal{A})$. For these cases, the off-diagonal element for the $\C_\mathcal{A}^*$'s were $c^* \geq 0.4172$, giving the designs a slight edge in recovering $\zvec_\mathcal{A}=\onevec$. Their $\phi_\lambda^\pm(\C_\mathcal{A} \, | \, \betavec_\mathcal{A})$ values, however, were still close to 0, due to the fact that $\zvec_\mathcal{A}=\onevec$ is only one of $128$ possible sign vectors. We observed similar results for $\beta=12$ for $\log(\lambda)=\{2.6,2.8,3\}$.

\subsection{Optimal $\C \in \mathcal{C}_{p,k}$ for relaxed criteria}


The numerical results in Section~3.1 provide a reasonable justification that the optimal $\C^*$ (or at least a nearly optimal $\C^*$) for approximate relaxed criteria can be found in $\mathcal{C}_{p,k}$.  
Therefore, this section focuses on identifying such a $\C^*$. This optimization can be done numerically and we investigate the localized derivatives of the criteria  at $c=0$ for fixed $\lambda$ values. We argue $\C^*=\I$ for summaries of $\Phi_{\lambda}^{\pm}(\C \, | \, k, \beta)$, and that there exists a $\C^*$ with all positive off-diagonals for summaries of the sign-dependent $\Phi_{\lambda}(\C \, | \, k, \beta)$. 



Events \eqref{eq:SignCheckUnscale} and \eqref{eq:InactiveCheck2} will be denoted here by $P(S_\lambda \, | \, c, \betavec_\mathcal{A})$ and $P(I_\lambda \, | \, c, \zvec)$.
The probabilities for the two events for $\C \in \mathcal{C}_{p,k}$ are given in the following lemma:
\begin{lemma}\label{lem:CS_events} Let $\V=\I$ and $\C \in \mathcal{C}_{p,k}$. Then for $\mathcal{A}$ where $|\mathcal{A}|=k$ and $\betavec_\mathcal{A}$ with sign vector $\zvec_\mathcal{A}$, $P(S_\lambda \, | \, c, \betavec_\mathcal{A})= P( \boldsymbol{u}< |\betavec_\mathcal{A}|)$ and $P(I_\lambda \, | \, c, \zvec)= P(|\boldsymbol{v}| \leq \lambda \onevec)$
where 
\begin{align*}
\boldsymbol{u}&\sim 
    N\left(\frac{\lambda}{1-c}\left[\onevec - z_\mathcal{A}\gamma\zvec_\mathcal{A} \right],\ \frac{1}{1-c} \left[\I_{k} - \gamma\zvec_\mathcal{A}\zvec_\mathcal{A}^T\right] \right)\ ,\    \\
    \boldsymbol{v}&\sim N\left(\lambda z_\mathcal{A}\gamma\onevec,(1-c)\left[\I + \gamma\J\right] \right)\ ,\
\end{align*}
 with $z_\mathcal{A}=\onevec^T\zvec_\mathcal{A}$ and $\gamma = c/(1+c(k-1))$. Moreover, if $\betavec_\mathcal{A}$ does not depend on $\mathcal{A}$, the probabilities of the two events are constant across all such $\mathcal{A}$.
\end{lemma} 
\noindent

\noindent
When $\betavec_\mathcal{A} = \beta \onevec$ the random vectors $\uvec$ and $\vvec$ are distributed as
\begin{align}
    \boldsymbol{u}\sim 
    N\left(\frac{\lambda}{1+c(k-1)}\onevec,\frac{1}{1-c} \left[\I - \gamma\J\right] \right), \qquad 
    \boldsymbol{v}\sim N\left(\lambda k\gamma\onevec,(1-c)\left[\I + \gamma\J\right] \right)\nonumber\ ,\
\end{align}
for all $\mathcal{A} \in \mathcal{A}_k$. The corresponding criterion defined with respect to $\mathcal{C}_{p,k}$ and for a fixed $\lambda$ is denoted
\begin{align}
\psi_\lambda(c \, | \, k, \beta)=P(S_\lambda \, | \, c, \betavec_\mathcal{A}=\beta\onevec) \times P(I_\lambda \, | \, c, \zvec) \label{eqn:CS_sign_crit}\ .\
\end{align}
The analog to $\Phi_{\lambda}^\pm(\X \, | \, k, \beta)$ is denoted
\begin{align}
\psi_\lambda^{\pm}(c \, | \, k, \beta)=\frac{1}{2^{k-1}}\sum_{\tilde{\zvec} \in \mathcal{Z}_\mathcal{A}^{\pm}} P(S_\lambda \, | \, c, \betavec_\mathcal{A}=\beta\tilde{\zvec}) \times P(I_\lambda \, | \, c, \zvec=\tilde{\zvec}) \label{eqn:CS_crit}\ .\
\end{align}
The summarized versions of \eqref{eqn:CS_sign_crit} and \eqref{eqn:CS_crit} across $\lambda$ will be denoted similarly to their exact design counterparts from Section~\ref{sec:loc}, replacing $\Phi$ with $\psi$.

\begin{figure}[ht]
    \centering
    \includegraphics[width=0.75\textwidth]{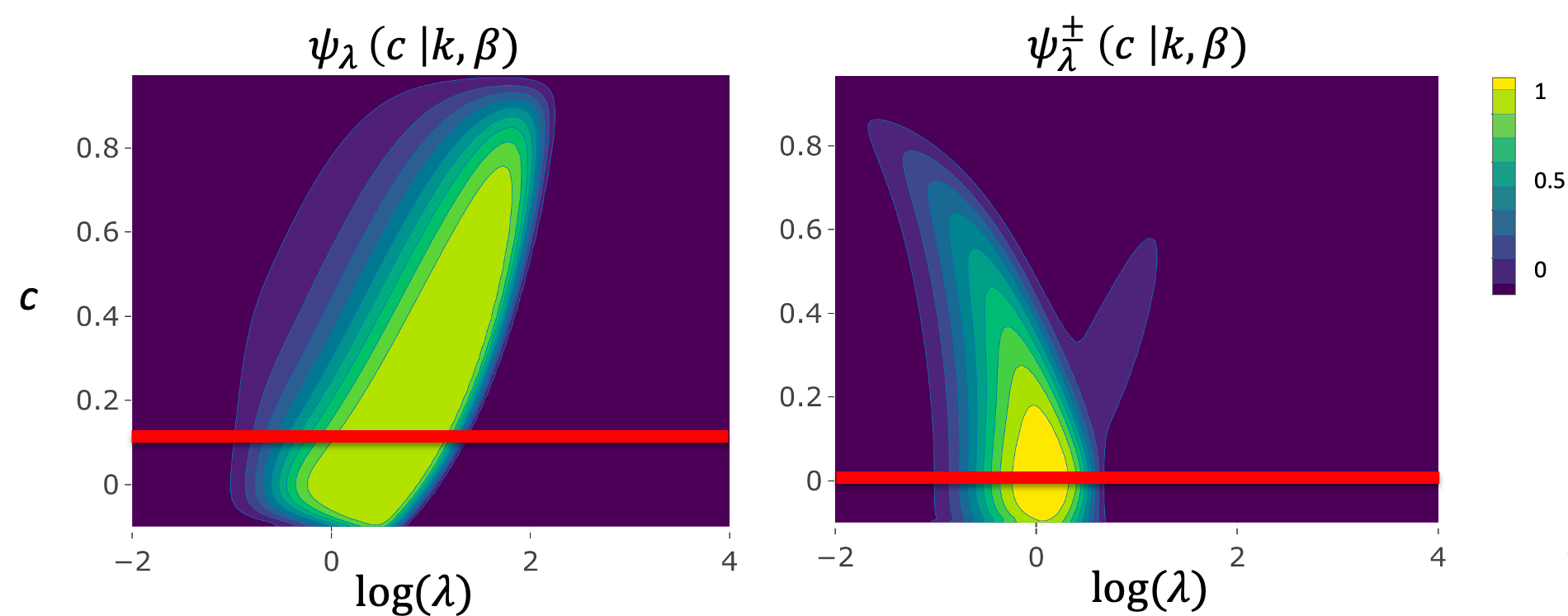}
    \caption{Contour plots of sign recovery probabilities  when $\C$ is completely symmetric across possible values of the off-diagonal $c$ and $\log(\lambda)$ where $k=4$ factors are active with $\beta=2\sqrt{10}$. The red lines correspond to values of $c=0.14$ for $\psi_\Lambda$ and $c=0$ for $\psi_\Lambda^\pm$. }
    \label{fig:CS_combined}
\end{figure}

The criteria $\psi_\lambda$ and $\psi_\lambda^{\pm}$, as well as their summarized versions, involve a single variable but are still challenging to optimize analytically.  Fortunately, numerical optimization of these criteria is straightforward and computationally efficient. We demonstrate numerical optimization of these new criteria for the situation of $p=10$, $k=4$, and $\beta=2\sqrt{10}$.   
Figure~\ref{fig:CS_combined} shows the contour plots of $\psi_\lambda$ and $\psi_\lambda^{\pm}$. The optimal $c$ values for the corresponding integral summary measures, $\psi_\Lambda$ and $\psi_\Lambda^\pm$, are $0.14$ and $0$, respectively. The resulting optimal $\C$ matrices then match the ideal forms for the $\text{Var}(s+)$ and $\text{UE}(s^2)$-criterion.  While these ideal $\C$ are not possible for SSDs, this example provides further justification for these two heuristic SSD criteria for sign recovery under the lasso. 

If $\psi^\pm_\lambda$ were concave in $\mathcal{C}$ across all $\lambda$, the matrix averaging technique would prove the optimality of $c=0$, as $\C=\I$ is the unique matrix after averaging across all permutations and sign transformations. Figure~\ref{fig:CS_combined} shows that there are some $\lambda$ for which $c=0$ does not maximize $\psi^\pm_\lambda$. For example $\log(\lambda) > 0.640$ shows a nearly 0 value for $\psi^\pm_\lambda$ when $c=0$ but a nonzero probability for $c \in (0.2, 0.6)$. This phenomenon is due to the improved sign recovery for $\zvec_\mathcal{A}=\onevec$ for these $c$ values in this range of $\log(\lambda)$, which was also observed and noted in Section~3.1.2. 

We now study the behavior of $\psi_\lambda$ and $\psi^\pm_\lambda$ in a neighborhood about $c=0$ to address the potential optimality of orthogonal designs. 
Direct application of the multivariate Leibniz integral rule gives the following lemma.
\begin{lemma}\label{lem:CS_event_derivs}
For all $\mathcal{A}$ where $|\mathcal{A}|=k$ and $\betavec_\mathcal{A}=\beta \onevec > \zerovec$, $\frac{d}{dc} P(I_\lambda | c, \zvec)\big|_{c=0} = 0$ for all $\lambda$ and $\frac{d}{dc} P(S_\lambda | c, \betavec_\mathcal{A}=\beta\onevec)\big|_{c=0} \geq 0$ for $\lambda$ satisfying
\begin{align}
2\lambda &\geq \frac{g(\tau)}{G(\tau)}\label{eqn:Lem3_ineq1}\ ,\ 
\end{align}
where $g(\cdot)$ and $G(\cdot)$ refer to the standard Normal probability density function and cumulative distribution function, respectively, and $\tau = \beta -\lambda$. 
\end{lemma}
\noindent
A direct consequence of Lemma~\ref{lem:CS_event_derivs} is that $\psi_\lambda(0| k, \beta)$ can often be improved upon by some $c > 0$. We state this important result as a theorem.
\begin{theorem}\label{thm:CS_local_opt_sign}
For all $\mathcal{A}$ where $|\mathcal{A}|=k$, $\frac{d}{dc} \psi_\lambda(c |k, \beta)\big|_{c=0}>0$ for $\lambda$ satisfying~\eqref{eqn:Lem3_ineq1}. Hence there exists some $c_\lambda>0$, where $\psi_\lambda(c_\lambda | k, \beta) > \psi_\lambda(0 | k, \beta)$.
\end{theorem}

\noindent
Applying this result to $\psi_\lambda(c\, | \, k,\beta)$ shown in Figure~\ref{fig:CS_combined}, inequality \eqref{eqn:Lem3_ineq1} holds for $\log(\lambda)\geq-22.763$, which covers the entire region of $\lambda$ values for which $\psi_\lambda(c\, | \, k,\beta)>0$. It follows then that some $c > 0$ will maximize the summary measures $\psi_\Lambda(c\, | \, k,\beta)$ or $\psi_{\max}(c\, | \, k,\beta)$. Theorem~\ref{thm:CS_local_opt_sign} thus provides a mathematically rigorous justification that the ideal designs sought after by the $\text{Var}(s+)$-criterion maximize the probability of sign recovery for $\betavec_\mathcal{A}=\beta \onevec$. The investigations in Section~\ref{sec:OptimalCA} further suggest such designs may be superior when $\zvec_\mathcal{A}=\onevec$, so long as the elements of $\betavec_\mathcal{A}$ are sufficiently large. 

The following theorem establishes a similar justification for heuristic orthogonality measures under unknown signs.
\begin{theorem}\label{thm:CS_locl_opt_allsign}
For all $\mathcal{A}$ where $|\mathcal{A}|=k$ and $\beta>0$, $\frac{d}{dc} \psi_\lambda^{\pm}(c\, | \, k,\beta)\big|_{c=0}=0$ for all $\lambda$. Moreover, $c=0$ is a local maximum for $\psi_\lambda^\pm(c\, | \, k,\beta)$ when
\begin{equation}
\begin{split}
\frac{q}{\binom{k}{2}} \ \frac{\lambda g(\lambda)}{1-2G(-\lambda)}&\left( k(1-\lambda^2)+(q-1)\frac{\lambda g(\lambda)}{1-2G(-\lambda)}\right) \leq\\
&\frac{g(\tau)}{G(\tau)}\left(\beta+\lambda+\lambda^2\tau-\left[\frac{\beta^2-\lambda^2}{2}+\beta\lambda\right]\frac{g(\tau)}{G(\tau)}\right)\label{eqn:Thm3_ineq}
\end{split}
\end{equation}
where $q=p-k$. 
\end{theorem}
\noindent
Applying this result to $\psi_\lambda^\pm(c\, | \, k,\beta)$ shown in Figure~\ref{fig:CS_combined}, inequality \eqref{eqn:Thm3_ineq} holds for $\log(\lambda)\in [-0.988,0.640]$. For $\log(\lambda)$ outside this region, there are some $c>0$ for which $\psi_\lambda^{\pm}(c)>0$, although the probabilities are relatively small. These small probabilities do not influence $\psi_{\max}^\pm$ and have minimal influence on $\psi_\Lambda^\pm$, making $c=0$ a likely global maximum for both criteria. We conjecture that generally $c^*=0$ and some $c^*>0$ are global maxima for $\psi_\lambda^\pm(c\, | \, k,\beta)$ and $\psi_\lambda(c\, | \, k,\beta)$, respectively.

\section{Fast Design Construction with HILS}\label{sec:EvalConstruct}

Section~\ref{sec:optimal_cs} identified optimal forms of completely symmetric $\C$ matrices under different assumptions about $\zvec$. Unfortunately, no SSD exists whose $\C$ achieves these forms because $n < p +1$. Ideally, one would implement a design search algorithm that ranks SSDs according to a summary criterion of the $\Phi_\lambda$- or $\Phi_\lambda^{\pm}$-criterion. However, these criteria demand intense computations to evaluate a single SSD and search algorithms require many evaluations. We now describe an algorithmic construction that compromises the rigorous but computationally-intensive criteria in Section~\ref{sec:loc} with new computationally-efficient heuristic criteria based on the theory in Section~\ref{sec:optimal_cs}.

The core approach to Section~\ref{sec:optimal_cs} was performing a joint optimization with respect to $\V$ and $\C$. The optimal form of $\V=\I$ holds for all scenarios, while the optimal form of $\C$ is assumed to be completely symmetric with off-diagonal element, say $c^*$. If $\zvec$ is unknown, we take $c^*=0$ and for known $\zvec=\onevec$, $c^*>0$ with the exact value depending on $\beta$ and $k$.  The value of $c^*$ is straightforward to calculate with numerical optimization methods. Let $\C^*$ denote the optimal $\C$. We propose first generating a candidate set of designs that simultaneously optimize two simple heuristic criteria based on Section~\ref{sec:optimal_cs}: 
\[
\text{tr}(\V) \quad \text{   and   } \quad ||\C-\C^*||_F\ ,\
\]
where $||\cdot||_F$ denotes the Frobenius norm. We want $ncand$ candidate designs that both maximize $\text{tr}(\V)$ and minimize $||\C-\C^*||_F$. These $ncand$ designs will then be ranked according to one of the lasso criteria. 

This approach differs from traditional techniques to construct SSDs, which optimize a single criterion measuring the size of the off-diagonals of $\X'\X$. The $\text{E}(s^2)$-criterion is a special case of this approach, as it first requires balanced columns, thus maximizing $\text{tr}(\V)$ and leads to $\C=\X^T\X$. For such designs, minimizing the sum of squared off-diagonals of $\X^T\X$ is equivalent to minimizing $||\C-\I||_F$. Another important contribution of this construction framework is that the optimal design will also come with an interpretable sign recovery probability value that better reflects a design's statistical value than simply comparing designs based on the off-diagonals of their $\X^T\X$ matrices. We demonstrate this new assessment for the examples in Section~5.

To construct designs based on multiple criteria, \citet{Caoetal2017} combined the coordinate exchange algorithm \citep{meyer_nachtsheim1995} with an elitist-like operator to produce $ncand$ designs lying near the Pareto front of the criteria. A set of $ntry$ initial designs are provided and each will generate its own Pareto set. That is, for each initial design, many exploratory designs are generated via coordinate exchanges. Only those exploratory designs that either dominate the initial design, or are non-dominated by the other exploratory designs form the initial design's Pareto set. The procedure continues until each non-dominated design cannot be further improved via coordinate exchange. The approximate Pareto fronts from each initial design are then combined into a single front comprised of the non-dominated designs across all $ntry$ Pareto sets, giving a set of, say, $ncand$ designs. 
Once the final approximate Pareto set is constructed, we identify which of its $ncand$ designs maximizes the lasso criterion of interest.


The proposed construction algorithm first identifies a set of $ncand$ candidate designs through the Pareto front optimization search algorithm from  \cite{Caoetal2017} based on the two lasso criteria. We then employ the exact lasso sign recovery criteria from  Section 2 as a secondary criterion, as it more accurately measures a design's screening capabilities, to rank the $ncand$ designs and choose the best one. We refer to this algorithm as the Heuristic-Initiated Lasso Sieve, or HILS.

\begin{figure}[ht]
    \centering
    \includegraphics[width=0.9\linewidth]{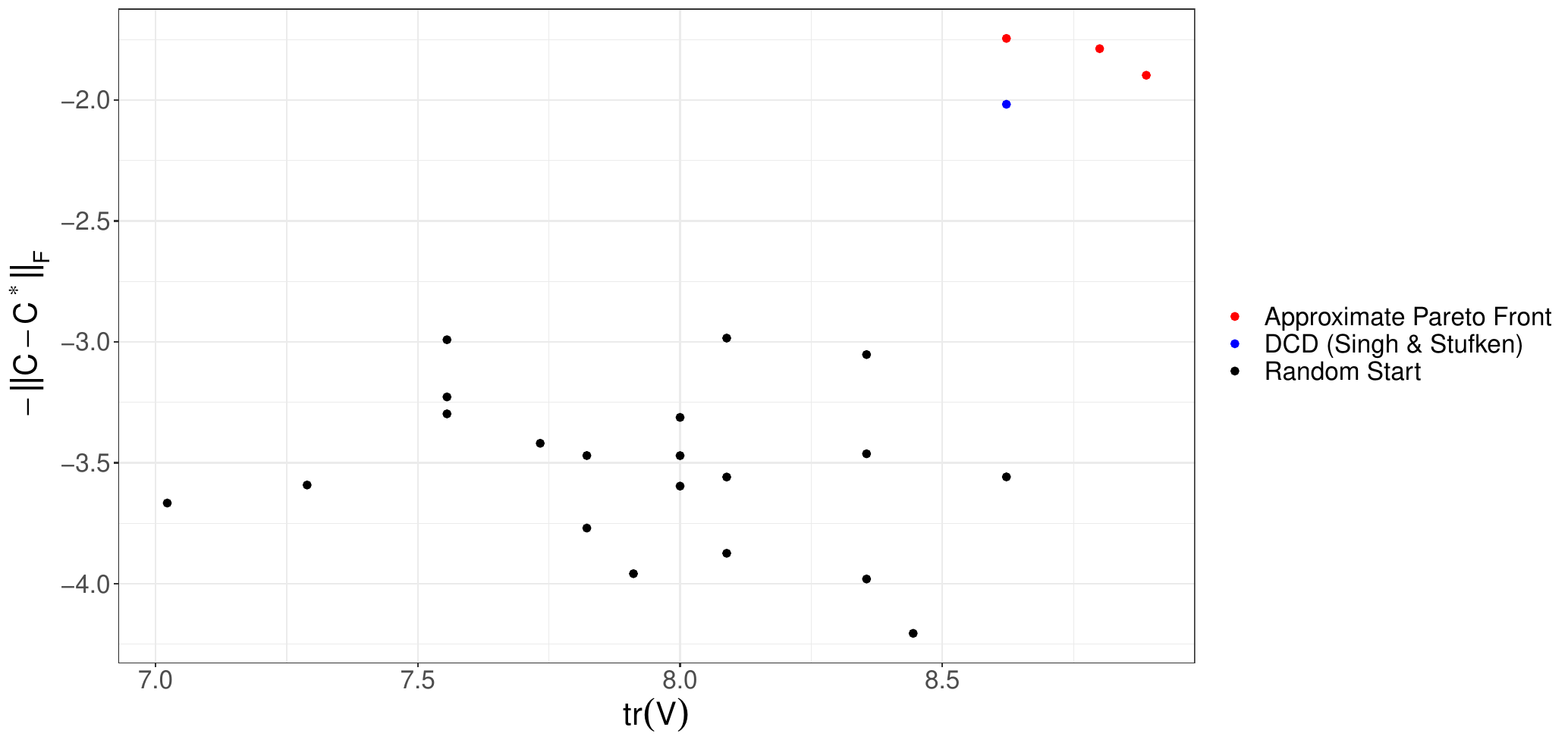}
    \caption{Plot of the approximate Pareto front for the three designs selected by the coordinate exchange algorithm for the $n=9$, $p=10$, unknown signs case. The starting criteria values for each of the $ntry=20$ starting designs from the algorithm are plotted in black. Each point represents a different design. The criteria values for the DCD are also plotted to show that the approximate Pareto Front designs the DCD in terms of the two criteria. }
    \label{fig:PF_example}
\end{figure}

To demonstrate, take $n=9$ and $p=10$ assuming unknown signs, making $c^{*}=0$. The Pareto front coordinate exchange algorithm with $ntry=20$ random starting designs concluded in 12 seconds on a Desktop with an Intel \textsuperscript{\textregistered} Core\textsuperscript{\texttrademark} i7-6700 CPU @ 3.40GHz processor with 16GB RAM. From this proceedure, three designs were found on the Pareto front maximizing $\text{tr}(\V)$ and minimizing $||\C -\I||_F$. Figure \ref{fig:PF_example} shows the Pareto front of the three designs. Note that, to present the Pareto front, we change to maximizing $-||\C -\I||_F$ instead of minimizing $||\C -\I||_F$. For comparison, the DCD from \cite{singh2023selection} is also plotted. Here, the DCD is dominated by the three designs selected by the coordinate exchange algorithm in terms of the $\text{tr}(\V)$ and $||\C-\I||_F$.

Although HILS can significantly reduce the number of evaluations of the relaxed lasso criteria compared to constructing an optimal design using, for instance, a coordinate exchange algorithm, there are settings when even a more modest number of evaluations is not computationally feasible. The same issue was encountered by \cite{singh2023selection}. We now discuss some ways to efficiently evaluate the $\Phi_\lambda$ and $\Phi_\lambda^{\pm}$-criteria, and hence their corresponding summarized versions. For a support size of $k$ and design $\X$, one full evaluation of $\Phi_{\lambda}$ and $\Phi_{\lambda}^{\pm}$ requires consideration of all $\binom{p}{k}$ supports of size $k$. For $\Phi_{\lambda}^\pm$, an additional $2^{k-1}$ computations are required for each $\mathcal{A}$ to consider all sign vectors. Evaluating either criteria would benefit from a reduction in the number of supports considered. Sampling from $\mathcal{A}_{k}$ completely at random is intuitive, but may require a large number of samples to be representative.  As high correlations between two factors can result in diminished lasso support recovery performance due to low $P(I_\lambda \, | \, \C, \, \zvec)$, a support sampling method that balances (or approximately balances) the number of times pairwise sets of factors are included in $\mathcal{A}$ together can be advantageous. \cite{smucker2015approximate} utilized nearly balanced incomplete block designs (NBIBDs) to approximate a model space with only a relatively small number of blocks, in the context of model-robust optimal designs. Hence, we recommend implementing the NBIBD sampling method to adequately represent $\mathcal{A}_{k}$ with between 64-128 supports for modest $p$ and $k$ values. 

From Lemma~\ref{lem:symmetry}, for a fixed $\mathcal{A}$, the probabilities are equal between $\zvec_\mathcal{A}$ and $-\zvec_{\mathcal{A}}$. Thus, only $2^{k-1}$ sign vectors need be considered for the $\Phi_\lambda^{\pm}$-criterion. While this cuts the computation in half, prior knowledge from the practitioner can be leveraged to select an even smaller set of representative sign vectors to evaluate. For example, if there is no knowledge on sign direction, the most likely sign vectors are those with an equal, or nearly equal, number of $\pm 1$. For $k=10$, one could reduce the number of sign vectors from $2^9=512$ to only the $126$ possible sign vectors having five $+1$'s and $-1$'s. There may also be strong prior belief about the signs of some or all of the factors that can reduce the set of sign vectors to a manageable number, or even to a single element. 

\section{Examples}\label{sec:NewDesigns}

We now demonstrate near-optimal SSDs found under the criteria via HILS under three example scenarios. For each case, the coordinate exchange algorithm was executed with $ntry=50$ random starting designs and identified $ncand$ final designs on the approximate Pareto front of the criteria described in Section \ref{sec:EvalConstruct}. Of these $ncand$ designs, the final HILS design was selected using the $\Phi_{\text{max}}^{\pm}$ criterion. For each scenario, we provide the computational time that it took to find the HILS-optimal design on a Desktop with an Intel \textsuperscript{\textregistered} Core\textsuperscript{\texttrademark} i7-6700 CPU @ 3.40GHz processor with 16GB RAM.  Additionally, we evaluate the performance of the HILS designs compared to competitors under $k$ and $\beta$ settings that are different to the $k$ and $\beta$ values used to construct the HILS designs. This evaluates the sensitivity of the HILS construction method to $k$ and $\beta$.

For a final comparison of competing designs, we recommend 
plotting the $\Phi_\lambda$ and $\Phi_\lambda^{\pm}$ curves as functions of $\log(\lambda)$. These  curves represent the average probability of sign recovery of the lasso solution path over a range of tuning parameter values. In comparing designs in this way, we eliminate the need for simulation studies that often require some tuning parameter selection procedure. The design with larger $\Phi_{\lambda}$ or $\Phi_{\lambda}^{\pm}$ over a larger range of $\log(\lambda)$ is considered to be the best performing. For context, we also compare these curves to that for the hypothetically optimal design with $\V=\I$ and completely symmetric $\C=\C^*$. 

\subsection{Scenario 1: $n=9$, $p=10$ Unknown Signs}

The HILS algorithm was performed for $k=3$ and $\beta=3$ assuming unknown signs. The HILS-optimal design was selected in approximately six minutes of computation time. Because $n \approx p$, this scenario is close to one in which an orthogonal design exists. Therefore, this represents an experiment in which optimal designs under traditional SSD criteria should perform rather well. For unknown signs, we set $c^*=0$ following Theorem \ref{thm:CS_locl_opt_allsign}. Since $\C^*$ is orthogonal and the $\text{UE}(s^2)$ criterion targets orthogonality directly, the $\text{UE}(s^2)$-optimal designs were deemed appropriate for comparison. Additionally, since the two criteria used to construct the DCDs are optimized under an orthogonal design, the DCD was also compared.  

\begin{figure}[ht]
    \centering
    \includegraphics[width=0.9\linewidth]{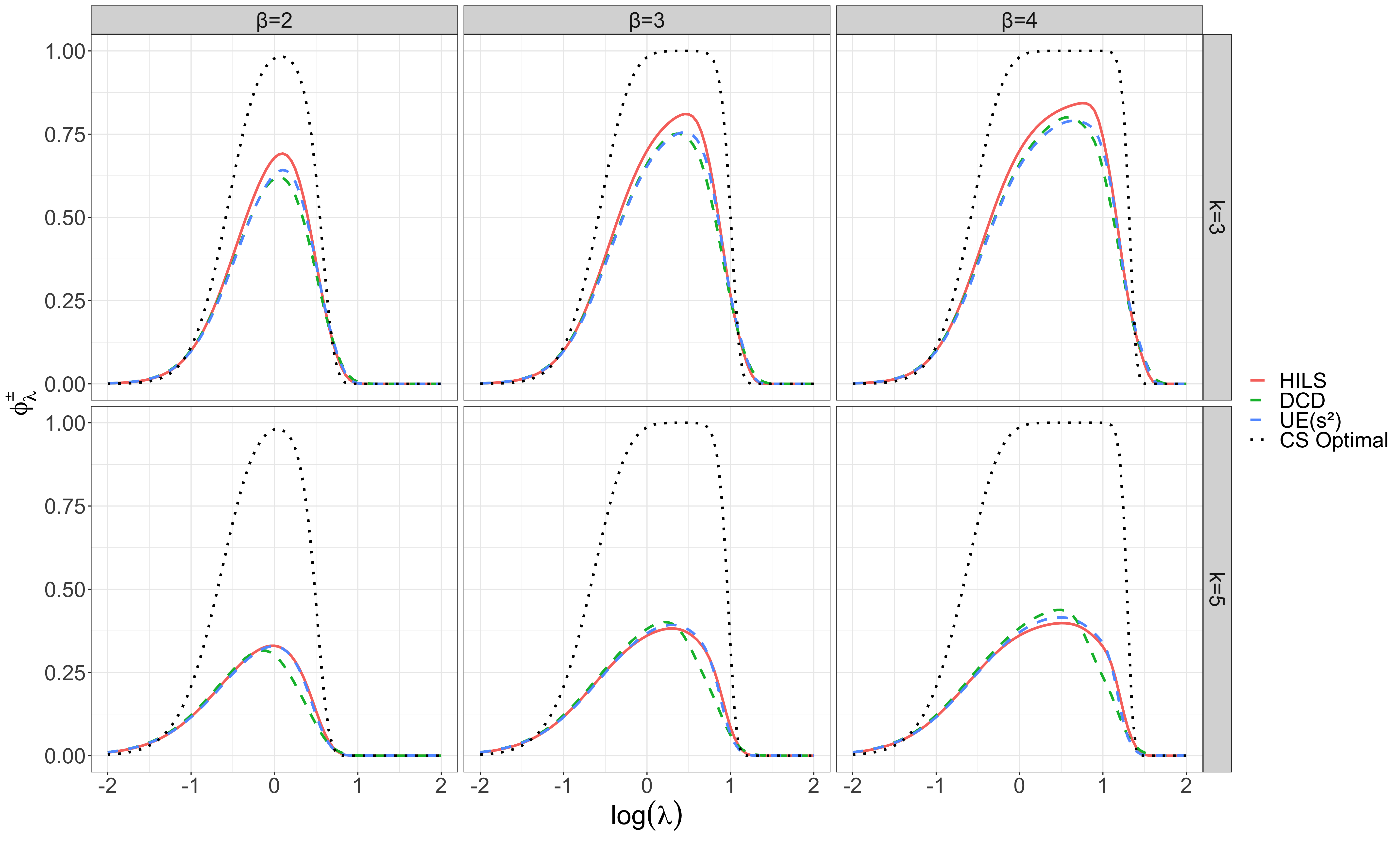}
    \caption{$\Phi_{\lambda}^{\pm}$ curves as a function of $\log(\lambda)$ for the HILS optimal design, DCD, and a $\text{UE}(s^2)$-optimal design under a variety of $k$ and $\beta$ settings for scenario 1. The HILS optimal design shown in all panels was chosen assuming $k=3$ and $\beta=3$. The CS Optimal curve refers to the hypothetical design having the optimal completely symmetric form for $k=3$ and $\beta=3$.}
    \label{fig:Scenario1AllSigns_sens}
\end{figure}

Figure \ref{fig:Scenario1AllSigns_sens} compares the $\Phi_{\lambda}^{\pm}$ curves as a function of $\log(\lambda)$ for the hypothetically optimal completely symmetric (CS optimal) design, HILS-optimal design, DCD, and a $\text{UE}(s^2)$-optimal design under a variety of $k$ and $\beta$ settings. The top-middle plot represents the $k$ and $\beta$ setting in which the HILS design was optimized for. The top row of Figure \ref{fig:Scenario1AllSigns_sens}, representing $k=3$, demonstrates that the HILS design out performs the others in terms of lasso sign recovery regardless of the $\beta$ value. When $k=5$, Figure \ref{fig:Scenario1AllSigns_sens} shows little difference between the 3 designs. The HILS and $\text{UE}(s^2)$- optimal designs show a larger $\Phi_{\lambda}^{\pm}$ for larger $\log(\lambda)$, but the DCD attains a slightly larger maximum of $\Phi_{\lambda}^{\pm}$ (for $\beta=3$ and $\beta=4$). Overall, the HILS-optimal design saw better performance than its competitors under the $k$ value that it was optimized for and comparable performance as the $k$ assumption was violated.


\subsection{Scenario 2: $n=9$, $p=10$ Known Signs}

To build on the $n$ and $p$ settings in scenario 1, the HILS algorithm was then performed for $k=3$ and $\beta=3$ assuming known signs. Here, the optimal $c$ value was found via numerical optimization to be taken as $c^*\approx 0.04$. The HILS algorithm selected the HILS-optimal design in approximately 4 minutes. 

\begin{figure}[ht]
    \centering
    \includegraphics[width=0.9\linewidth]{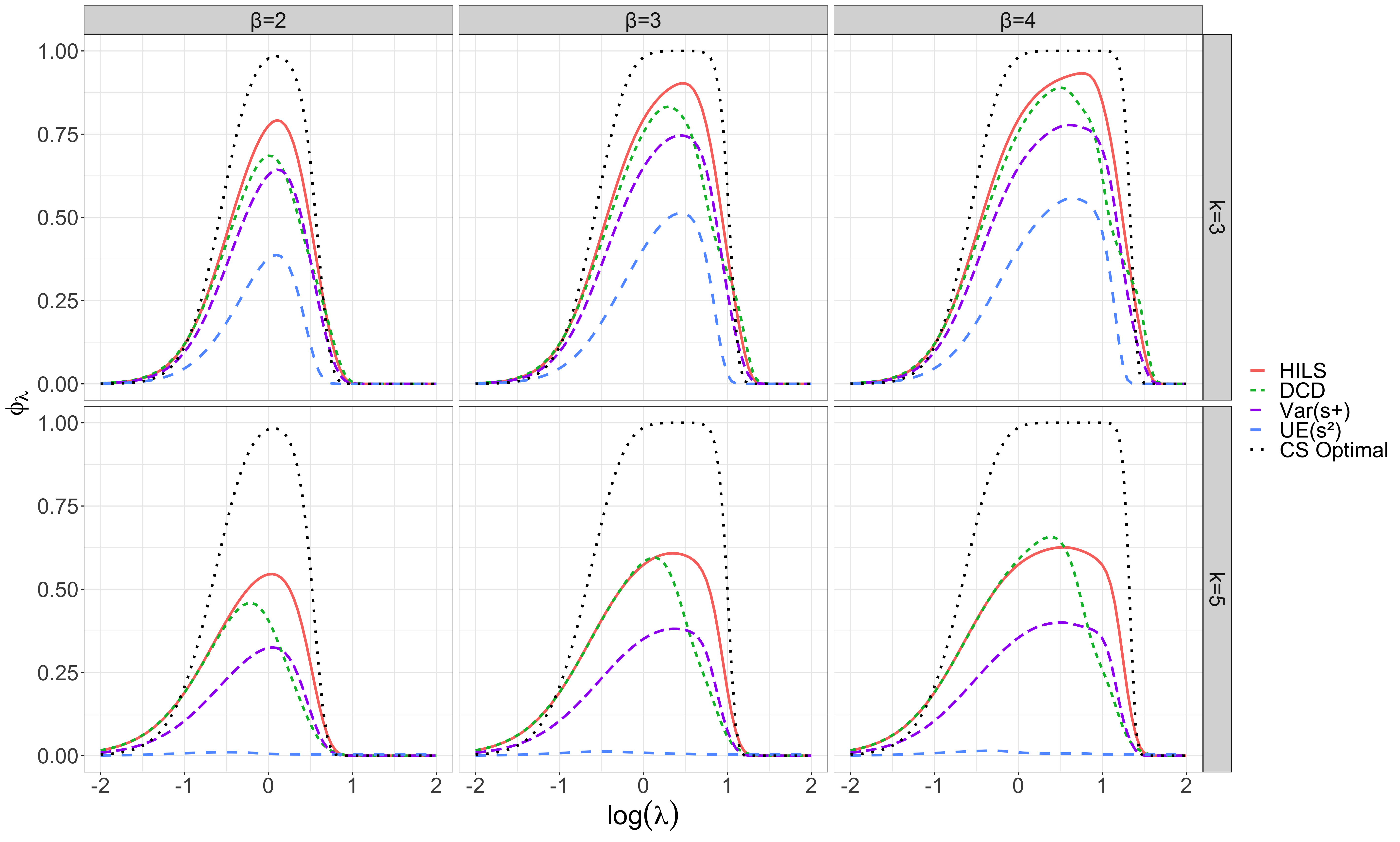}
    \caption{$\Phi_{\lambda}$ curves as a function of $\log(\lambda)$ for the HILS-optimal design, DCD, and a $\text{Var}(s+)$-optimal design under a variety of $k$ and $\beta$ settings under scenario 2. The HILS optimal design shown in all panels was chosen assuming $k=3$ and $\beta=3$. The CS Optimal curve refers to the hypothetical design having the optimal completely symmetric form for $k=3$ and $\beta=3$.}
    \label{fig:Scenario1AllPos_sens}
\end{figure}

    Figure \ref{fig:Scenario1AllPos_sens} compares the $\Phi_{\lambda}$ curves as a function of $\log(\lambda)$ for the CS-optimal design, HILS-optimal design, DCD \citep{singh2023selection}, and a $\text{Var}(s+)$-optimal design under a variety of $k$ and $\beta$ settings. To demonstrate the added value these three designs have by incorporating sign information, the figure also has the $\Phi_{\lambda}$ curve for a $\text{UE}(s^2)$-optimal design. The top-middle plot represents the $k$ and $\beta$ setting in which the HILS design was optimized for. For all $\beta$ and $k$ settings except $\beta=4$ and $k=5$, the HILS-optimal design shows a larger $\Phi_{\text{max}}$ and area under the curve compared to DCD and $\text{Var}(s+)$. For $\beta=4$ and $k=5$, Figure \ref{fig:Scenario1AllPos_sens} shows a slightly larger $\Phi_{\text{max}}$ for DCD compared to the HILS-optimal design. However, the HILS-optimal design sees higher $\Phi_\lambda$ for larger values of $\log(\lambda)$ (around $\log(\lambda)=1$) and has a larger area under its $\Phi_\lambda$ curve. Thus, we deem that the HILS-optimal design is superior to the competitors in all settings demonstrated in Figure \ref{fig:Scenario1AllPos_sens}. Finally, we note the poor performance of the $\text{UE}(s^2)$-optimal design compared to all designs, especially when $k=5$ with $\Phi_{\lambda}\approx 0$. This is most likely due to the prevalence of negative correlations leading to singular $C_\mathcal{A}$ for some supports.


\subsection{Scenario 3: $n=14$, $p=20$ Known Signs}

This scenario assumed known effect directions with $k=5$ and $\beta=2$. The optimal $c$ value was found via numerical optimization to be taken as $c^*\approx 0.045$. With a larger difference between $n$ and $p$, the Pareto front found by the coordinate exchange algorithm consisted of 14 designs (the Pareto fronts consisted of 3 designs for Scenarios 1 and 2). The HILS algorithm selected the HILS-optimal design in approximately 1 hour and 21 minutes.

\begin{figure}[ht]
    \centering
    \includegraphics[width=0.9\linewidth]{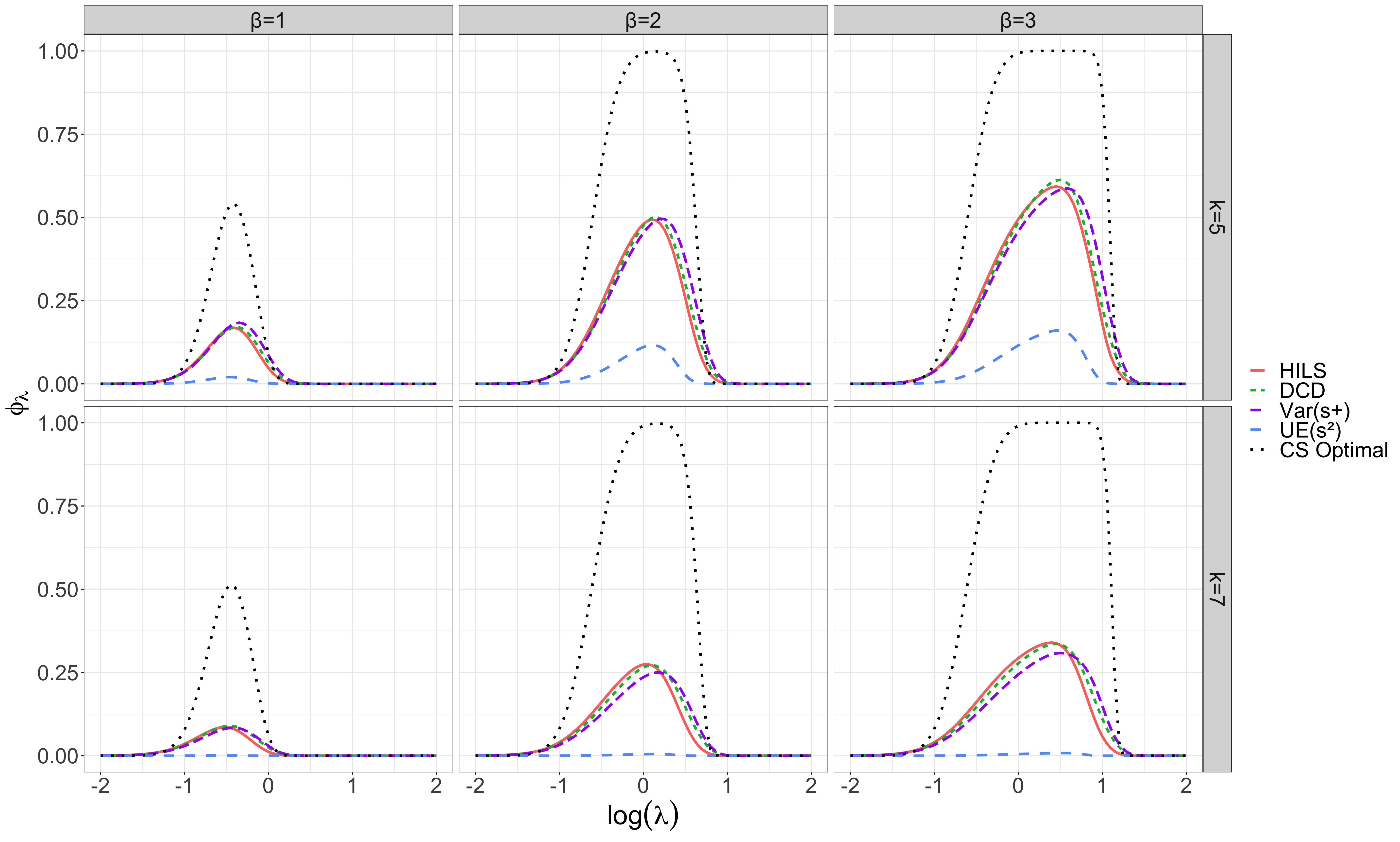}
    \caption{$\Phi_{\lambda}$ curves as a function of $\log(\lambda)$ for the HILS-optimal design, DCD, and a $\text{Var}(s+)$-optimal design under a variety of $k$ and $\beta$ settings under scenario 3. The HILS optimal design shown in all panels was chosen assuming $k=5$ and $\beta=2$. The CS Optimal curve refers to the hypothetical design having the optimal completely symmetric form for $k=5$ and $\beta=2$.}
    \label{fig:Scenario3AllPos_sens}
\end{figure}

Figure \ref{fig:Scenario3AllPos_sens} compares the $\Phi_{\lambda}$ curves as a function of $\log(\lambda)$ for the HILS-optimal design, DCD \citep{singh2023selection}, and a $\text{Var}(s+)$-optimal design under a variety of $k$ and $\beta$ settings. We again included the $\Phi_{\lambda}$ curve for a $\text{UE}(s^2)$-optimal design. The top-middle plot represents the $k$ and $\beta$ setting in which the HILS design was optimized for. Across all $\beta$ and $k$ settings there is very little difference between the performances of the three designs that include sign information.  However, all of these designs had significantly higher $\Phi_\lambda$ values than the $\text{UE}(s^2)$-optimal design.


The lack of separation between the HILS-optimal design from the other designs was not unexpected given $n$ is nearly half of $p$. It is challenging to find a design that simultaneously minimizes correlations. Nonetheless, the HILS-optimal design performed comparably to the DCD and $\text{Var}(s+)$-optimal designs, with all three significantly outperforming the $\text{UE}(s^2)$-optimal design. The HILS-optimal design was found much faster than the DCD, which took us dozens of hours to construct. Since the HILS-optimal design can be constructed more quickly and show comparable performance to the DCDs, we view this algorithmic design construction approach as worthwhile.

\section{Summary}\label{sec:Discussion}

The SSD literature has predominately constructed designs by optimizing heuristic criteria that measure a design's proximity to a (nonexistent) orthogonal design. The heuristic criteria are tractable in their optimization, in that optimal designs can generally be constructed directly or algorithmically in a reasonable amount of time. However, these heuristic criteria are not directly tied to a screening analysis method, so there is no guarantee that the resulting analysis under an optimal design will have good statistical properties. This article resolves this disconnect by optimizing criteria based on the probability of sign recovery under the lasso, which is well-defined even when $n-1 < p$. The major contributions are:
\begin{enumerate}
    \item \textbf{A local optimality criterion assuming known $\betavec$ and fixed $\lambda$}. A trivial design that confounds all inactive factors with the intercept is shown to be optimal for sign recovery. An exact design construction is given that can improve the probability of sign recovery over an orthogonal design for some $\lambda$ . The design has positive and nearly constant pairwise correlations, following the ideal structure for the $\text{Var}(s+)$-criterion.
    \item \textbf{More practical criteria that relax the assumptions about $\betavec$ and $\lambda$}. Such criteria are computationally intensive and hence difficult to optimize both analytically and algorithmically, requiring computations across all supports of a given size and potentially many sign vectors.
    \item \textbf{A study of the optimal form of $\C$ with and without known $\zvec$}. By conditioning on completely symmetric matrices, we arrive at a univariate optimization problem. The framework mimics the approximate design approach for least-squares analyses.
    \item \textbf{A rigorous proof that across all completely symmetric matrices, an orthogonal correlation matrix is suboptimal for nearly all $\lambda >0$ when $\zvec$ is known}. The optimal form instead takes on positive, constant off-diagonal elements. In the case of unknown sign information, an orthogonal correlation matrix is shown to be a local maximum for a range of $\lambda$. These optimal forms rigorously justify the $\text{UE}(s^2)$- and $\text{Var}(s+)$-criterion in the cases of unknown and known $\zvec$, respectively. This at least partially solves the two open problems from Section 1.1.
    \item \textbf{Heuristic-Initiated Lasso Sieve: a computationally efficient exact design construction algorithm that leverages new heuristic criteria motivated by the theory in Section~\ref{sec:optimal_cs}}. A Pareto front optimization is done with a coordinate exchange algorithm using two simple heuristic criteria targeting $\text{tr}(\V)$ and $||\C-\C^*||_F$. The proposed lasso criteria from Section~\ref{sec:loc} are then employed as a secondary ranking of many high-quality designs lying on the Pareto front. 
\end{enumerate}

The proposed optimal design framework provides an alternative to using simulation to compare SSDs, which can be tedious and difficult to reproduce independently. In particular, there are at least two reasons to be skeptical of conclusions drawn from simulations in the context of SSDs. First, simulations can be misleading when subtly different versions of complicated statistical procedures are used. In our attempts to reproduce and compare the simulation studies of \cite{singh2023selection} and \cite{weese2021strategies}, we discovered the regularization methods used were sensitive to two different aspects of the procedure implementation. In particular, these papers used the Gauss-Dantzig selector with different ways of exploring the tuning parameter space. They also differed in how they implemented a secondary thresholding of the estimates which determines a subset of the active factors that should remain active. Typically this threshold removes estimates whose absolute magnitude are less than $\sigma^2$, but thresholding in this way is of dubious reliability because there is no natural way to estimate $\sigma^{2}$ and thus threshold levels are more or less arbitrary.
The approach we present in this paper avoids these difficulties by eliminating the need for simulation at all, in a similar way that closed-form power analysis procedures are routinely used in simple, replicated experimental design settings.

Another danger of the existing heuristic criteria is that there is no guarantee the optimal design has any statistical value. For example, there technically exists a $\text{UE}(s^2)$-optimal design with $n=2$ runs and $p=100$ factors, but this design has no statistical value. The lasso criteria, however, do reflect statistical value and so can be used to give objective information about design quality, akin to a power analysis. Experimenters can specify a minimum effect size of interest, $\beta$, as well as an educated guess regarding the number of factors that are likely to have such an effect size and investigate the sign-recovery probabilities as a function of $\lambda$. If the maximum average probability is close to 0, it is unlikely that the design can provide reliable sign recovery for the specified effect size. Such a low value could lead to a reconsideration of the runs budget, expected sparsity, and/or the size of effects that are deemed of interest to detect. 

The proposed framework allows an even deeper investigation, if desired. For example, we can determine whether a low average probability is due to an inability to reliably exclude all inactive effects (the $I_\lambda$ event) or to reliably include all active effects (the $S_\lambda$ event). In screening experiments, we are much more concerned with identifying the true effects than with allowing inactive effects through the filter. This suggests that the $S_\lambda$ event is more important to the experimenter than the $I_\lambda$ event.

There are many avenues of future research to be explored. First, we noted in multiple sections and examples that for larger $\lambda$, $\Phi^\pm_\lambda(\C\, | \, \betavec)$ sometimes is optimized by a $\C$ with large positive off-diagonals, but its improvement is due solely to its ability to recover $\zvec_\mathcal{A}=\onevec$. This is practically undesirable because we would like to do well across many possible sign vectors. A similar criticism also holds for averaging across all $\mathcal{A} \in \mathcal{A}_k$. That is, the best design may do great for only some supports, and poorly for others. It would be interesting to consider such criteria that consider the minimum probability rather than the average. However, we expect that such criteria would have low probabilities that they would not seem practically useful, as it focuses on the worst possible scenario.

Another area of future work is to explore direct construction techniques for designs with all positive correlations, as such designs may be strong candidates to optimize the known-sign criteria. Continuing with exact design construction, while the HILS algorithm is faster relative to \cite{singh2023selection}, it is much slower than that for the heuristic criteria. More work is needed to make faster design construction algorithms. 

Finally, it is important to extend the criteria to incorporate more practical assumptions. Rather than assuming a known magnitude of the active effects, it would be better to incorporate uncertainty of the magnitudes using a prior distribution. This may require even heavier computation, as one would have to integrate the proposed criteria over the prior distribution. We are also considering criteria where the inactive factors have small, but nonzero, effects. This requires a definition of sign recovery with respect to the thresholded lasso in which estimated effects below some threshold value are set to 0.\\

\bigskip
\begin{center}
{\large\bf SUPPLEMENTARY MATERIAL}
\end{center}

Following the References are proofs of all results as well as computational details on the evaluation of our criteria. R files that construct designs following \eqref{eqn:LODbetter}; replicate the numerical studies in Section~3.1; identify optimal $c$ values for the criteria in Section~\ref{sec:optimal_cs}; and construct designs with the HILS algorithm may be found at
\url{https://github.com/hkyoung361/Lasso_Optimal_SSD}.












\bibliographystyle{asa} 
\bibliography{power.bib}

\begin{thebibliography}{59}
\newcommand{\enquote}[1]{``#1''}
\expandafter\ifx\csname natexlab\endcsname\relax\def\natexlab#1{#1}\fi

\bibitem[{Abraham et~al.(1999)Abraham, Chipman, and Vijayan}]{Abraham}
Abraham, B., Chipman, H., and Vijayan, K. (1999), \enquote{Some risks in the
  construction and analysis of supersaturated designs,} \textit{Technometrics},
  41, 135--141.

\bibitem[{Asif and Romberg(2010)}]{asif2010lasso}
Asif, M.~S. and Romberg, J. (2010), \enquote{On the lasso and Dantzig selector
  equivalence,} in \textit{2010 44th Annual Conference on Information Sciences
  and Systems (CISS)}, IEEE, pp. 1--6.

\bibitem[{Bickel et~al.(2009)Bickel, Ritov, and
  Tsybakov}]{bickel2009simultaneous}
Bickel, P.~J., Ritov, Y., and Tsybakov, A.~B. (2009), \enquote{Simultaneous
  analysis of lasso and Dantzig selector,} \textit{The Annals of statistics},
  37, 1705--1732.

\bibitem[{Booth and Cox(1962)}]{Booth62}
Booth, K. H.~V. and Cox, D.~R. (1962), \enquote{Some systematic supersaturated
  designs,} \textit{Technometrics}, 4, 489--495.

\bibitem[{Box and Hunter(1961{\natexlab{a}})}]{box1961}
Box, G.~E. and Hunter, J. (1961{\natexlab{a}}), \enquote{The {$2^{k-p}$}
  fractional factorial designs,} \textit{Technometrics}, 3, 311--351.

\bibitem[{Box and Hunter(1961{\natexlab{b}})}]{box19612}
--- (1961{\natexlab{b}}), \enquote{The {$2^{k-p}$} fractional factorial designs
  Part II.} \textit{Technometrics}, 3, 449--458.

\bibitem[{Candes and Tao(2007)}]{candes_tao2007}
Candes, E. and Tao, T. (2007), \enquote{The Dantzig selector: statistical
  estimations when $p$ is much larger than $n$,} \textit{The Annals of
  Statistics}, 35, 2313--2351.

\bibitem[{Cao et~al.(2017)Cao, Smucker, and Robinson}]{Caoetal2017}
Cao, Y., Smucker, B.~J., and Robinson, T.~J. (2017), \enquote{A hybrid elitist
  pareto-based coordinate exchange algorithm for constructing multi-criteria
  optimal experimental designs,} \textit{Statistics and Computing}, 27,
  423--437.

\bibitem[{Chen and Lin(1998)}]{CHEN199899}
Chen, J. and Lin, D.~K. (1998), \enquote{On the identifiability of a
  supersaturated design,} \textit{Journal of Statistical Planning and
  Inference}, 72, 99--107.

\bibitem[{Deng et~al.(1996)Deng, Lin, and Wang}]{deng1996measurement}
Deng, L.-Y., Lin, D. K.~J., and Wang, J. (1996), \enquote{A measurement of
  multi-factor orthogonality,} \textit{Statistics \& Probability Letters}, 28,
  203--209.

\bibitem[{Deng et~al.(1999)Deng, Lin, and Wang}]{deng1999resolution}
--- (1999), \enquote{A resolution rank criterion for supersaturated designs,}
  \textit{Statistica Sinica}, 9, 605--610.

\bibitem[{Deng et~al.(2013)Deng, Lin, and Qian}]{deng2013lasso}
Deng, X., Lin, C.~D., and Qian, P.~Z. (2013), \enquote{The lasso with nearly
  orthogonal Latin hypercube designs,} .

\bibitem[{Dragulji\'c et~al.(2014)Dragulji\'c, Woods, Dean, Lewis, and
  Vine}]{draguljic_etal2014}
Dragulji\'c, D., Woods, D.~C., Dean, A.~M., Lewis, S.~M., and Vine, A.-J.~E.
  (2014), \enquote{Screening strategies in the presence of interactions,}
  \textit{Technometrics}, 56, 1--16.

\bibitem[{Gai et~al.(2013)Gai, Zhu, and Lin}]{gai2013model}
Gai, Y., Zhu, L., and Lin, L. (2013), \enquote{Model selection consistency of
  Dantzig selector,} \textit{Statistica Sinica}, 615--634.

\bibitem[{Goos and Jones(2011)}]{goos2011optimal}
Goos, P. and Jones, B. (2011), \textit{Optimal Design of Experiments: A Case
  Study Approach}, John Wiley \& Sons.

\bibitem[{Hastie et~al.(2019)Hastie, Tibshirani, and
  Wainwright}]{hastie2019statistical}
Hastie, T., Tibshirani, R., and Wainwright, M. (2019), \textit{Statistical
  learning with sparsity: the lasso and generalizations}, Chapman and Hall/CRC.

\bibitem[{Huang et~al.(2020)Huang, Kong, and Ai}]{huang2020optimal}
Huang, Y., Kong, X., and Ai, M. (2020), \enquote{Optimal designs in sparse
  linear models,} \textit{Metrika}, 83, 255--273.

\bibitem[{James et~al.(2009)James, Radchenko, and Lv}]{james2009dasso}
James, G.~M., Radchenko, P., and Lv, J. (2009), \enquote{DASSO: connections
  between the Dantzig selector and lasso,} \textit{Journal of the Royal
  Statistical Society: Series B (Statistical Methodology)}, 71, 127--142.

\bibitem[{Javanmard and Montanari(2014)}]{javanmard2014confidence}
Javanmard, A. and Montanari, A. (2014), \enquote{Confidence intervals and
  hypothesis testing for high-dimensional regression,} .

\bibitem[{Jia and Rohe(2015)}]{jia2015preconditioning}
Jia, J. and Rohe, K. (2015), \enquote{Preconditioning the lasso for sign
  consistency,} \textit{Electronic Journal of Statistics}, 9, 1150--1172.

\bibitem[{Jones et~al.(2009)Jones, Li, Nachtsheim, and Ye}]{Jones2009}
Jones, B., Li, W., Nachtsheim, C.~J., and Ye, K.~Q. (2009), \enquote{{Model
  robust supersaturated and partially supersaturated designs},} \textit{Journal
  of Statistical Planning and Inference}, 139, 45--53.

\bibitem[{Jones and Majumdar(2014)}]{jones2014optimal}
Jones, B. and Majumdar, D. (2014), \enquote{Optimal supersaturated designs,}
  \textit{Journal of the American Statistical Association}, 109, 1592--1600.

\bibitem[{Khuri et~al.(2006)Khuri, Mukherjee, Sinha, and Ghosh}]{Khuri2006}
Khuri, A.~I., Mukherjee, B., Sinha, B.~K., and Ghosh, M. (2006),
  \enquote{{Design issues for generalized linear models: a review},}
  \textit{Statistical Science}, 21, 376 -- 399.

\bibitem[{Li and Lin(2002)}]{Li2002}
Li, R. and Lin, D. K.~J. (2002), \enquote{Data analysis in supersaturated
  designs,} \textit{Statistics and Probablity Letters}, 59, 135--144.

\bibitem[{Li and Wu(1997)}]{LiWu97}
Li, W. and Wu, C. F.~J. (1997), \enquote{Columnwise-pairwise alogrithims with
  applications to the construction of supersaturated designs,}
  \textit{Technometrics}.

\bibitem[{Lin(1993)}]{Lin93}
Lin, D. K.~J. (1993), \enquote{A new class of supersaturated designs,}
  \textit{Technometrics}, 35, 28--31.

\bibitem[{Lin(1995)}]{Lin95}
--- (1995), \enquote{Generating systematic supersaturated designs,}
  \textit{Technometrics}, 37, 213--225.

\bibitem[{Liu et~al.(2007)Liu, Ruan, and Dean}]{liu2007construction}
Liu, Y., Ruan, S., and Dean, A.~M. (2007), \enquote{Construction and analysis
  of {$E(s^2)$} efficient supersaturated designs,} \textit{Journal of
  Statistical Planning and Inference}, 137, 1516--1529.

\bibitem[{Lounici(2008)}]{lounici2008sup}
Lounici, K. (2008), \enquote{Sup-norm convergence rate and sign concentration
  property of lasso and Dantzig estimators,} \textit{Electronic Journal of
  statistics}, 2, 90--102.

\bibitem[{Marley and Woods(2010)}]{MarleyWoods10}
Marley, C.~J. and Woods, D.~C. (2010), \enquote{{A comparison of design and
  model selection methods for supersaturated experiments},}
  \textit{Computational Statistics and Data Analysis}, 54, 3158--3167.

\bibitem[{Mee(2009)}]{mee}
Mee, R. (2009), \textit{A Comprehensive Guide to Factorial Two-Level
  Experimentation}, Springer.

\bibitem[{Mee et~al.(2017)Mee, Schoen, and Edwards}]{mee2017selecting}
Mee, R.~W., Schoen, E.~D., and Edwards, D.~J. (2017), \enquote{Selecting an
  orthogonal or nonorthogonal two-level design for screening,}
  \textit{Technometrics}, 59, 305--318.

\bibitem[{Meinshausen et~al.(2007)Meinshausen, Rocha, and
  Yu}]{meinshausen2007discussion}
Meinshausen, N., Rocha, G., and Yu, B. (2007), \enquote{Discussion: A tale of
  three cousins: lasso, L2Boosting and Dantzig,} \textit{The Annals of
  Statistics}, 35, 2373--2384.

\bibitem[{Meyer and Nachtsheim(1995)}]{meyer_nachtsheim1995}
Meyer, R.~K. and Nachtsheim, C.~J. (1995), \enquote{The coordinate-exchange
  algorithm for constructing exact optimal experimental designs,}
  \textit{Technometrics}, 37, 60--69.

\bibitem[{Nadarajah and Kotz(2008)}]{maxk2prob}
Nadarajah, S. and Kotz, S. (2008), \enquote{Exact Distribution of the Max/Min
  of Two Gaussian Random Variables,} \textit{IEEE Transactions on Very Large
  Scale Integration (VLSI) Systems}, 16, 210--212.

\bibitem[{Nguyen(1996)}]{Nguyen96}
Nguyen, N. (1996), \enquote{An alogrithmic approach to constructing
  supersaturated designs,} \textit{Technometrics}, 38, 69--73.

\bibitem[{Owen(1980)}]{Owen01011980}
Owen, D.~B. (1980), \enquote{A table of normal integrals,}
  \textit{Communications in Statistics - Simulation and Computation}, 9,
  389--419.

\bibitem[{Phoa et~al.(2009)Phoa, Pan, and Xu}]{phoa2009analysis}
Phoa, F.~K., Pan, Y.-H., and Xu, H. (2009), \enquote{Analysis of supersaturated
  designs via the Dantzig selector,} \textit{Journal of Statistical Planning
  and Inference}, 139, 2362--2372.

\bibitem[{Prékopa(1971)}]{LogConcaveIntegral}
Prékopa, A. (1971), \enquote{Logarithmic concave measures with application to
  stochastic programming,} \textit{Acta Scientiarum Mathematicarum}, 32,
  301--316.

\bibitem[{Pukelsheim(2006)}]{pukelsheim2006optimal}
Pukelsheim, F. (2006), \textit{Optimal design of experiments}, SIAM.

\bibitem[{Pukelsheim and Rieder(1992)}]{Pukelsheim1992}
Pukelsheim, F. and Rieder, S. (1992), \enquote{Efficient Rounding of
  Approximate Designs,} \textit{Biometrika}, 79, 763--770.

\bibitem[{Ryan and Bulutoglu(2007)}]{Ryan2007}
Ryan, K.~J. and Bulutoglu, D.~A. (2007), \enquote{$E(s^{2})$-optimal
  supersaturated designs with good minimax properties,} \textit{Journal of
  Statistical Planning and Inference}, 137, 2250--2262.

\bibitem[{Sarkar et~al.(2009)Sarkar, Lin, and Chatterjee}]{SARKAR20091224}
Sarkar, A., Lin, D.~K., and Chatterjee, K. (2009), \enquote{Probability of
  correct model identification in supersaturated design,} \textit{Statistics \&
  Probability Letters}, 79, 1224--1230.

\bibitem[{Silvey(1980)}]{Silvey1980}
Silvey, S.~D. (1980), \textit{Optimal Design}, Chapman and Hall.

\bibitem[{Singh and Stufken(2023)}]{singh2023selection}
Singh, R. and Stufken, J. (2023), \enquote{Selection of two-level
  supersaturated designs for main effects models,} \textit{Technometrics}, 65,
  96--104.

\bibitem[{Smucker and Drew(2015)}]{smucker2015approximate}
Smucker, B.~J. and Drew, N.~M. (2015), \enquote{Approximate model spaces for
  model-robust experiment design,} \textit{Technometrics}, 57, 54--63.

\bibitem[{Tang and Wu(1997)}]{tang1997method}
Tang, B. and Wu, C. (1997), \enquote{A method for constructing supersaturated
  designs and its {$E(s^2)$} optimality,} \textit{Canadian Journal of
  Statistics}, 25, 191--201.

\bibitem[{Tibshirani(2012)}]{tibshirani2012lasso}
Tibshirani, R.~J. (2012), \enquote{The lasso problem and uniqueness,}
  \textit{Electronic Journal of Statistics}, 7, 1456--1490.

\bibitem[{Wainwright(2009)}]{wainwright2009sharp}
Wainwright, M.~J. (2009), \enquote{Sharp thresholds for high-dimensional and
  noisy sparsity recovery using $ell \_ {1}$ -constrained quadratic programming
  (lasso),} \textit{IEEE transactions on information theory}, 55, 2183--2202.

\bibitem[{Weese et~al.(2017)Weese, Edwards, and Smucker}]{weese2017powerful}
Weese, M.~L., Edwards, D.~J., and Smucker, B.~J. (2017), \enquote{Powerful
  supersaturated designs when effect directions are known,} \textit{Journal of
  Quality Technology}, 49, 265--277.

\bibitem[{Weese et~al.(2015)Weese, Smucker, and Edwards}]{weese_etal2015}
Weese, M.~L., Smucker, B.~J., and Edwards, D.~J. (2015), \enquote{Searching for
  powerful supersaturated designs,} \textit{Journal of Quality Technology}, 47,
  66--84.

\bibitem[{Weese et~al.(2021)Weese, Stallrich, Smucker, and
  Edwards}]{weese2021strategies}
Weese, M.~L., Stallrich, J.~W., Smucker, B.~J., and Edwards, D.~J. (2021),
  \enquote{Strategies for supersaturated screening: group orthogonal and
  constrained Var(s) designs,} \textit{Technometrics}, 63, 443--455.

\bibitem[{Westfall et~al.(1998)Westfall, Young, and Lin}]{westfall1998forward}
Westfall, P.~H., Young, S.~S., and Lin, D.~K. (1998), \enquote{Forward
  selection error control in the analysis of supersaturated designs,}
  \textit{Statistica Sinica}, 101--117.

\bibitem[{Wu(1993)}]{Wu93}
Wu, C. F.~J. (1993), \enquote{Construction of supersaturated designs through
  partially aliased intercations,} \textit{Biometrika}, 80, 661--669.

\bibitem[{Xing(2015)}]{phdthesis}
Xing, D. (2015), \enquote{Lasso-optimal supersaturated design and analysis For
  factor screening in simulation experiments,} Ph.D. thesis, Purdue University.

\bibitem[{Xu et~al.(2009)Xu, Phoa, and Wong}]{xu2009recent}
Xu, H., Phoa, F.~K., and Wong, W.~K. (2009), \enquote{Recent developments in
  nonregular fractional factorial designs,} \textit{Statistics Surveys}, 3,
  18--46.

\bibitem[{Yang and Stufken(2009)}]{YangStufken2009}
Yang, M. and Stufken, J. (2009), \enquote{{Support points of locally optimal
  designs for nonlinear models with two parameters},} \textit{The Annals of
  Statistics}, 37, 518 -- 541.

\bibitem[{Zhang and Huang(2008)}]{zhang2008sparsity}
Zhang, C.-H. and Huang, J. (2008), \enquote{The sparsity and bias of the lasso
  selection in high-dimensional linear regression,} \textit{The Annals of
  Statistics}, 36, 1567--1594.

\bibitem[{Zhao and Yu(2006)}]{zhao2006model}
Zhao, P. and Yu, B. (2006), \enquote{On model selection consistency of lasso,}
  \textit{The Journal of Machine Learning Research}, 7, 2541--2563.

\end{thebibliography}

\newpage

\setcounter{equation}{0}

\footnotesize

\section{Proofs and Examples}

The equations have been renumbered here starting with (1). When referring to previous equation numbers, we will say the equation refers to one from the main document.

\subsection{Proof of Lemma 1}
First, the random vectors $\uvec$ and $\vvec$ in the $S_\lambda$ and $I_\lambda$ events, respectively are multivariate normal random vectors with covariance $0$. Thus, $\uvec$ and $\vvec$ are independent, meaning the $S_\lambda$ and $I_\lambda$ events themselves are independent. Therefore, 
\begin{align}
    \phi_\lambda(\X \, | \, \betavec)=P(S_\lambda \cap I_\lambda \, | \, \C, \, \V_{\mathcal{A}}, \, \betavec) = P(S_\lambda \, | \, \C_\mathcal{A}, \, \V_\mathcal{A}, \, \betavec_{\mathcal{A}} ) \times P(I_\lambda \, | \, \C, \, \zvec)\ .\ 
\end{align}
\noindent
Showing that $\phi_\lambda(\X \, | \, \betavec)=\phi_\lambda(\X \, | \, -\betavec)$ is equivalent to showing that both $P(S_\lambda \, | \, \F_\mathcal{A}, \, \betavec_{\mathcal{A}})= P(S_\lambda \, | \, \F_\mathcal{A}, \, -\betavec_{\mathcal{A}})$ and $P(I_\lambda \, | \, \F , \zvec )=P(I_\lambda \, | \, \F , -\zvec )$.

The normal random vector $\uvec$ in the $S_\lambda$ event and the $\vvec$ in the $I_\lambda$ event can we written as functions of the standard normal vector $\evec$ in the following ways:

\begin{align}
\uvec &=  -\frac{1}{\sqrt{n}}\Z_{\mathcal{A}}\C_\mathcal{A}^{-1} \F_{\mathcal{A}}^{T}\evec + \sqrt{n}\lambda \Z_\mathcal{A}\C_\mathcal{A}^{-1}\zvec_\mathcal{A}\ ,\ \label{eq:Sevent_rewrite}\\  
\boldsymbol{v} &= \frac{1}{\sqrt{n}}\F_{ \mathcal{I}}^{T}(\I-\Pmat_{\mathcal{A}})\evec +  \sqrt{n}\lambda \C_{\mathcal{I}\mathcal{A}}\C_{\mathcal{A}}^{-1}\zvec_{\mathcal{A}} \label{eq:Ievent_rewrite}\ .\
\end{align}

\noindent
where $\Pmat_{\mathcal{A}}=\F_\mathcal{A}(\F_\mathcal{A}^{T}\F_\mathcal{A})^{-1}\F_\mathcal{A}^{T}$. Both $\evec$ and $-\evec$ follow a $N(\zerovec,\I)$ so
 \begin{align*}
     P(I_\lambda \, | \, \C , \zvec )&=P\left( \ \left|\frac{1}{\sqrt{n}}\F_{ \mathcal{I}}^{T}(\I-\Pmat_{\mathcal{A}})\evec +  \sqrt{n}\lambda \C_{\mathcal{I}\mathcal{A}}\C_{\mathcal{A}}^{-1}\zvec_{\mathcal{A}} \right| \leq \sqrt{n}\lambda \onevec \right)\\
     &=P\left( \ \left|\frac{1}{\sqrt{n}}\F_{ \mathcal{I}}^{T}(\I-\Pmat_{\mathcal{A}})(-\evec) +  \sqrt{n}\lambda \C_{\mathcal{I}\mathcal{A}}\C_{\mathcal{A}}^{-1}(-\zvec_{\mathcal{A}}) \right| \leq \sqrt{n}\lambda \onevec \right)\\
     &=P(I_\lambda \, | \, \C , -\zvec )
     \end{align*}
     \noindent
    From a similar argument and the identities $\Z_\mathcal{A}\V_\mathcal{A}^{1/2}\betavec_\mathcal{A}=-\Z_\mathcal{A}\V_\mathcal{A}^{1/2}(-\betavec_\mathcal{A})$ and $\Z_{\mathcal{A}}\C_\mathcal{A}^{-1}\zvec_{\mathcal{A}}^T=(-\Z_{\mathcal{A}})\C_\mathcal{A}^{-1}(-\zvec_{\mathcal{A}})^T$, we see that $P(S_\lambda \, | \, \C_\mathcal{A}, \, \V_\mathcal{A}, \, \betavec_{\mathcal{A}})= P(S_\lambda \, | \, \C_\mathcal{A}, \, \V_\mathcal{A}, \, -\betavec_{\mathcal{A}})$. Thus, $\phi_\lambda(\X \, | \, \betavec)=\phi_\lambda(\X \, | \, -\betavec)$.   
    
    \subsection{Proof of Theorem 1}

Without loss of generality, let $\X= (\X_{\mathcal{A}}| \X_{\mathcal{I}})$ where $\X_{\mathcal{A}}$ represents the columns of active effects and $\X_{\mathcal{I}}$ represents the columns of inactive effects. Furthermore, let $\Fmat$ be the centered and scaled version of $\X$ and write $\widetilde{\Z} = \text{Diag}(
\widetilde{\Z}_{\mathcal{A}}, \widetilde{\Z}_{\mathcal{I}})$. Then $\X\widetilde{\Z} = 
(\X_{\mathcal{A}}\widetilde{\Z}_{\mathcal{A}} \, | \, \X_{\mathcal{I}}\widetilde{\Z}_{\mathcal{I}})$.
From Lemma 1, if this theorem holds for $\widetilde{\Z}_{\mathcal{A}}$, then it will hold for $-\widetilde{\Z}_{\mathcal{A}}$. Thus, this only needs to be proven for $\widetilde{\Z}_{\mathcal{A}}$. Also, note that, since post multiplying by $\widetilde{\Z}$ only changes the sign of the columns of $\X$, the centered and scaled version of $\X\widetilde{\Z}$ is $\Fmat\widetilde{\Z}$. Therefore, Theorem 1 is proven if 
\begin{equation*}
    \begin{split}
       P\left( I_\lambda \, | \, \C, \, \onevec\right) &= P\left(I_\lambda \, | \, \widetilde{\Z}\C\widetilde{\Z}, \,  \tilde{\zvec}\right)\\
        P\left( S_\lambda \, | \, \C_{\mathcal{A}}, \,\V_{\mathcal{A}}, \, \betavec_{\mathcal{A}}\right) &= P\left( S_\lambda \, | \, \widetilde{\Z}_{\mathcal{A}}\C_{\mathcal{A}}\widetilde{\Z}_{\mathcal{A}}, \, \V_{\mathcal{A}}, \,\tilde{\betavec}_{\mathcal{A}}\right)
    \end{split}
\end{equation*}
Denoting
\begin{equation*}
    \begin{split}
        \C&= \begin{bmatrix}
        \C_{\mathcal{A}} & \C_{\mathcal{A}\mathcal{I}}\\
        \C_{\mathcal{I}\mathcal{A}}& \C_{\mathcal{I}}\\
        \end{bmatrix}\\
        \widetilde{\Z}\C\widetilde{\Z}&= \begin{bmatrix}
        \widetilde{\Z}_{\mathcal{A}}\C_{\mathcal{A}}\widetilde{\Z}_{\mathcal{A}} & \widetilde{\Z}_{\mathcal{A}}\C_{\mathcal{A}\mathcal{I}}\widetilde{\Z}_{\mathcal{I}}\\
        \widetilde{\Z}_{\mathcal{I}}\C_{\mathcal{I}\mathcal{A}}\widetilde{\Z}_{\mathcal{A}}& \widetilde{\Z}_{\mathcal{I}}\C_{\mathcal{I}}\widetilde{\Z}_{\mathcal{I}}\\
        \end{bmatrix}\ ,\
    \end{split}
\end{equation*}
the probability of the $I_{\lambda}$ event may be written as
\begin{equation*}
    \begin{split}
        P\left( I_\lambda \, | \, \C, \, \onevec\right) &=P\left\{ \ \left|\frac{1}{\sqrt{n}}\F_{ \mathcal{I}}^{T}(\I-\Pmat_{\mathcal{A}})\evec +  \sqrt{n}\lambda \C_{\mathcal{I}\mathcal{A}}\C_{\mathcal{A}}^{-1}\onevec_{\mathcal{A}} \right| \leq \sqrt{n}\lambda \onevec \ \right\} \,\
    \end{split}
\end{equation*}
where $\Pmat_{\mathcal{A}}=\F_\mathcal{A}(\F_\mathcal{A}^{T}\F_\mathcal{A})^{-1}\F_\mathcal{A}^{T}$. Consider the event
\begin{equation*}
    \begin{split}
        P\left( I_\lambda \, | \, \widetilde{\Z}\C\widetilde{\Z}, \, \tilde{\zvec} \right)&=P\left\{ \ \left|\frac{1}{\sqrt{n}}\widetilde{\Z}_{\mathcal{I}}\F_{\mathcal{I}}^{*T}(\I-\Pmat_{\mathcal{A}})\evec +  \sqrt{n}\lambda \widetilde{\Z}_{\mathcal{I}}\C_{\mathcal{I}\mathcal{A}}\widetilde{\Z}_{\mathcal{A}}\widetilde{\Z}_{\mathcal{A}}\C_{\mathcal{A}}^{-1}\widetilde{\Z}_{\mathcal{A}}\tilde{\zvec}_{\mathcal{A}} \right| \leq \sqrt{n}\lambda \onevec \ \right\} \\
        &=P\left\{ \ \left|\widetilde{\Z}_{\mathcal{I}}\left[\frac{1}{\sqrt{n}}\F_{\mathcal{I}}^{T}(\I-\Pmat_{\mathcal{A}})\evec +  \sqrt{n}\lambda \C_{\mathcal{I}\mathcal{A}}\C_{\mathcal{A}}^{*-1}\onevec_{\mathcal{A}}\right] \right| \leq \sqrt{n}\lambda \onevec \ \right\}\ .\
    \end{split}
\end{equation*}

\noindent
The expression inside the absolute value in $I_\lambda \, | \, \widetilde{\Z}\C\widetilde{\Z}, \, \tilde{\zvec}$ is simply the product of $\widetilde{\Z}_{\mathcal{I}}$ and the expression inside the absolute value in $I_\lambda \, | \, \C, \, \onevec$. Since $\widetilde{\Z}_{\mathcal{I}}$ is diagonal with $\pm 1$ on the diagonal, the absolute value of the expressions are equal, meaning that we can conclude that $P(I_\lambda \, | \, \C, \, \onevec) = P(I_\lambda \, | \, \widetilde{\Z}\C\widetilde{\Z}, \, \tilde{\zvec})$. 

For the $S_{\lambda}$ event, note it depends on $\C_{\mathcal{A}}$ and 
\begin{equation*}
    P\left( S_\lambda \, | \, \C_{\mathcal{A}}, \, \V_{\mathcal{A}}, \, \betavec_{\mathcal{A}}\right)
    =P\left\{  -\frac{1}{\sqrt{n}}\C_{\mathcal{A}}^{-1} \F_{\mathcal{A}}^{T}\evec < \sqrt{n} \V^{1/2}\betavec_{\mathcal{A}} - \sqrt{n}\lambda\C_{\mathcal{A}}^{-1}\onevec_\mathcal{A}\right\} \ .\ 
\end{equation*}
Consider the event
\begin{equation*}
\begin{split}
       P( S_\lambda \, | \, &\widetilde{\Z}_{\mathcal{A}}\C_{\mathcal{A}}\widetilde{\Z}_{\mathcal{A}}, \, \V_{\mathcal{A}}, \,\tilde{\betavec}_{\mathcal{A}} ) =\\
       &P\left\{  -\frac{1}{\sqrt{n}}\widetilde{\Z}_{\mathcal{A}}\left[\widetilde{\Z}_{\mathcal{A}}\C_{\mathcal{A}}^{ -1}\widetilde{\Z}_{\mathcal{A}}\right]\widetilde{\Z}_{\mathcal{A}} \F_{\mathcal{A}}^{T}\evec < \sqrt{n} \widetilde{\Z}_{\mathcal{A}}\V^{1/2} \tilde{\betavec}_{\mathcal{A}} - \sqrt{n}\lambda\widetilde{\Z}_{\mathcal{A}}\widetilde{\Z}_{\mathcal{A}}\C_{\mathcal{A}}^{ -1}\widetilde{\Z}_{\mathcal{A}}\tilde{\zvec}_{\mathcal{A}} 
       \right\} \\
       = & P\left\{  -\frac{1}{\sqrt{n}}\C_{\mathcal{A}}^{ -1} \F_{\mathcal{A}}^{T}\evec < \sqrt{n} \widetilde{\Z}_{\mathcal{A}}\widetilde{\Z}_{\mathcal{A}}\V^{1/2} \betavec_{\mathcal{A}} - \sqrt{n}\lambda \C_{\mathcal{A}}^{ -1}\onevec_{\mathcal{A}}
       \right\} \\
       = & P\left\{  -\frac{1}{\sqrt{n}}\C_{\mathcal{A}}^{ -1} \F_{\mathcal{A}}^{T}\evec < \sqrt{n} \V^{1/2} \betavec_{\mathcal{A}} - \sqrt{n}\lambda \C_{\mathcal{A}}^{ -1}\onevec_{\mathcal{A}}
       \right\} \\
\end{split}
\end{equation*}
Thus, we have $ P(S_\lambda \, | \, \C_{\mathcal{A}}, \,\V_{\mathcal{A}}, \, \betavec_{\mathcal{A}}) = P(S_\lambda \, | \, \widetilde{\Z}_{\mathcal{A}}\C_{\mathcal{A}}\widetilde{\Z}_{\mathcal{A}}, \, \V_{\mathcal{A}}, \,\tilde{\betavec}_{\mathcal{A}})$, and therefore we have $\phi_\lambda(\X \, | \, \betavec)=\phi_\lambda(\X\tilde{\Z} \, | \, \tilde{\betavec})$. This completes the proof.



\subsection{Proof of Corollary 1}

This follows directly from Lemma 1. Since $\phi_\lambda(\X \, | \, \betavec)=\phi_\lambda(\X \, | \, -\betavec)$, the $\phi_\lambda ^{\pm}$-criterion does not need to consider both $\tilde{\zvec}$ and $-\tilde{\zvec}$. Thus, only the sign vectors in $\mathcal{Z}_\mathcal{A}^\pm$ need to be considered in $\phi_\lambda ^{\pm}$. 

\subsection{Proof of Proposition 1}
For any design $\X_1=(\X_{\mathcal{A}} \, | \, \X_{\mathcal{I}})$ consider the alternative design $\X_2=(\X_{\mathcal{A}} \, | \J)$, with obvious definitions for their model matrices $\F_1=(\F_\mathcal{A} | \F_\mathcal{I})$ and $\F_2=(\F_\mathcal{A} | \zerovec)$, respectively. Let $\C_1$ and $\C_2$ represent the correlation matrices for $\X_1$ and $\X_2$, respectively. Note that, by construction, $\C_1$ and $\C_2$ have the same $\C_{\mathcal{A}}$ submatrix. Clearly $P(I_\lambda \, | \, \C_2 , \zvec )=1$ as the left hand side of the event's inequality equals $\zerovec$.  Then
\[
\phi_\lambda(\X_1 \, | \, \betavec)= P(S_\lambda \, | \, \C_\mathcal{A}, \,\V_{\mathcal{A}}, \, \betavec_{\mathcal{A}}) \times P(I_\lambda \, | \, \C_1 , \zvec ) < P(S_\lambda \, | \,\C_\mathcal{A}, \, \V_{\mathcal{A}}, \, \betavec_\mathcal{A}) = \phi_\lambda(\X_2 \, | \, \betavec)\ .\
\]
This completes the proof.

\subsection{Derivation of $\phi_\lambda(\C \, |\, \betavec)=P(S_\lambda \, | \, \C_\mathcal{A}, \, \I, \, \betavec_{\mathcal{A}}=\beta \onevec)$ for $k=2$}

This derivation corresponds to Section~3.1.1 in the main document which assumes $\V=\I$ and absorbs $\sqrt{n}$ into $\betavec$ and $\lambda$. Rewrite $P(S_\lambda \, | \, \C_\mathcal{A}, \, \I, \, \betavec_{\mathcal{A}}=\beta \onevec)$ from equation (7) of the main document as
$P(\uvec^* < \zerovec)$ where $\uvec^*$  has the same covariance as $\uvec$, being $\C_\mathcal{A}^{-1}$ but with mean $\lambda \C_\mathcal{A}^{-1}\onevec-\beta\onevec$. Then $P(\uvec^* < \zerovec)=P(\max(U_1^*,U_2^*) < 0)$. Following \cite{maxk2prob} and some algebra, the density function of $U^*=\max(U_1^*,U_2^*)$ is
\begin{align}
    f_{U^*}(u)=2\sqrt{1-c^2}\times g\left(\sqrt{\frac{1-c}{1+c}}[(1+c)(\beta+u)-\lambda]\right)\times G((1+c)(\beta+u)-\lambda)\ ,\
\end{align}
and so $P(S_\lambda \, | \, \C_\mathcal{A}, \, \I, \, \betavec_{\mathcal{A}}=\beta \onevec)=\int_{-\infty}^0 f_{U^*}(u) \ du$. With the change of variables $t=\sqrt{\frac{1-c}{1+c}}[(1+c)(\beta+u)-\lambda]$, $dt=\sqrt{1-c^2}du$ and the bounds of integration with respect to $t$ is $-\infty < t < \mu(c)=\sqrt{\frac{1-c}{1+c}}[(1+c)\beta-\lambda]$, we have the expression
\begin{align}
    \int_{-\infty}^0 f_{U^*}(u) \ du=2 \int_{-\infty}^{\mu(c)} g(t)G\left(t\sqrt{\frac{1+c}{1-c}} \right) dt\ .\ \label{eqn:SmaxU}
\end{align}
where $g(\cdot)$ and $G(\cdot)$ refer to the probability density function and cumulative distribution function, respectively, for a standard normal distribution. Using \cite{Owen01011980}, you can show that \eqref{eqn:SmaxU} can be written as
\begin{align}
    [2G(\mu(c)-1]\times G((1+c)\beta-\lambda)+2T\left((1+c)\beta-\lambda, \sqrt{\frac{1-c}{1+c}}\right)\ .\
\end{align}

\subsection{Example violating bivariate log concavity of integrand of $P(S_\lambda \, | \, \C_\mathcal{A}, \, \I, \, \betavec_{\mathcal{A}}=\beta \onevec)$ for $k=2$}

One technique to establish log concavity of $P(S_\lambda \, | \, \C_\mathcal{A}, \, \I, \, \betavec_{\mathcal{A}}=\beta \onevec)$ using the Prékopa–Leindler inequality \citep{LogConcaveIntegral} requires establishing the integrand on the right hand side of \eqref{eqn:SmaxU} is bivariate log concave. This fails when $t=0.25$, as shown in the following figure.

\begin{figure}[h]
    \centering
    \includegraphics[width=0.5\linewidth]{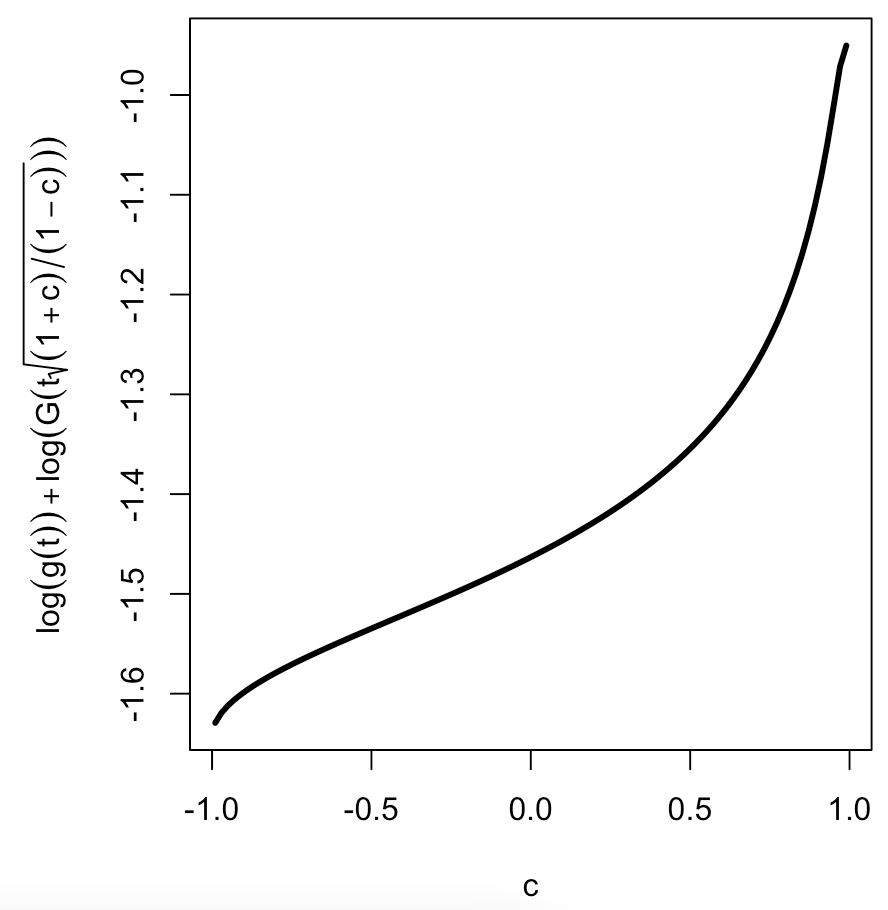}
    \caption{Plot of the log of $g(t)G\left(t\sqrt{\frac{1+c}{1-c}} \right)$ as a function of $c$ with $t=0.25$. This shows the function is not bivariate log concave.}
    \label{fig:Counterexample}
\end{figure}

\subsection{Symmetry of $\log[\phi_\lambda^\pm(\C_\mathcal{A} \, | \, \betavec_\mathcal{A})]$ for known $\betavec$ for $k=2$}

The only two sign vectors to consider for $\phi_\lambda^\pm(\C_\mathcal{A} \, | \, \betavec_\mathcal{A})$ are $\zvec_1^T=(1,1)^T$ and $\zvec_2^T=(1, -1)$. Therefore 
\begin{align}
    \phi_\lambda^\pm(\C_\mathcal{A} \, | \, \betavec)=\frac{1}{2}(\phi_\lambda(\C_\mathcal{A} | \betavec_\mathcal{A}=\beta\onevec_2)+\phi_\lambda(\C_\mathcal{A} | \betavec_\mathcal{A}=\beta\zvec_2))\ .\ \label{eqn:Sk2plusminus}
\end{align}
Following Theorem~1, $\phi_\lambda(\C_\mathcal{A} | \betavec_\mathcal{A}=\beta\zvec_2)=\phi_\lambda(\Z_2\C_\mathcal{A}\Z_2 | \betavec_\mathcal{A}=\beta\onevec_2)$ where $\Z_2=\text{Diag}(\zvec_2)$. Note that $\Z_2\C_\mathcal{A}\Z_2$ replaces the off-diagonal of $\C_\mathcal{A}$ with $-c$. Hence the first and second summands in \eqref{eqn:Sk2plusminus} consider the $2 \times 2$ correlation matrices with off-diagonals $\pm c$. Hence the function is symmetric about $c=0$.

\subsection{Proof of Lemma 2}

The distribution details of random normal vectors $\uvec$ and $\vvec$ are analogues of the random normal vectors in equations (7) and (8) of the main document, respectively, when $\C = (1-c)\I + c\J$ for $-(k-1)^{-1}<c<1$. Additionally, since $\V = \I$, the right-hand side of equation (7) from the main document becomes $\Z_\mathcal{A}\betavec_\mathcal{A} = |\betavec_\mathcal{A}|$. (Reminder: the approximate lasso criteria in this section absorb $\sqrt{n}$ into the values of $\betavec_\mathcal{A}$ and $\lambda$) Also, $\C_\mathcal{A}$, $\C_\mathcal{I}$, $\C_{\mathcal{A}\mathcal{I}}$, and $\C_{\mathcal{I}\mathcal{A}}$ are the same for any $\mathcal{A}$. Thus, if $\betavec_\mathcal{A}$ does not depend on $\mathcal{A}$ (e.g., $\betavec_\mathcal{A} = \beta\onevec$), then the probability of the two events are constant for any $\mathcal{A}$.

\subsection{Proof of Lemma 3}
We will prove this lemma in two parts. First we will show that $\frac{d}{dc}P(I_\lambda | c, \zvec)\big|_{c=0} = 0$ for any $\zvec \in \mathcal{Z}$, which includes $\zvec_\mathcal{A} =\onevec$. Then we derive conditions where $\frac{d}{dc} P(S_\lambda | c, \betavec_\mathcal{A}=\beta\onevec)\big|_{c=0} > 0$. 

\textbf{Part 1: Proving $\frac{d}{dc} P(I_\lambda | c, \zvec)\big|_{c=0} = 0$ for all $\lambda$ and all $\zvec$}. For an arbitrary sign vector for active effects, $\zvec_A$, let $\boldsymbol{v}\sim N\left(\lambda z_\mathcal{A}\gamma\onevec,(1-c)\left[\I + \gamma\J\right] \right)$ with pdf $f(\vvec,c,\zvec_{\mathcal{A}})$, where $z_\mathcal{A}= \onevec^T\zvec_\mathcal{A}$ and $\gamma = c/(1+c(k-1))$. Recall that 
$$
P(I_\lambda \, | \, c, \zvec) = P(|\boldsymbol{v}| \leq \lambda \onevec) = P(-\lambda \onevec\leq \boldsymbol{v} \leq \lambda \onevec) \ .\
$$

\noindent
We will use Lebnitz's rule iteratively over the $q$ dimensional integral. Note that, since $c$ is not in the bounds of the probability, using Lebnitz's rule is quite simple:

\begin{equation}
    \frac{d}{dc}P(I_\lambda \, | \, c, \zvec)\rvert_{c=0} = \int_{-\lambda \onevec}^{\lambda \onevec}\frac{d}{dc}\left[f(\vvec,c,\zvec_\mathcal{A})\right]d\vvec\big |_{c=0} \ .\
\end{equation}

\noindent
 When $c=0$, $\boldsymbol{v}\sim N(0,\I_q)$ which does not depend on $\zvec_\mathcal{A}$, so $f(\boldsymbol{v},c=0, \zvec_\mathcal{A})= \prod_{i=1}^q g(v_i)$ where $g(v_i)$ is the standard normal PDF. For simplicity, denote $f(\boldsymbol{v},c=0, \zvec_\mathcal{A})$ with $f_ v(\vvec,0)$. After taking the derivative and simplifying, we see that 
 \begin{equation}
     \begin{split}
         \frac{d}{dc}P(I_\lambda \, | \, c, \zvec)\rvert_{c=0}& =  -\frac{1}{2}\int_{-\lambda \onevec}^{\lambda \onevec}\boldsymbol{v}^T(\I_q-\J_q)\boldsymbol{v}f_v(\boldsymbol{v},0)d\boldsymbol{v}\\
         & \qquad +\lambda z_\mathcal{A} \int_{-\lambda \onevec}^{\lambda \onevec} \boldsymbol{v}^T\onevec f_v(\boldsymbol{v},0)d\boldsymbol{v}\ .\
     \end{split}
 \end{equation}
 
Intuitively, $\lambda z_{\mathcal{A}} \int_{-\lambda \onevec}^{\lambda \onevec} \boldsymbol{v}^T\onevec f(\boldsymbol{v},0)d\boldsymbol{v}$ is the symmetric integral of an odd function centered at 0, so it evaluates to zero. Additionally, note $\boldsymbol{v}^T(\I_q-\J_q)\boldsymbol{v}=-2 \sum_{i<j}v_iv_j$, so
 \begin{equation}
     \begin{split}
         \int_{-\lambda \onevec}^{\lambda \onevec}\boldsymbol{v}^T(\I_q-\J_q)\boldsymbol{v}f_v(\boldsymbol{v},0)d\boldsymbol{v} & = \int_{-\lambda \onevec}^{\lambda \onevec}-2 \sum_{i<j}v_iv_j f_v(\boldsymbol{v},0)d\boldsymbol{v}\\
         &= -2\sum_{i<j}\int_{-\lambda \onevec}^{\lambda \onevec}v_iv_j f_v(\boldsymbol{v},0)d\boldsymbol{v}\\
         &\varpropto -2\sum_{i<j}\int_{-\lambda}^{\lambda}\int_{-\lambda}^{\lambda} v_i v_j \exp\left(-\frac{1}{2}v_i^2-\frac{1}{2}v_j^2\right)dv_idv_j\\
         &= -2\sum_{i<j}\int_{-\lambda}^{\lambda} v_je^{-1/2v_j^2} \bigg(\int_{-\lambda}^{\lambda} v_i e^{-\frac{1}{2}v_i^2}dv_i\bigg)dv_j\\
         & = 0\ ,\
     \end{split}
 \end{equation}
 since $\int_{-\lambda}^{\lambda} v_i e^{-\frac{1}{2}v_i^2}dv_i=0$. So for an arbitrary sign vector $\zvec_A$, $\frac{d}{dc}P(I_\lambda \, | \, c, \zvec)\rvert_{c=0}=0$.

\textbf{Part 2: Proving $\frac{d}{dc} P(S_\lambda | c, \betavec_\mathcal{A}=\beta\onevec)\big|_{c=0} > 0$}. First, when $\zvec_\mathcal{A}=\onevec$ we can rewrite equation (7) in the main document by subtracting the mean:
\begin{equation}\label{eqn:Prob_S}
 P\left( \boldsymbol{u}^*<  \left[\beta - \frac{\lambda}{1+c(k-1)}\right]\onevec\right) \ ,\
\end{equation}
\noindent
where $\uvec^* \sim N\left( \zerovec, 1/(1-c)\left[ \I - \gamma \J\right]\right)$ with pdf $f(\uvec^*,c)$. The probability in~\eqref{eqn:Prob_S} can be written as the following integral
\begin{equation} \label{eqn:Prob_S_integ}
    \lim_{a\to -\infty}\bigg\{\int_{a}^{b(c)}\underset{k \text{ times}}{....}\int_{a}^{b(c)} f(\uvec^*,c)du_k... du_1\bigg\} \,\
\end{equation}
\noindent
where, $b(c) = \beta - \frac{\lambda}{1+c(k-1)}$ and $\uvec^{*\,T} = (u_1,u_2,\cdot \cdot \cdot u_{k-1},u_k)$. So then the derivative of the $S_\lambda$ event with respect to $c$ can be written as:

\begin{align}
    \frac{d}{dc}[P(S_\lambda \, | \, c, \beta)] & = \frac{d}{dc}\left[\lim_{a\to -\infty}\bigg\{\int_{a}^{b(c)}\underset{k \text{ times}}{....}\int_{a}^{b(c)} f(\uvec^*,c)du_k... du_1\bigg\}\right] \label{eq:PS_der_of_lim} \,\ \\
    & = \lim_{a\to -\infty}\frac{d}{dc}\bigg\{\int_{a}^{b(c)}\underset{k \text{ times}}{....}\int_{a}^{b(c)} f(\uvec^*,c)du_k... du_1\bigg\} \label{eq:PS_limit_of_der} \ .\ 
\end{align}

Utilizing an iterated application of the Liebnitz's Integral Rule and the fact that $f(\uvec^*, c=0)$ is the pdf of a multivariate independent standard normal vector , \eqref{eq:PS_limit_of_der} evaluated at $c=0$ becomes
\begin{equation}\label{eq:integral_simp}
    \begin{split}
         \lambda k(k-1)G(\tau)^{k-1}g(\tau)
         + \bigg\{\int_{-\infty}^{b(c)}\underset{k}{...}\int_{-\infty}^{b(c)}\bigg( \frac{d}{dc} f(\uvec^*,c)\bigg)du_k... du_2du_1\bigg \}\bigg\rvert_{c=0} \ ,\ \\
    \end{split}
\end{equation}
 where $g(\cdot)$ and $G(\cdot)$ represent the univariate standard normal pdf and CDF, respectively, with  $\tau = \beta-\lambda$. After much tedious calculus, we see that:
 
 \begin{equation}
\begin{split}
     \frac{d}{dc} f(\uvec^*,c)\rvert_{c=0}
    &= \frac{d}{dc} \left\{\frac{1}{(2\pi)^{k/2}|\C_{\mathcal{A}}|^{-1/2}} \exp\left(-\frac{1}{2} \uvec^{*\, T} ((1-c)\I + c\J)\uvec^* \right) 
    \right\} \bigg\rvert_{c=0}\\
    & = \frac{1}{2} \uvec^{* \, T}(\I-\J)\uvec^* f(\uvec^*,0)\\
    &= \frac{1}{2} \uvec^{* \, T}(\I-\J)\uvec^* \prod_{i=1}^k g(u_i)\\
    &= -\sum_{i<j} \left[\big(u_iu_j)\prod_{l=1}^k g(u_l)\right]\\
    & = - \sum_{i<j} \left(\prod_{l\neq i,j}^k g(u_l)\right) u_ig(u_i)u_jg(u_j) \ .\
    \end{split}
\end{equation}
 
Thus, the last term in Equation~\ref{eq:integral_simp} becomes:

\begin{equation}\label{eq:int_S_der}
    \begin{split}
        \bigg\{\int_{-\infty}^{b(c)}\underset{k}{...}&\int_{-\infty}^{b(c)}\bigg( \frac{d}{dc} f(\uvec,c)\bigg)du_k... du_2du_1\bigg \}\bigg\rvert_{c=0}\\
        & = - \sum_{i<j} \int_{-\infty}^{\tau}\underset{k}{...} \int_{-\infty}^{\tau} \left(\prod_{l\neq i,j}^k g(u_l)\right) u_ig(u_i)u_jg(u_j) 
        \ du_k... du_2du_1 \\
        & =- \sum_{i<j} G(\tau)^{k-2} \int_{-\infty}^{\tau}\int_{-\infty}^{\tau}u_ig(u_i)u_jg(u_j)du_idu_j\\
        & = -\sum_{i<j} G(\tau)^{k-2} \bigg( \frac{\exp\big\{ \frac{-\tau^2}{2}\big\}}{\sqrt{2\pi}} \bigg)^2\\
        & =-\frac{k(k-1)}{2}G(\tau)^{k-2} g\left(\tau\right)^2\ ,\
    \end{split}
\end{equation}
 where we have used the identity $\int xg(x)dx = -g(x) + D$, where $D$ is a constant. Lastly, by plugging \eqref{eq:int_S_der} into \eqref{eq:integral_simp}, it follows
 \begin{equation}
     \begin{split}
 \frac{d}{dc} P(S_\lambda \, | \, c, \beta) \rvert_{c=0} &=
 \lambda k(k-1)G(\tau)^{k-1}g(\tau) \\
 & \qquad - \frac{k(k-1)}{2}G(\tau)^{k-2} g\left(\tau\right)^2
     \end{split}
 \end{equation}
 \noindent
 Hence $\frac{d}{dc} P(S_\lambda \, | \, c, \beta) \rvert_{c=0} >0$ if and only if
 \begin{equation}
     2\lambda > \frac{g\left(\tau\right)}{G(\tau)}\ .\
 \end{equation} 
 \subsection{Proof of Theorem 2}
By the product rule, $\frac{d}{dc} \psi_\lambda(c |k, \beta)\big|_{c=0}$ equals
 \begin{equation*}
      \frac{d}{dc} P(S_\lambda \, | \, c,  \beta ) \big|_{c=0} P(I_\lambda \, | \, c=0, \zvec=\onevec) + \frac{d}{dc} P(I_\lambda \, | \, c, \zvec=\onevec) \big|_{c=0}P(S_\lambda \, | \, c=0,  \beta) \ .\
 \end{equation*}
 Since $P(S_\lambda \, | \, c=0,  \beta)$ and $P(I_\lambda \, | \, c=0, \zvec=\onevec)$ are non-negative and from the results in Lemma 3, $\frac{d}{dc} \psi_\lambda(c |k, \beta)\big|_{c=0} > 0$ when the conditions for Lemma 3 hold. Thus, Theorem 2 is proven. $\square$
 \subsection{Proof of Theorem 3}
This proof will contain two parts. First we show that $\frac{d}{dc} \psi_\lambda^{\pm}(c |k, \beta)\big|_{c=0}=0$. Then we develop the condition for which  $\psi_\lambda^{\pm}(c \, | \, k, \beta)$ is locally concave at $c=0$. When the concavity condition holds, $c=0$ is a local maximum for $\psi_\lambda^{\pm}(c \, | \, k, \beta)$.

\textbf{Part 1: Proving $\frac{d}{dc} \psi_\lambda^{\pm}(c |k, \beta)\big|_{c=0}=0$}. Note that if the individual entries of  $\tilde{\zvec}$, denoted $\tilde{z}_i$, are i.i.d. with probability $0.5$ to $\tilde{z}_i=1$ and $\tilde{z}_i=-1$, we can rewrite $\psi_\lambda^{\pm}(c \, | \, k, \beta)$ as:

\begin{align}
\psi_\lambda^{\pm}(c \, | \, k, \beta)=\mathbb{E}_{\tilde{\zvec}} \left\{P(S_\lambda \, | \, c, \betavec_\mathcal{A}=\beta\tilde{\zvec}) \times P(I_\lambda \, | \, c, \zvec=\tilde{\zvec})\right\} \label{eqn:CS_EZ}\ .\
\end{align}
\noindent
So, in taking the derivative, it is clear that:
\begin{align}
\frac{d}{dc} \psi_\lambda^{\pm}(c |k, \beta)\big|_{c=0}=\mathbb{E}_{\tilde{\zvec}} \left\{\frac{d}{dc}\left[P(S_\lambda \, | \, c, \betavec_\mathcal{A}=\beta\tilde{\zvec}) \times P(I_\lambda \, | \, c, \zvec=\tilde{\zvec})\right]\big |_{c=0}\right\} \label{eqn:CS_EZ_der}\ .\
\end{align}

We can rewrite and simplify the derivative on  the right-hand side of (\ref{eqn:CS_EZ_der}) using the product rule and the fact that $\frac{d}{dc} P(I_\lambda \, | \, c, \zvec=\tilde{\zvec}) \big |_{c=0} =0$ for any $\tilde{\zvec}$, we see that: 
\begin{equation}\label{eqn:CS_product_rule_gensign}
    \begin{split}
     &\frac{d}{dc}\left[P(S_\lambda \, | \, c, \betavec_\mathcal{A}=\beta\tilde{\zvec}) \times P(I_\lambda \, | \, c, \zvec=\tilde{\zvec})\right]\big |_{c=0}\\
     & =\frac{d}{dc}\left[P(S_\lambda \, | \, c, \betavec_\mathcal{A}=\beta\tilde{\zvec})\right ]\big |_{c=0} P(I_\lambda \, | \, 0, \zvec=\tilde{\zvec}) + 0\times P(S_\lambda \, | \, 0, \betavec_\mathcal{A}=\beta\tilde{\zvec}) \\
     & = \frac{d}{dc}\left[P(S_\lambda \, | \, c, \betavec_\mathcal{A}=\beta\tilde{\zvec}\right) ]\big |_{c=0} P(I_\lambda \, | \, 0, \zvec=\tilde{\zvec})
    \end{split}
\end{equation}
\noindent
Also note that for $c=0$, the $I_\lambda$ event does not depend on $\tilde{\zvec}$, so $P(I_\lambda \, | \, 0, \zvec=\tilde{\zvec})$ is a constant denoted simply $P(I_\lambda)$. By plugging \eqref{eqn:CS_product_rule_gensign} into \eqref{eqn:CS_EZ_der}, and using the fact that, when $c=0$, the $I_\lambda$ event does not depend on $\tilde{\zvec}$:

\begin{equation}
    \begin{split}
        \frac{d}{dc} \psi_\lambda^{\pm}(c |k, \beta)\big|_{c=0} &= \mathbb{E}_{\tilde{\zvec}} \left\{\frac{d}{dc}\left[P(S_\lambda \, | \, c, \betavec_\mathcal{A}=\beta\tilde{\zvec}\right) ]\big |_{c=0} P(I_\lambda \, | \, 0, \zvec=\tilde{\zvec})\right\}\\
        &= P(I_\lambda)\mathbb{E}_{\tilde{\zvec}} \left\{\frac{d}{dc}\left[P(S_\lambda \, | \, c, \betavec_\mathcal{A}=\beta\tilde{\zvec}\right) ]\big |_{c=0}\right\} \ .\ \\
    \end{split}
\end{equation}
Thus, $\frac{d}{dc} \psi_\lambda^{\pm}(c |k, \beta)\big|_{c=0} =0$ if and only if $\mathbb{E}_{\tilde{\zvec}} \left\{\frac{d}{dc}\left[P(S_\lambda \, | \, c, \betavec_\mathcal{A}=\beta\tilde{\zvec}\right) ]\big |_{c=0}\right\}=0$

When $\zvec_\mathcal{A}=\tilde{\zvec}$ we can rewrite $P(S_\lambda \, | \, c, \betavec_\mathcal{A})$ by subtracting the mean as:

\begin{equation}\label{eqn:Prob_S_gen_sign}
 P\left( \boldsymbol{u}< \beta\onevec - \frac{\lambda}{1-c}\left[\onevec - \tilde{z}^*\gamma\tilde{\zvec} \right]\right) \ ,\
\end{equation}
where $\uvec \sim N\left( \zerovec,\frac{1}{1-c} \left[\I_{k} - \gamma\tilde{\zvec}\tilde{\zvec}^T\right]\right)$ and $\tilde{z}^*=\onevec^T\tilde{\zvec}$. Let $f(\uvec,c, \tilde{\zvec})$ be the corresponding pdf. Thus, \eqref{eqn:Prob_S_gen_sign} can be expressed as

\begin{equation} \label{eqn:Prob_S_integ_gen_sign}
    \lim_{a\to -\infty}\bigg\{\int_{a}^{b_1(c)}\underset{k \text{ times}}{....}\int_{a}^{b_k(c)} f(\uvec,c, \tilde{\zvec})du_k... du_1\bigg\} \ ,\
\end{equation}
\noindent
where $b_i(c) = \beta - \frac{\lambda}{1-c}\left[1-\tilde{z}^*\gamma\tilde{z}_i \right]$ and $\uvec^T = (u_1,u_2,\cdot \cdot \cdot u_{k-1},u_k)$. So then the derivative of the $S_\lambda$ event with respect to $c$ can be written as:

\begin{align}
    \frac{d}{dc}[P(S_\lambda \, | \, c, \betavec_\mathcal{A}=\beta\tilde{\zvec})] & = \frac{d}{dc}\left[\lim_{a\to -\infty}\bigg\{\int_{a}^{b_1(c)}\underset{k \text{ times}}{....}\int_{a}^{b_k(c)} f(\uvec,c, \tilde{\zvec})du_k... du_1\bigg\}\right] \label{eq:PS_der_of_lim_all_signs} \,\ \\
    & = \lim_{a\to -\infty}\frac{d}{dc}\bigg\{\int_{a}^{b_1(c)}\underset{k \text{ times}}{....}\int_{a}^{b_k(c)} f(\uvec,c, \tilde{\zvec})du_k... du_1\bigg\} \label{eq:PS_limit_of_der_allsigns} \ .\ 
\end{align}

Let $\uvec_{-i}^T= (u_1, u_2, \cdot 
\cdot \cdot u_{i-1}, u_{i+1}, \cdot \cdot \cdot u_{k-1},u_k)$ , and let $\uvec_{-i}(b)^T= (u_1, u_2, \cdot 
\cdot \cdot u_{i-1},b_i(c), u_{i+1}, \cdot \cdot \cdot u_{k-1},u_k)$. Lastly let $\boldsymbol{b}_{-i}(c)^T = (b_1(c), b_2(c), \cdot 
\cdot \cdot b_{i-1}(c), b_{i+1}(c), \cdot \cdot \cdot b_{k-1}(c),b_k(c)) $. Using Liebnitz's rule once on $$\frac{d}{dc}\bigg\{\int_{a}^{b(c)}\underset{k \text{ times}}{....}\int_{a}^{b(c)} f(\uvec,c, \tilde{\zvec})du_k... du_1\bigg\}\ ,\ $$ we see that $\frac{d}{dc}\bigg\{\int_{a}^{b_1(c)}{....}\int_{a}^{b_k(c)} f(\uvec,c, \tilde{\zvec})du_k... du_1\bigg\} \bigg | _{c=0}$
has the following expression:
\begin{equation}\label{eqn:leb_one_step}
    \begin{split}
        &\frac{d}{dc}\left(b_1(c)\right)\int_{a}^{b_2(c)}{....}\int_{a}^{b_k(c)} f(\uvec_{-1}(b),c, \tilde{\zvec})du_k... du_2\bigg | _{c=0}\\
        &\qquad + \int_{a}^{b_1(c)}\frac{d}{dc}\bigg\{\int_{a}^{b_2(c)}{....}\int_{a}^{b_k(c)} f(\uvec,c, \tilde{\zvec})du_k... du_2\bigg\}\bigg|_{c=0}du_1\\
        & = \frac{d}{dc}\left(b_1(c)\right)\big|_{c=0}\int_{a}^{b_2(0)}{....}\int_{a}^{b_k(0)} f(\uvec_{-1}(b),c=0, \tilde{\zvec})du_k... du_2\\
        &\qquad + \int_{a}^{b_1(c)}\frac{d}{dc}\bigg\{\int_{a}^{b_2(c)}{....}\int_{a}^{b_k(c)} f(\uvec,c, \tilde{\zvec})du_k... du_2\bigg\}\bigg|_{c=0}du_1\ .\ \\
    \end{split}
\end{equation}

Note that, $\frac{d}{dc}\left(b_1(c)\right)\big|_{c=0} = \lambda\left[\tilde{z}_i\tilde{z}^*-1\right]$. By applying Lebnitz's rule again to the derivative in the last line of (24) and  iterating, we see that \eqref{eq:PS_limit_of_der_allsigns} becomes
\begin{equation}\label{eqn:limit_leb}
    \begin{split}
        \lim_{a\to -\infty}\bigg\{&\sum_{i=1}^{k}\left[\lambda\left(\tilde{z}_i\tilde{z}^*-1\right) \int_{-\infty}^{\boldsymbol{b}_i(0)}f(\uvec_{-i}(b), c=0, \tilde{z})d\uvec_{-i}\right]\\
        &+ \int_{a}^{b_1(c)}{....}\int_{a}^{b_k(c)} \frac{d}{dc}\left[f(\uvec,c, \tilde{\zvec})\right]\big |_{c=0}du_k... du_2du_1\bigg\} \ .\ \\
    \end{split}
\end{equation}
\noindent
When $c=0$, $f(\uvec,c=0, \tilde{\zvec})$ is the pdf of $N(\zerovec,\I_k)$ and does not depend on $\tilde{\zvec}$ (and similar for $f(\uvec_{-i},c=0, \tilde{\zvec})$). For brevity, let $f(\uvec,0) = f(\uvec,c=0, \tilde{\zvec})$. Additionally, letting $g(\cdot)$ represent the pdf of the univariate standard normal distribution, it is clear that $f(\uvec_{-i}(b), 0) = g(b(0))f(\uvec_{-i},0)$. Thus, we can simplify by the first term in the summand of (25) by applying the limit as $a \to -\infty$, giving the expression
\begin{equation}
    \begin{split}
        \frac{d}{dc}[P(S_\lambda \, | \, c, \betavec_\mathcal{A}=\beta\tilde{\zvec})] & = g(\tau)G\left(\tau\right)^{k-1}\sum_{i=1}^k \left[\lambda\left(\tilde{z}_i\tilde{z}^*-1\right)\right]\\
        & \qquad + \int_{a}^{\tau}\underset {k \text{ times}}{....}\int_{a}^{\tau} \frac{d}{dc}\left[f(\uvec,c, \tilde{\zvec})\right]\big | _{c=0}du_k... du_2du_1 \ .\
    \end{split}
\end{equation}
\noindent
After much tedious calculus, we see that:
 
 \begin{equation}
\begin{split}
     \frac{d}{dc} f(&(\uvec,c, \tilde{z})\rvert_{c=0}\\
    &= \frac{d}{dc} \left\{\frac{1}{(2\pi)^{k/2}\left|\frac{1}{1-c}\left(I-\gamma\tilde{\zvec}\tilde{\zvec}^T\right)\right|^{-1/2}} \exp\left[-\frac{1}{2} \uvec^T ((1-c)\I + c\tilde{\zvec}\tilde{\zvec}^T)\uvec\right] \right\}\bigg \rvert_{c=0}\\
    & = \frac{1}{2} \uvec^T(\I-\tilde{\zvec}\tilde{\zvec}^T)\uvec f(\uvec,0)\\
    &= -\sum_{i<j}u_i\tilde{z}_i\tilde{z}_ju_jf(\uvec,0) \ .\
\end{split}
\end{equation}

\noindent
Taking the expectation over equally likely sign vectors, we have
\begin{equation}\label{eqn:Expectation_simple}
    \begin{split}
        \mathbb{E}_{\tilde{\zvec}} \bigg\{ \frac{d}{dc}
        P(S_\lambda \, &| \, c, \betavec_\mathcal{A}=\beta\tilde{\zvec}) |_{c=0} \bigg\}\\ &= \mathbb{E}_{\tilde{\zvec}} \bigg\{ g(\tau)G\left(\tau\right)^{k-1}\sum_{i=1}^k \left[\lambda\left(\tilde{z}_i\tilde{z}^*-1\right)\right]\\
        & \qquad - \int_{a}^{\tau}\underset {k \text{ times}}{....}\int_{a}^{\tau} \sum_{i<j}u_i\tilde{z}_i\tilde{z}_ju_jf(\uvec,0) du_k... du_2du_1\bigg\} \ .\
    \end{split}
\end{equation}
\noindent
By linearity of expectation, \eqref{eqn:Expectation_simple} becomes:
\begin{equation}\label{eqn:expectation_brokendown}
    \begin{split}
        &\mathbb{E}_{\tilde{\zvec}} \bigg\{ g(\tau)G\left(\tau\right)^{k-1}\sum_{i=1}^k \left[\lambda\left(\tilde{z}_i\tilde{z}^*-1\right)\right]\bigg \}\\
        & \qquad -\mathbb{E}_{\tilde{\zvec}} \bigg\{\int_{a}^{\tau}\underset {k \text{ times}}{....}\int_{a}^{\tau} \sum_{i<j}u_i\tilde{z}_i\tilde{z}_ju_jf(\uvec,0) du_k... du_2du_1\bigg\}\\
        & = g(\tau)G\left(\tau\right)^{k-1}\sum_{i=1}^k \mathbb{E}_{\tilde{\zvec}}\left \{\left[\lambda\left(\tilde{z}_i\tilde{z}^*-1\right)\right]\right\}\\
        & \qquad - \int_{a}^{\tau}\underset {k \text{ times}}{....}\int_{a}^{\tau}f(\uvec,0) \left(\sum_{i<j}u_i\mathbb{E}_{\tilde{\zvec}}[\tilde{z}_i\tilde{z}_j]u_j\right ) du_k... du_2du_1 \ .\
    \end{split}
\end{equation}

Note the following expectation properties of $\tilde{z}_i$ and $\tilde{z}^*=\sum_i \tilde{z}_i$: 
\begin{equation}\label{eqn:Expectation_Facts}
    \begin{split}
    \mathbb{E}_{\tilde{\zvec}}\tilde{z}_i &= 0\\
    \mathbb{E}_{\tilde{\zvec}}(\tilde{z}_i\tilde{z}_j) &= \mathbb{E}_{\tilde{\zvec}}(\tilde{z}_i)\mathbb{E}_{\tilde{\zvec}}(\tilde{z}_j) = 0\\
    \mathbb{E}_{\tilde{\zvec}}(\tilde{z}_i^2) &= 1\\
    \mathbb{E}_{\tilde{\zvec}}(\tilde{z}_i\tilde{z}^*) &=\mathbb{E}_{\tilde{\zvec}}(\tilde{z}_i^2) + \sum_{i \neq j}\mathbb{E}_{\tilde{\zvec}}(\tilde{z}_i\tilde{z}_j) = 1 \ .\
    \end{split}
\end{equation}
Then $\mathbb{E}_{\tilde{\zvec}}\left \{\left[\lambda\left(\tilde{z}_i\tilde{z}^*-1\right)\right]\right\} = 0$ and $\sum_{i<j}u_i\mathbb{E}_{\tilde{\zvec}}[\tilde{z}_i\tilde{z}_j]u_j =0$, proving that $\mathbb{E}_{\tilde{\zvec}} \left\{\frac{d}{dc}\left[P(S_\lambda \, | \, c, \betavec_\mathcal{A}=\beta\tilde{\zvec}\right) ]\big |_{c=0}\right\}=0$ which happens if and only if $\frac{d}{dc} \psi_\lambda^{\pm}(c |k, \beta)\big|_{c=0} =0$.

\textbf{Part 2: Concavity condition for $\psi_\lambda^{\pm}(c |k, \beta)$ at $c=0$}. We aim to show what conditions must hold for $\psi_\lambda^{\pm}(c \, | \, k, \beta)$ to be locally concave at $c=0$. Clearly, 
\begin{align}
\frac{d^2}{dc^2} \psi_\lambda^{\pm}(c |k, \beta)\big|_{c=0}=\mathbb{E}_{\tilde{\zvec}} \left\{\frac{d^2}{dc^2}\left[P(S_\lambda \, | \, c, \betavec_\mathcal{A}=\beta\tilde{\zvec}) \times P(I_\lambda \, | \, c, \zvec=\tilde{\zvec})\right]\big |_{c=0}\right\} \ .\
\end{align}
and the second derivative of the product can be expressed as:
\begin{equation}\label{eq:2nd_der_phi}
    \begin{split}
        \frac{d^2}{dc^2} [ P(S_\lambda \, | \, c, \betavec_\mathcal{A}=\beta\tilde{\zvec}) &\times P(I_\lambda \, | \, c, \zvec=\tilde{\zvec}) ] \big |_{c=0}\\
        &= \frac{d^2}{dc^2}\left\{P(I_\lambda \, | \, c, \zvec=\tilde{\zvec}) \right \}\big |_{c=0}P(S_\lambda \, | \, c=0, \betavec_\mathcal{A}=\beta\tilde{\zvec})\\
        & + 2\frac{d}{dc}\left\{P(I_\lambda \, | \, c, \zvec=\tilde{\zvec}) \right \}\big |_{c=0}\frac{d}{dc}\left\{P(S_\lambda \, | \, c, \betavec_\mathcal{A}=\beta\tilde{\zvec}) \right \}\big |_{c=0}\\
        &+ \frac{d^2}{dc^2}\left\{P(S_\lambda \, | \, c=0, \betavec_\mathcal{A}=\beta\tilde{\zvec}) \right \}\big |_{c=0}P(I_\lambda \, | \, c=0, \zvec=\tilde{\zvec}) \ .\ \\
    \end{split}
\end{equation}
Since $\frac{d}{dc}\left\{P(I_\lambda \, | \, c, \zvec=\tilde{\zvec}) \right \}\big |_{c=0} = 0$, and both $I_\lambda$ and $S_\lambda$ do not depend on $\tilde{\zvec}$ when$c=0$, it follows from \eqref{eq:2nd_der_phi} that
\begin{equation}\label{eq:simplified_2nd_der_EZ}
    \begin{split}
    \frac{d^2}{dc^2} \psi_\lambda^{\pm}(c |k, \beta)\big|_{c=0}&=\mathbb{E}_{\tilde{\zvec}} \bigg\{\frac{d^2}{dc^2}\left\{P(I_\lambda \, | \, c, \zvec=\tilde{\zvec}) \right \}\big |_{c=0}P(S_\lambda \, | \, c=0, \betavec_\mathcal{A}=\beta\tilde{\zvec}) \\
    & \qquad+ \frac{d^2}{dc^2}\left\{P(S_\lambda \, | \, c=0, \betavec_\mathcal{A}=\beta\tilde{\zvec}) \right \}\big |_{c=0}P(I_\lambda \, | \, c=0, \zvec=\tilde{\zvec}) \bigg \} \\
    & = \mathbb{E}_{\tilde{\zvec}} \bigg\{\frac{d^2}{dc^2}\left\{P(I_\lambda \, | \, c, \zvec=\tilde{\zvec}) \right \}\big |_{c=0}\bigg \}P(S_\lambda \, | \, c=0) \\
    & \qquad+ \mathbb{E}_{\tilde{\zvec}} \bigg\{\frac{d^2}{dc^2}\left\{P(S_\lambda \, | \, c, \betavec_\mathcal{A}=\beta\tilde{\zvec}) \right \}\big |_{c=0}\bigg \}P(I_\lambda \, | \, c=0) \ .\ \\
    \end{split}
\end{equation}
If \eqref{eq:simplified_2nd_der_EZ} is less than or equal to 0, we have concavity (and thus a local maximum) at $c=0$. We will provide expressions for each of the second derivatives separately, and then combine the expressions to get a condition for concavity to hold.


First note that 
\begin{equation}\label{eq:second_der_expect}
    \begin{split}
        \mathbb{E}_{\tilde{\zvec}} \bigg\{\frac{d^2}{dc^2}\left\{P(S_\lambda \, | \, c, \betavec_\mathcal{A}=\beta\tilde{\zvec}) \right \}\big |_{c=0}\bigg \}& =  \mathbb{E}_{\tilde{\zvec}} \bigg\{\frac{d^2}{dc^2}\left\{\int_{-\infty}^{|\betavec|}f_u(\uvec,c,\tilde{\zvec}) d\uvec \right \}\bigg |_{c=0}\bigg \}\\
        & = \mathbb{E}_{\tilde{\zvec}} \bigg\{\int_{-\infty}^{|\betavec|}\frac{d^2}{dc^2}\left\{f_u(\uvec,c,\tilde{\zvec}) \right \}\bigg |_{c=0}d\uvec\bigg \}\\
        &= \int_{-\infty}^{|\betavec|}\mathbb{E}_{\tilde{\zvec}} \bigg\{\frac{d^2}{dc^2}\left\{f_u(\uvec,c,\tilde{\zvec}) \right \}\big |_{c=0}\bigg \}d\uvec\ ,\
    \end{split}
\end{equation}
 where $f_u(\uvec, c, \tilde{\zvec})$ is the pdf of $N\left(\frac{\lambda}{1-c}\left[\onevec - \tilde{z}^*\gamma\tilde{\zvec} \right],\frac{1}{1-c} \left[\I_{k} - \gamma\tilde{\zvec}\tilde{\zvec}^T\right] \right) $. Write 
     \begin{align}
         f_u(\uvec, c, \tilde{\zvec}) &= t(c)\exp\{h(\uvec, c, \tilde{\zvec}) \}\\
     t(c) &= \frac{1}{(2\pi)^{k/2}}\left| \frac{1}{1-c} \left[\I_{k} - \gamma\tilde{\zvec}\tilde{\zvec}^T\right]
     \right|^{-1/2}\\
     h(\uvec, c, \tilde{\zvec}) & = -\frac{1}{2}\left(\uvec -\frac{\lambda}{1-c}\left[\onevec - \tilde{z}^*\gamma\tilde{\zvec} \right]  \right)^T \left[(1-c)\I_{k} + c\tilde{\zvec}\tilde{\zvec}^T\right] \left(\uvec -\frac{\lambda}{1-c}\left[\onevec - \tilde{z}^*\gamma\tilde{\zvec} \right]  \right)\ .\
 \end{align}
Now let $t'(c) \equiv \frac{d}{dc}t(c)$ and $h'(\uvec,c, \tilde{\zvec}) \equiv \frac{d}{dc}h(\uvec,c,\tilde{\zvec})$ and likewise for the second derivative and double prime notation. Then, we can write
\begin{equation}
    \begin{split}
        \frac{d^2}{dc^2}f_u(\uvec, c, \tilde{\zvec}) & = t''(c) \exp \{h(\uvec, c, \tilde{\zvec})\} + 2t'(c)h'(\uvec,c,\tilde{\zvec}) \exp \{h(\uvec, c, \tilde{\zvec})\}\\
        & \qquad+ t(c) \exp \{h(\uvec, c, \tilde{\zvec})\}\left\{ h'(\uvec,c,\tilde{\zvec})^2 + h''(\uvec,c,\tilde{\zvec}) \right\}\\
    \end{split}
\end{equation}
It follows that $t'(0) = 0$, $t(0) = \frac{1}{(2\pi)^{k/2}}$, $t''(0) = -\frac{k(k-1)}{2(2\pi)^{k/2}}$, and 
\begin{equation}
    \begin{split}
        h'(\uvec, c, \tilde{\zvec}) & = -\frac{1}{2}\left[-\uvec^T\uvec + \uvec^T\tilde{\zvec}\tilde{\zvec}^T\uvec+   \frac{\lambda^2  k}{(1-c)^2} - \frac{(1+(k-1)c^2)\lambda^2  \tilde{z}^{*\,2} }{(1-c)^2 (1+c(k-1))^2}\right]\\
        h'(\uvec, c, \tilde{\zvec})|_{c=0} & = -\frac{1}{2}\left[-\uvec^T(I-\tilde{\zvec}\tilde{\zvec}^T) \uvec +   \lambda^2  k - \lambda^2 \tilde{z}^{*\,2} \right]\\
        h''(\uvec, c, \tilde{\zvec}) & = --\frac{1}{2}\left\{\frac{2\lambda^2  k}{(1-c)^3} - \frac{2\left[(k-1)^2 c^3 + 3(k-1) c -k +2 \right]\lambda^2  \tilde{z}^{*\,2} }{(1-c)^3 (1+c(k-1))^3}\right\}\\
        h''(\uvec, c, \tilde{\zvec})|_{c=0} & = -\lambda^2  \left[k + \tilde{z}^{*\,2}(k-2) \right]
    \end{split}
\end{equation}
This means that $t''(0)\exp \{h(\uvec, 0, \tilde{\zvec})\} = -\frac{k(k-1)}{2} f_u(\uvec, 0, \tilde{\zvec})$ and 
\begin{equation}
    \begin{split}
        \frac{d^2}{dc^2}f_u(\uvec, c&, \tilde{\zvec})\bigg |_{c=0}\\
        &= \frac{-(k-1)k}{2} f_u(\uvec, 0, \tilde{\zvec})\\
        & \qquad +f_u(\uvec, 0, \tilde{\zvec})\left\{ \frac{1}{4}\left[-\uvec^T(\I-\tilde{\zvec}\tilde{\zvec}^T)\uvec + \lambda^2(k-\tilde{z}^{*\,2})\right]^2 - \lambda^2  \left[k + \tilde{z}^{*\,2}(k-2)\right]\right\}\\
        & = f_u(\uvec, 0, \tilde{\zvec}) \bigg\{-\frac{k(k-1)}{2} +  \frac{1}{4}\left(\uvec^T(\I-\tilde{\zvec}\tilde{\zvec}^T)\uvec\right)^2 -\frac{1}{2}\lambda^2(k-\tilde{z}^{*\,2})\uvec^T(\I-\tilde{\zvec}\tilde{\zvec}^T)\uvec\\
        & \qquad +\frac{1}{4}\lambda^4(k-\tilde{z}^{*\,2})^2 - \lambda^2 \left[k + \tilde{z}^{*\,2}(k-2) \right] \bigg\} \ .\
    \end{split}
\end{equation}
Since $f_u(\uvec, 0, \tilde{\zvec})$ does not depend on the sign vector, we can rewrite:
\begin{equation}
    \begin{split}
        \mathbb{E}_{\tilde{\zvec}} \bigg\{\frac{d^2}{dc^2}\left\{f_u(\uvec,c,\tilde{\zvec}) \right \}\big |_{c=0}\bigg \} & = f_u(\uvec, 0, \tilde{\zvec}) \bigg\{\frac{-(k-1)k}{2} \\
        & \qquad + \frac{1}{4}\mathbb{E}_{\tilde{\zvec}}\left[\left(\uvec^T(\I-\tilde{\zvec}\tilde{\zvec}^T)\uvec\right)^2\right] \\
        &\qquad -\frac{1}{2}\lambda^2\mathbb{E}_{\tilde{\zvec}}\left[(k-\tilde{z}^{*\,2})\uvec^T(\I-\tilde{\zvec}\tilde{\zvec}^T)\uvec\right]\\
        &\qquad +\frac{1}{4}\lambda^4\mathbb{E}_{\tilde{\zvec}}\left[(k-\tilde{z}^{*\,2})^2\right] - \lambda^2  \left[k + \mathbb{E}_{\tilde{\zvec}}\left[\tilde{z}^{*\,2}\right](k-2) \right]\bigg\}
    \end{split}
\end{equation}
Note that $\uvec^T(\I-\tilde{\zvec}\tilde{\zvec}^T)\uvec = -2\sum_{i<j}u_i\tilde{z}_i\tilde{z}_ju_j$ and consider the following expectations
\begin{equation}\label{eqn:Expectation_Facts_second}
    \begin{split}
    \mathbb{E}_{\tilde{\zvec}}\tilde{z}_i &= 0\\
    \mathbb{E}_{\tilde{\zvec}}(\tilde{z}_i\tilde{z}_j) &= \mathbb{E}_{\tilde{\zvec}}(\tilde{z}_i)\mathbb{E}_{\tilde{\zvec}}(\tilde{z}_j) = 0\\
    \mathbb{E}_{\tilde{\zvec}}\tilde{z}_i^a\tilde{z}_j & =0 \text{ for all } a>0\\
    \mathbb{E}_{\tilde{\zvec}}(\tilde{z}_i^2) &= 1\\
    \mathbb{E}_{\tilde{\zvec}}(\tilde{z}_i\tilde{z}^*) &=\mathbb{E}_{\tilde{\zvec}}(\tilde{z}_i^2) + \sum_{i \neq j}\mathbb{E}_{\tilde{\zvec}}(\tilde{z}_i\tilde{z}_j) = 1\\
    \mathbb{E}_{\tilde{\zvec}}(\tilde{z}^{*\,2}) &=\mathbb{E}_{\tilde{\zvec}}(\sum_{i=1}^k\tilde{z}_i^2) + 2\sum_{i \neq j}\mathbb{E}_{\tilde{\zvec}}(\tilde{z}_i\tilde{z}_j) = k\\
    \mathbb{E}_{\tilde{\zvec}}(\tilde{z}^{*\,4}) &=k +3k(k-1)= 3k^2-2k\\
    \mathbb{E}_{\tilde{\zvec}}[(k-\tilde{z}^{*\,2})^2] &=k^2 - 2k\mathbb{E}_{\tilde{\zvec}}(\tilde{z}^{*\,2}) + \mathbb{E}_{\tilde{\zvec}}(\tilde{z}^{*\,4}) =2k(k-1)\\
    \mathbb{E}_{\tilde{\zvec}}[(k-\tilde{z}^{*\,2})\uvec^T(\I-\tilde{\zvec}\tilde{\zvec}^T)\uvec] &=4 \sum_{i<j}u_iu_j\\
    \mathbb{E}_{\tilde{\zvec}}[(\uvec^T(\I-\tilde{\zvec}\tilde{\zvec}^T)\uvec)^2] &=4 \sum_{i<j}u_i^2u_j^2\\
    \end{split}
\end{equation}
So substituting these in, we see that

\begin{equation}
\begin{split}
    \mathbb{E}_{\tilde{\zvec}} \bigg\{\frac{d^2}{dc^2}\left\{f_u(\uvec,c,\tilde{\zvec}) \right \}\big |_{c=0}\bigg \} & = f_u(\uvec, 0, \tilde{\zvec}) \left\{\frac{-(k-1)k}{2} - \lambda^2k(k-1) +\frac{\lambda^4  k(k-1)}{2}  \right\}\\
        &+ f_u(\uvec, 0, \tilde{\zvec})\bigg\{\sum_{i<j}u_i^2u_j^2 -2\lambda^2\sum_{i<j}u_iu_j \bigg\}\\
\end{split}
\end{equation}
\noindent
When $c=0$, $f_u(\uvec,c, \tilde{\zvec})= \prod_{i=1}^kf_{u_i}(u_i)$ where $f_{u_i}$ is the pdf for $N(\lambda, 1)$. So we can split the densities up after adding in the integration:

\begin{equation}
    \begin{split}
        \int_{-\infty}^{|\betavec|}\mathbb{E}_{\tilde{\zvec}} \bigg\{\frac{d^2}{dc^2}&[f_u(\uvec,c,\tilde{\zvec})]\big |_{c=0}\bigg \}d\uvec\\ & = \left\{\frac{-(k-1)k}{2} - \lambda^2nk(k-1) +\frac{\lambda^4 n^2 k(k-1)}{2} \right\}\int_{-\infty}^{|\betavec|}f_u(\uvec, 0, \tilde{\zvec})d\uvec \\
        & \qquad + \int_{-\infty}^{|\betavec|}f_u(\uvec, 0, \tilde{\zvec})\sum_{i<j}u_i^2u_j^2d\uvec\\
        &\qquad -2\lambda^2n \int_{-\infty}^{|\betavec|}f_u(\uvec, 0, \tilde{\zvec})\sum_{i<j}u_iu_jd\uvec \ .\
    \end{split}
\end{equation}
\noindent
Writing $\int_{-\infty}^{|\betavec|}f_u(\uvec, 0, \tilde{\zvec})d\uvec = G(\tau)^k$ where $G(\cdot)$ is the standard normal CDF gives
\begin{equation}
    \begin{split}
        \int_{-\infty}^{|\betavec|}\mathbb{E}_{\tilde{\zvec}} \bigg\{\frac{d^2}{dc^2} [f_u(\uvec,c&,\tilde{\zvec}) ]\big |_{c=0}\bigg \}d\uvec\\ & = \left\{\frac{-(k-1)k}{2} - \lambda^2k(k-1) +\frac{\lambda^4  k(k-1)}{2} \right\} G(\tau)^k\\
        &\qquad + G(\tau)^{k-2}\frac{k(k-1)}{2}\left\{\int_{-\infty}^{\beta}u_i^2f_{u_i}(u_i)du_i\right\}^2\\
        &\qquad -2\lambda^2 G(\tau)^{k-2}\frac{k(k-1)}{2} \left\{\int_{-\infty}^{\beta}u_if_{u_i}(u_i) du_i\right\}^2\ .\
    \end{split}
\end{equation}
We can express $\int_{-\infty}^{|\betavec|}\mathbb{E}_{\tilde{\zvec}} \bigg\{\frac{d^2}{dc^2}\left\{f_u(\uvec,c,\tilde{\zvec}) \right \}\big |_{c=0}\bigg \}d\uvec$ as: 
\begin{equation}
    \begin{split}
       G(\tau)^{k-2}\bigg \{&  \left(\frac{-(k-1)k}{2} - \lambda^2k(k-1) +\frac{\lambda^4 k(k-1)}{2} \right) G(\tau)^2 \\
       &+ \frac{k(k-1)}{2}\left\{\int_{-\infty}^{\beta}u_i^2f_{u_i}(u_i)du_i\right\}^2 -2\lambda^2n \frac{k(k-1)}{2} \left\{\int_{-\infty}^{\beta}u_if_{u_i}(u_i)du_i\right\}^2 \bigg\} \ .\ \\
    \end{split}
\end{equation}
Note that 
\begin{equation}
    \int_{-\infty}^{\beta}u_if_{u_i}(u_i)du_i = g(\tau) + \lambda G(\tau) \ ,\
\end{equation}
and (51) squared is
\begin{equation}
    \begin{split}
      \left\{\int_{-\infty}^{\beta}u_if_{u_i}(u_i)du_i\right\}^2  = g(\tau)^2-2\lambda g(\tau)G(\tau) + \lambda^2 G(\tau)^2\ .\
    \end{split}
\end{equation}
By substituting $v=u_i-\lambda$ and by the fact that, for a standard normal density $g(\cdot)$, that $\int v^2g(v)dv =G(v)-vg(v)+D$ for some constant $D$,
\begin{equation*}
    \begin{split}
       \int_{-\infty}^{\beta}u_i^2f_{u_i}(u_i)du_i & = \int_{-\infty}^{\tau}v^2g(v)dv
       + 2\lambda \int_{-\infty}^{\tau}vg(v)dv
       + \lambda^2 \int_{-\infty}^{\tau}g(v)dv\\
       &= (1+\lambda^2)G(\tau) -(\beta+\lambda)g(\tau)\ .\
    \end{split}
\end{equation*}
We then get the expression for $\mathbb{E}_{\tilde{\zvec}} \bigg\{\frac{d^2}{dc^2}\left\{P(S_\lambda \, | \, c, \betavec_\mathcal{A}\beta\tilde{\zvec}) \right \}\big |_{c=0}\bigg \}$: 
\begin{equation}\label{eqn:simplified_condition_2nd_der_PS}
\begin{split}
       k(k-1) g(\lambda) G(\tau)^{k-2} &\bigg\{-\left\{ (\beta +\lambda) + \lambda^2  (\beta-\lambda) \right\}G(\tau)
    + \left\{ \frac{\beta^2-\lambda^2}{2} + \beta\lambda\right\}g(\tau)\bigg\}\ .\ 
\end{split}
\end{equation}


Now we turn to the $I_\lambda$ event. Again, note that
\begin{equation}\label{eq:I_second_der_expect}
    \begin{split}
        \mathbb{E}_{\tilde{\zvec}} \bigg\{\frac{d^2}{dc^2}\left\{P(I_\lambda \, | \, c, \zvec=\tilde{\zvec}) \right \}\big |_{c=0}\bigg \}& =  \mathbb{E}_{\tilde{\zvec}} \bigg\{\frac{d^2}{dc^2}\left\{\int_{-\lambda \onevec}^{\lambda \onevec}f_v(\vvec,c,\tilde{\zvec}) d\vvec \right \}\bigg |_{c=0}\bigg \}\\
        & = \mathbb{E}_{\tilde{\zvec}} \bigg\{\int_{-\lambda \onevec}^{\lambda \onevec}\frac{d^2}{dc^2}\left\{f_v(\vvec,c,\tilde{\zvec}) \right \}\bigg |_{c=0}d\vvec\bigg \}\\
        & =  \int_{-\lambda \onevec}^{\lambda \onevec}\mathbb{E}_{\tilde{\zvec}} \left\{\frac{d^2}{dc^2}\left\{f_v(\vvec,c,\tilde{\zvec}) \right \}\bigg |_{c=0}\right\}d\vvec\\
    \end{split}
\end{equation}
 where $f_v(\vvec, c, \tilde{\zvec})$ is the pdf of $N\left(\lambda \tilde{z}\gamma \onevec,(1-c) \left[\I_{q} + \gamma\J_q\right] \right) $ for $\tilde{z}=\onevec^T\tilde{\zvec}$. Again, we write 
 \begin{align}
         f_v(\vvec, c, \tilde{\zvec}) &= t(c)\exp\{h(\vvec, c, \tilde{\zvec}) \} \ ,\ \\
     t(c) &= \frac{1}{(2\pi)^{q/2}}\left| (1-c) \left( \I_{q} + \gamma\J_q\right)
     \right|^{-1/2}\\
     h(\vvec, c, \tilde{\zvec}) & = -\frac{1}{2(1-c)}\left(\vvec -\lambda \tilde{z}\gamma \onevec  \right)^T \left[\I_{q} - \frac{c}{1+c(q+k-1)}\J_q\right] \left(\vvec -\lambda \tilde{z}\gamma \onevec  \right)\ .\
 \end{align}
Then $\frac{d^2}{dc^2}f_v(\vvec, c, \tilde{\zvec})$ has the same expression as (42) but with $t'(0) = 0$, $t(0) = \frac{1}{(2\pi)^{q/2}}$, $t''(0) = \frac{q(2q +4k -2)}{4}(2\pi)^{-q/2}$, and
\begin{equation}
    \begin{split}
        h'(\vvec, c, \tilde{\zvec})|_{c=0} & = - \frac{1}{2}\left[\vvec^T(\I_q-\J_q) \vvec - 2\lambda\tilde{z}\vvec^T\onevec\right]\\
        h''(\vvec, c, \tilde{\zvec})|_{c=0} & = -\frac{1}{2}\left[2\vvec^T\vvec + 2(q+k-2)\vvec \J_q\vvec -4\lambda \tilde{z}(q-k+2)\vvec^T\onevec +2\lambda^2\tilde{z}^2q \right]\ .\
    \end{split}
\end{equation}
As $t''(0)\exp \{h(\vvec, 0, \tilde{\zvec})\} = \frac{q(2q +4k -2)}{4} f_v(\vvec, 0, \tilde{\zvec})$, we get the expression
\begin{equation}
    \begin{split}
        \frac{d^2}{dc^2}&f_v(\vvec, c, \tilde{\zvec})\bigg |_{c=0}\\&= \frac{q(2q +4k -2)}{4} f_v(\vvec, 0, \tilde{\zvec})
        +f_v(\vvec, 0, \tilde{\zvec})\left\{ \frac{1}{4}\left[\vvec^T(\I_q-\J_q)\vvec - 2\lambda\tilde{z}\vvec^T\onevec\right]^2 \right\}\\
        & \quad +f_v(\vvec, 0, \tilde{\zvec})\left\{-\frac{1}{2}\left[2\vvec^T\vvec + 2(q+k-2)\vvec \J_q\vvec -4\lambda \tilde{z}(q-k+2)\vvec^T\onevec +2\lambda^2\tilde{z}^2q\right]\right\} \ .\ \\
    \end{split}
\end{equation}


\noindent
Since $f_v(\vvec, 0, \tilde{\zvec})$ does not depend on the sign vector, we can rewrite:
\begin{equation}
    \begin{split}
        \mathbb{E}_{\tilde{\zvec}} \bigg\{\frac{d^2}{dc^2}[f_v(\vvec,c&,\tilde{\zvec})  ]\big |_{c=0}\bigg \}\\  &= \frac{q(2q +4k -2)}{4} f_v(\vvec, 0, \tilde{\zvec})\\
        &\qquad +f_v(\vvec, 0, \tilde{\zvec})\bigg\{ \frac{1}{4}\big[(\vvec^T(\I_q-\J_q)\vvec)^2 -4 \lambda\mathbb{E}_{\tilde{\zvec}}\{\tilde{z}\}\vvec^T\onevec\vvec^T(\I_q-\J_q)\vvec\\
        &\qquad \qquad \qquad + 4\lambda^2 \mathbb{E}_{\tilde{\zvec}}\{\tilde{z}^2\} (\vvec^T\onevec)^2\big] \bigg\}\\
        &\qquad +f_v(\vvec, 0, \tilde{\zvec})\bigg\{-\frac{1}{2}\big[2\vvec^T\vvec + 2(q+k-2)\vvec \J_q\vvec\\
        &\qquad \qquad \qquad -4\lambda \mathbb{E}_{\tilde{\zvec}}\{\tilde{z}\}(q-k+2)\vvec^T\onevec +2\lambda^2\mathbb{E}_{\tilde{\zvec}}\{\tilde{z}^2\}q\big]\bigg\} \ .\ \\
    \end{split}
\end{equation}

By plugging in the expectation facts from equation (46), equation (60) becomes:
\begin{equation}
    \begin{split}
        \mathbb{E}_{\tilde{\zvec}} \bigg\{\frac{d^2}{dc^2}\left\{f_v(\vvec,c,\tilde{\zvec}) \right \}\big |_{c=0}&\bigg \}  = \frac{q(2q +4k -2)}{4} f_v(\vvec, 0, \tilde{\zvec})\\
        &\qquad +f_v(\vvec, 0, \tilde{\zvec})\left\{ \frac{1}{4}\left[(\vvec^T(\I_q-\J_q)\vvec)^2  + \lambda^2 k (\vvec^T\onevec)^2\right] \right\}\\
        &\qquad  +f_v(\vvec, 0, \tilde{\zvec})\bigg\{\frac{-1}{2}\big[2\vvec^T\vvec + 2(q+k-2)\vvec \J_q\vvec +2\lambda^2kq\big]\bigg\} \ .\ \\
    \end{split}
\end{equation}
We can write without loss of generality
\begin{equation}
    \begin{split}
       \int_{-\lambda \onevec}^{\lambda\onevec}v_if_v(\vvec, 0, \tilde{\zvec})d\vvec &= [1-2G(-\lambda)]^{q-1} \int_{-\lambda }^{\lambda}v_ig(v_i)dv_i\\
       &=0 \ ,\ \\
    \end{split}
\end{equation}
where $ G(\lambda)-G(-\lambda)=1-2G(-\lambda)$. 
Furthermore, we can generalize this result for $i \neq j$:
\begin{equation}
    \begin{split}
       \int_{-\lambda \onevec}^{\lambda\onevec}v_i^av_j f_v(\vvec, 0, \tilde{\zvec})d\vvec &= [1-2G(-\lambda)]^{q-2} \left(\int_{-\lambda }^{\lambda}v_i^a g(v_i) dv_i\right)\left(\int_{-\lambda }^{\lambda}v_j g(v_j)dv_j\right)\\
       &=0 \ ,\ \\
    \end{split}
\end{equation}
for any $a>0$. Integrating (61) gives the expression
\begin{equation}
    \begin{split}
        \frac{q(2q +4k -2)}{4} &\int_{-\lambda \onevec}^{\lambda\onevec}f_v(\vvec, 0, \tilde{\zvec})d\vvec\\
        &+\int_{-\lambda \onevec}^{\lambda\onevec}f_v(\vvec, 0, \tilde{\zvec})\left\{ \frac{1}{4}\left[(\vvec^T(\I_q-\J_q)\vvec)^2  + 4\lambda^2 k (\vvec^T\onevec)^2\right] \right\}d\vvec\\
        & +\int_{-\lambda \onevec}^{\lambda\onevec}f_v(\vvec, 0, \tilde{\zvec})\bigg\{-\frac{1}{2}\big[2\vvec^T\vvec + 2(q+k-2)\vvec \J_q\vvec +2\lambda^2kq\big]\bigg\} d\vvec \ .\ \\
    \end{split}
\end{equation}

First, $\int_{-\lambda \onevec}^{\lambda\onevec}f_v(\vvec, 0, \tilde{\zvec})d\vvec = [1-2G(-\lambda)]^q$. As $\vvec^T(\I_q-\J_q)\vvec = -2\sum_{i<j}v_iv_j$,  $[\vvec^T(\I_q-\J_q)\vvec]^2$ has $q(q-1)/2$ terms with $v_i^2v_j^2$ and other terms with $v_i$ in them. By (62), these latter terms will integrate to $0$, giving
\begin{equation}
    \begin{split}
        \frac{1}{4}\int_{-\lambda \onevec}^{\lambda\onevec} (\vvec^T(\I_q-\J_q)\vvec)^2f_v(\vvec, 0, \tilde{\zvec})d\vvec & = \int_{-\lambda \onevec}^{\lambda\onevec} \sum_{i<j}(v_i^2v_j^2)f_v(\vvec, 0, \tilde{\zvec})d\vvec\\
        & = \frac{q(q-1)}{2}[1-2G(-\lambda)]^{q-2}\left(\int_{-\lambda }^{\lambda}v_i^2g(v_i)dv_i\right )^2\ .\\\
    \end{split}
\end{equation}
Next, note that $(\vvec^T\onevec)^2 = \vvec^T\J_q\vvec = \sum_{i=1}^qv_i^2 + 2\sum_{i<j}v_iv_j$. From (63) the $v_iv_j$ terms will integrate to $0$. Thus, we can write:
\begin{equation}
    \begin{split}
        \frac{4\lambda^2k}{4}\int_{-\lambda \onevec}^{\lambda\onevec} (\vvec^T\onevec)^2f_v(\vvec, 0, \tilde{\zvec})d\vvec & = \int_{-\lambda \onevec}^{\lambda\onevec} \sum_{i=1}^q f_v(\vvec, 0, \tilde{\zvec})d\vvec\\
        & = q \lambda^2k [1-2G(-\lambda)]^{q-1}\int_{-\lambda }^{\lambda}v_i^2g(v_i)dv_i \ .\ \\
    \end{split}
\end{equation}
\noindent
Lastly, we can rewrite the last term in (64) as follows
\begin{equation}
    \begin{split}
        -\frac{1}{2}\int_{-\lambda \onevec}^{\lambda\onevec} \big[2\vvec^T\vvec &+ 2(q+k-2)\vvec \J_q\vvec +2\lambda^2kq\big]f_v(\vvec, 0, \tilde{\zvec})d\vvec\\
        & = - \int_{-\lambda \onevec}^{\lambda\onevec} \big[(1+q+k-2)\sum_{i=1}^q v_i^2 +2\lambda^2k q\big]f_v(\vvec, 0, \tilde{\zvec})d\vvec\\
        &=-q(q+k-1)[1-2G(-\lambda)]^{q-1}\int_{-\lambda }^{\lambda} v_i^2 g(v_i)dv_i
        -\lambda^2kq [1-2G(-\lambda)]^q \ .\ \\
    \end{split}
\end{equation}
With these simplifications, we arrive at the expression
\begin{equation}
    \begin{split}
        \int_{-\lambda \onevec}^{\lambda\onevec}\mathbb{E}_{\tilde{\zvec}} \bigg\{\frac{d^2}{dc^2}\left\{f_v(\vvec,c,\tilde{\zvec}) \right \}\big |_{c=0}\bigg \}d\vvec  &= \left(\frac{q(2q +4k -2)}{4}-q\lambda^2k \right)[1-2G(-\lambda)]^q\\
        &\quad +\frac{q(q-1)}{2}[1-2G(-\lambda)]^{q-2}\left(\int_{-\lambda }^{\lambda}v_i^2g(v_i)dv_i\right )^2\\
        &\quad +q \lambda^2k [1-2G(-\lambda)]^{q-1}\int_{-\lambda }^{\lambda}v_i^2g(v_i)dv_i\\
        &\quad -q(q+k-1)[1-2G(-\lambda)]^{q-1}\int_{-\lambda }^{\lambda} v_i^2 g(v_i)dv_i\\
    \end{split}
\end{equation}
Using the Gaussian integral identity 
\begin{equation}
    \int_{-\lambda }^{\lambda} v_i^2 g(v_i)dv_i= [1-2G(-\lambda)]-2\lambda g(\lambda) \ ,\
\end{equation}
we arrive at a final expression:
\begin{equation}
\begin{split}
\mathbb{E}_{\tilde{\zvec}} \bigg\{\frac{d^2}{dc^2} [P(I_\lambda \, | \, c&, \zvec=\tilde{\zvec}) ] \big |_{c=0}\bigg \}\\ &= 2q\lambda [1-2G(-\lambda)]^{q-2}g(\lambda)\bigg\{ \lambda (q-1)g(\lambda) + k(1-\lambda^2)[1-2G(-\lambda)]\bigg\} \ .\
\end{split}
\end{equation}

Recall that $\phi^{\pm}$ is locally concave at $c=0$ when equation \eqref{eq:simplified_2nd_der_EZ} is less than or equal to 0, or,
\begin{equation}
    \begin{split}
        \mathbb{E}_{\tilde{\zvec}} \bigg\{\frac{d^2}{dc^2} [P(I_\lambda \, | \, c, \zvec&=\tilde{\zvec}) ] \big |_{c=0}\bigg \}P(S_\lambda \, | \, c=0) \\ &\leq - \mathbb{E}_{\tilde{\zvec}} \bigg\{\frac{d^2}{dc^2}\left\{P(S_\lambda \, | \, c, \betavec_\mathcal{A}=\beta\tilde{\zvec}) \right \}\big |_{c=0}\bigg \}P(I_\lambda \, | \, c=0)\ .\
    \end{split}
\end{equation}
When $c=0$, $P(I_\lambda \, | \, c=0)=[1-2G(-\lambda)]^q$ and $P(S_\lambda \, | \, c=0)= G(\tau)^k$. Then, after tedious algebra, the inequality (71) can be written as
\begin{equation}
\begin{split}
   & 2\lambda q \,g(\lambda)G(\tau)^2\left[k(1-\lambda^2)[1-2G(-\lambda)] + (q-1)\lambda g(\lambda) \right]\leq \\
    &\qquad -k(k-1)g(\tau)[1-2G(-\lambda)]^2 \left\{ -\left[ (\beta +\lambda) +n\lambda^2(\beta-\lambda)\right]G(\tau) + \left[ \frac{\beta^2-\lambda^2}{2} + \beta\lambda\right]g(\tau)\right\} \ ,\ 
\end{split}
\end{equation}
or the slightly more condensed inequality
\begin{equation}
\begin{split}
\frac{q}{\binom{k}{2}} \ \frac{\lambda g(\lambda)}{[1-2G(-\lambda)]}&\left( k(1-\lambda^2)+\frac{(q-1)\lambda g(\lambda)}{[1-2G(-\lambda)]}\right) \leq\\
&\frac{g(\tau)}{G(\tau)}\left(\beta +\lambda+\lambda^2\tau-\left[\frac{\beta^2-\lambda^2}{2}+\beta\lambda\right]\frac{g(\tau)}{G(\tau)}\right) \ ,
\ 
\end{split}
\end{equation}
This completes the proof.

\section{Evaluating $\phi_{\max}$ and $\phi_\Lambda$}
 
This section describes the evaluation of $\phi_{\max}$ and $\phi_\Lambda$ in the software. These implementations are generally software independent, meaning it can be done in a wide variety of languages and utilizing different software packages. 

For $\phi_{\max}$, any optimization package can be used to find the $\lambda$ that maximizes $\phi_\lambda(\X, \betavec)$. However, it should be noted that, for $\Phi_{\max}^{\pm}(\X,k, \betavec)$, maximizing $\Phi_\lambda^{\pm}(\X,k, \betavec)$ over $\lambda$ could be computationally costly due to the number of evaluations of $\Phi_\lambda^{\pm}(\X,k, \betavec)$. While computational burden can be reduced by implementing some of the strategies given in Section 4, the implementation of $\Phi_{\max}^{\pm}(\X,k, \betavec)$ also begins with a warm-start to select a high-quality starting point for the optimization algorithm. The maximization is implemented over a user specified range of $\log(\lambda)$, and the warm- start evaluates $\Phi_{\max}^{\pm}(\X,k, \betavec)$ over a grid of evenly spaced $\log(\lambda)$ values. The $\log(\lambda)$ from this grid that maximizes $\Phi_{\max}^{\pm}(\X,k, \betavec)$ is then used as the starting value in the optimization algorithm. Our implementation uses the \texttt{optim} function in \texttt{R}.

For $\phi_{\Lambda}$, an adaptive Riemann sum integration is utilized over a range of $\log(\lambda)$. The adaptive Riemann sum begins with a user-specified lower bound for $\log(\lambda)$ (default is $-5$). The Riemann sum is then calculated over evenly spaced $\log(\lambda)$ values increasing from the lower bound. For each new $\log(\lambda)$ value, if $\phi_\lambda$ is non-zero or above some user-specified $\epsilon$, it is added to the Reimann sum and the sum continues with the next $\log(\lambda)$. If $\phi_\lambda < \epsilon$,  the Riemann sum is stopped and $\phi_\Lambda$ is given by the final value of the sum.

\end{document}